\documentclass[envcountsame,envcountresetsect]{svmult}

\usepackage{makeidx}         % allows index generation
\usepackage{graphicx}        % standard LaTeX graphics tool
\usepackage{multicol}        % used for the two-column index
\usepackage[bottom]{footmisc}% places footnotes at page bottom
\makeindex             % used for the subject index
\usepackage{mathrsfs,amssymb,amsmath,xspace}
\let\mathcal\mathscr
\usepackage[all]{xy}

\def\thetheorem{\thesection.\arabic{theorem}}
\def\thesubsection{\thesection.\Alph{subsection}}

\makeatletter
\DeclareRobustCommand*\textsubscript[1]{%
  \@textsubscript{\selectfont#1}}
\def\@textsubscript#1{%
  {\m@th\ensuremath{_{\mbox{\fontsize\sf@size\z@#1}}}}}

\let\c@equation\c@theorem
\def\qedsymbol{\ensuremath{\square}}
\smartqed
\def\endproof{\qed\@endtheorem }
\makeatother

\spdefaulttheorem{remarks}{Remarks}{\itshape}{\rmfamily}

\def\hdeg{\mathop{\widehat{\mathrm {deg}}}}

\def\ddc{\mathop{\mathrm d\mathrm d^{c}}}
\def\hZ^#1{\widehat{\mathrm{Z}}\vphantom{\operatorname{Z}}^{#1}}
\def\CH^#1{\operatorname{CH}}
\def\hCH^#1{\widehat{\mathrm{CH}}\vphantom{\operatorname{CH}}^{#1}}
\def\hdiv{\widehat{\mathrm{div}}}

\def\can{{\text{\upshape can}}}
\def\CAP{{\text{\upshape cap}}}
\def\Id{\operatorname{Id}}

\def\Gal{\operatorname{Gal}}
\let\mathbb\mathbf

\def\pr{\operatorname{pr}}

\def\norm#1{\mathopen\| {#1} \mathclose\|}
\def\abs#1{\mathopen| {#1} \mathclose|}

\def\div{\operatorname{div}}
\let\ra\rightarrow
\let\hra\hookrightarrow

\let\bar\overline

\let\geq\geqslant
\let\leq\leqslant
\let\epsilon\varepsilon
\let\phi\varphi
\let\hat\widehat

\def\A{\mathbf{A}}
\def\N{\mathbf{N}}
\def\P{\mathbf{P}}
\def\Q{\mathbf{Q}}
\def\R{\mathbf{R}}

\def\C{\mathbf{C}}
\def\F{\mathbf{F}}
\def\Z{\mathbf{Z}}
\def\Sym{\operatorname{Sym}}
\def\rank{\operatorname{rank}}

\def\an{{\operatorname{an}}}

\def\ord{{\operatorname{ord}}}
\def\Pic{\operatorname{Pic}}
\def\Spec{\operatorname{Spec}}
\def\Aut{\operatorname{Aut}}
\def\GL{\operatorname{GL}}
\def\Div{\operatorname{Div}}

\def\tube#1{\mathopen]{#1}\mathclose[}
\def\sozat{\,;\,}
\let\eps\varepsilon
\def\jet{\mathrm j}
\def\cf{\emph{cf.}\nobreak\space\ignorespaces}

\def\ie{\emph{i.e.}\nobreak\xspace}
\def\eg{\emph{e.g.}\nobreak\xspace}

\newcommand{\cC}{{\mathcal  C}}

\newcommand{\cE}{{\mathcal  E}}
\newcommand{\cF}{{\mathcal  F}}

\newcommand{\cL}{{\mathcal L}}
\newcommand{\cO}{{\mathcal O}}
\newcommand{\cP}{{\mathcal P}}

\def\cX{{\mathcal X}}
\def\cY{{\mathcal Y}}

\newcommand{\OK}{{\mathcal O}_K}

\def\labelenumi{\textup{\theenumi)}}

\begin{document}
\title* {Analytic curves in algebraic varieties over number fields}

\author{Jean-Beno\^{\i}t Bost\inst{1} \and Antoine Chambert-Loir\inst{2}}
\institute{D\'epartement de math\'ematiques,
   Universit\'e de Paris-Sud (Orsay),
   B\^atiment~425,
  F-91405 Orsay Cedex, France. \texttt{Jean-Benoit.Bost@math.u-psud.fr}
\and
Universit\'e de Rennes~1, Institut de recherche
math\'ematique de Rennes,
Campus de Beaulieu, F-35042 Rennes Cedex, France.
\texttt{antoine.chambert-loir@univ-rennes1.fr}}

\maketitle

\begin{center}\itshape
To Yuri Ivanovich Manin
\end{center}

\begin{abstract}
We establish algebraicity criteria for formal germs of curves in
algebraic varieties over number fields and apply them to derive a
rationality criterion for formal germs of functions on algebraic
curves, which extends the classical rationality theorems of Borel-Dwork
and P\'olya-Bertrandias valid over the projective line to arbitrary
algebraic curves over a number field. The formulation and the proof
of these criteria involve some basic notions in Arakelov geometry,
combined with complex and rigid analytic geometry (notably, potential
theory over complex and $p$-adic curves). We also discuss geometric
analogues, pertaining to the algebraic geometry of projective
surfaces, of these arithmetic criteria.
\end{abstract}


\par\addvspace\baselineskip
\noindent\subclassname\enspace
{14G40 (Primary); 14G22, 31A15, 14B20 (Secondary)}

\setcounter{tocdepth}{2}

\def\theequation{\thesection.\arabic{equation}}

%% \begin{center}
%% \emph{To Yuri Ivanovich Manin}
%% \end{center}
%% 
\section{Introduction}
The purpose of this article is to  establish \emph{algebraicity 
criteria} 
for \emph{formal germs of  curves} in algebraic varieties over number 
fields 
and to apply them to derive   a 
\emph{rationality criterion} for \emph{formal germs of functions}, 
which
extends the classical 
rationality theorems of  Borel-Dwork (\cite{eborel1894}, 
\cite{dwork1960}) and P\'olya-Bertrandias (\cite{polya1928},
\cite[Chapter 5]{amice75}; see also~\cite{cantor1980}),
valid over the projective line, 
to arbitrary algebraic curves over  a  number field.

Our algebraicity criteria improve on the ones in~\cite{bost2001} and 
\cite{bost2004}, which themselves were inspired by the 
papers~\cite{chudnovsky1985b}
and~\cite{chudnovsky1985a} of D.~V. and G.~V.~Chudnovsky and by the 
subsequent works by Andr\'e~\cite{andre1989} and Graftieaux 
\cite{graftieaux2001a,graftieaux2001b}.    As in~\cite{bost2001} and 
\cite{bost2004}, our results will be proved by means of a geometric 
version of  ``transcendence techniques'', which 
avoids the traditional constructions of  ``auxiliary polynomials'' and 
the explicit use of Siegel's Lemma,  replacing them by a few basic  
concepts of Arakelov geometry. In the proofs, our main objects of 
interest will be  some evaluation maps, defined on the spaces of global 
sections of powers of an ample line bundle on a projective 
variety by restricting these sections to formal subschemes or to 
subschemes of finite lengths. Arakelov geometry enters through the estimates satisfied by  
the heights of the evaluation maps, and the 
slopes and Arakelov degrees of the hermitian vector bundles defined 
by spaces of 
sections (see~\cite{chambert-loir2002} and~\cite{bost2007} for 
more details and references on this approach).

Our main motivation in investigating the algebraicity and rationality 
criteria presented in this article has been the desire to obtain 
 theorems respecting the classical principle of number theory that ``all
  places of number fields should appear on 
an equal footing'' ---
which actually is not the case in ``classical'' Arakelov 
geometry and in its applications in~\cite{bost2001}. 
A closely 
related aim has been to establish arithmetic theorems whose geometric 
counterparts (obtained through the analogy between number fields and 
function fields) have simple formulations and proofs.  These concerns led 
us to two technical developments in this paper: the use of (rigid) 
analytic geometry over~$p$-adic fields to define and estimate local 
invariants of formal curves over number fields\footnote{%
Since the first version of this paper was written, the relevance of rigid analytic geometry \emph{\`a la} Berkovich to develop a 
non-archimedean potential theory on~$p$-adic curves, and
consequently a ``modern''
version of Arakelov geometry of arithmetic surfaces satisfying the
above principle of ``equality of places'', has been largely demonstrated
by A.~Thuillier  in his thesis~\cite{thuillier2005}.},
and the derivation of a rationality criterion 
from an algebraicity criterion by means of the Hodge index theorem on 
(algebraic or arithmetic) 
surfaces.

%This paper is organized as follows. 

Let us describe the content of this article in more details.

In Section~\ref{sec.geometric},
we discuss geometric analogues of our arithmetic 
theorems. Actually, these are classical results in algebraic 
geometry, 
going back to Hartshorne~\cite{hartshorne1969} and  
Hironaka-Matsumura 
\cite{hironaka-m1968}. 
For instance, our algebraicity result admits as analogue the 
following fact.
 \emph{Let~$X$ be a quasi-projective variety over a
field~$k$, and let~$Y$ be a smooth projective integral curve in~$X$; 
let~$\hat{S}$ be  a smooth formal germ of 
surface through~$Y$} (that is, a smooth formal subscheme of dimension~$2$,
containing~$Y$, of the completion~$\hat{X}_Y$). \emph{If the 
degree~$\deg_Y  N_Y \hat{S}$ of the normal bundle to~$Y$ in~$\hat{S}$
is positive, then $\hat{S}$ is algebraic.}

Our point is 
that, transposed to a geometric setting, the arguments leading to 
our algebraicity and rationality criteria in the arithmetic setting 
--- which rely on 
the consideration of suitable evaluation maps, and on the Hodge index 
theorem --- provide simple proofs of these non-trivial 
algebro-geometric results, in which the geometric punch-line of the 
arguments appears more clearly.

In Section~\ref{sec.A-analytic},
we introduce the notion of \emph{$A$-analytic curve} 
 in an algebraic variety~$X$ over a number field~$K$.  By 
definition, 
this will be a smooth formal curve~$\hat{C}$ 
through a rational point~$P$ in~$X(K)$ --- that is, a 
smooth formal subscheme of dimension~$1$ in  the completion $\hat{X}_P$ --- 
which, firstly, is analytic at every place of~$K$, finite or infinite. 
Namely, if 
$v$ denotes any such place and $K_v$ the corresponding completion of~$K$,
the formal curve $\hat{C}_{K_v}$ in~$X_{K_v}$ deduced from 
$\hat{C}$ by the extension of scalars $K \hookrightarrow K_v$ is the 
formal curve defined by a $K_v$-analytic curve in~$X(K_v)$. 
Moreover the~$v$-adic radius~$r_v$ (in~$\mathopen]0, 1]$) of the open ball in 
$X(K_v)$ in which $\hat{C}_{K_v}$ ``analytically extends'' is required 
to ``stay close to~$1$ when~$v$ varies'', in the sense that the series 
$\sum_v \log r_v^{-1}$ has to be convergent. The precise formulation of 
this condition relies on the notion of \emph{size} of a smooth 
analytic germ in an algebraic variety over a~$p$-adic field.  This notion was
introduced in~\cite[\S3.1]{bost2001};  we review it
 in Section~\ref{subsec.size}, adding some complements. 

With the above notation, if $\cX$ is a model of~$X$ over the ring of 
integers~$\OK$ of~$K$, and if  $P$ extends to an integral point~$\cP$ 
in~$\cX(\OK)$, then a formal curve $\hat{C}$ through~$P$ is 
% an~$A$-germ 
$A$-analytic if it is analytic at each archimedean place of~$K$ and 
extends to a smooth formal surface~$\hat{\cC}$ in~$\hat{\cX}_\cP$. For 
a general formal curve~$\hat{C}$ that is analytic at archimedean places, 
being an~$A$-analytic germ may be seen as a weakened form of the 
existence of such a smooth extension $\hat{\cC}$ of~$\hat{C}$ along 
$\cP$. In this way, an~$A$-analytic  curve through the point~$P$ 
appears 
as an arithmetic counterpart of the smooth formal surface $\hat{S}$ 
along the curve~$Y$ in the geometric algebraicity criterion above.

The tools needed to formulate the arithmetic counterpart of the 
positivity condition  $\deg_Y N_Y \hat{S} > 0$ are developed in 
Sections~\ref{sec.can} and~\ref{sec.capa}.
We first show in Section~\ref{sec.can} how, 
for any germ of analytic curve~$\hat{C}$ 
through a rational point~$P$ in some algebraic variety~$X$ over a 
local field~$K$, one is led to introduce the so-called 
\emph{canonical semi-norm}  $\|\cdot\|^{\can}_{X, \hat{C}} $
on the~$K$-line $T_P \hat{C}$ through the 
consideration of the metric properties of the evaluation maps 
involved in our geometric version of the method of auxiliary 
polynomials. This extends a definition introduced in 
\cite{bost2004} when $K= \C$. 
In Section~\ref{sec.capa},
we discuss the construction of
Green functions and capacities on rigid analytic curves over~$p$-adic
fields. We then extend the comparison
of ``canonical semi-norms'' and ``capacitary metrics'' in  
\cite[\S3.4]{bost2004} to the non-archimedean setting.

In Section~\ref{sec.algebraic},
we apply these notions to formulate and establish our 
algebraicity results. If  $\hat{C}$ is an~$A$-analytic curve 
through a rational point~$P$ in an algebraic variety~$X$ over some 
number field
$K$, then the~$K$-line  $T_P \hat{C}$ may be equipped with  a 
``$K_v$-adic semi-norm'' for every place~$v$ by the above 
construction --- namely, the semi-norm $\|.\|^{\rm can}_{X_{K_v}, 
\hat{C}_{K_v}}$ on 
\[ T_P \hat{C}_{K_v} \simeq T_P \hat{C} \otimes_K K_v.\]
The so-defined metrized~$K$-line $\overline{T_P \hat{C}}$ has a 
well defined Arakelov degree 
%%% $\hdeg \overline{T_P \hat{V}}$
 in~$\mathopen]-\infty, +\infty]$, 
and our main  algebraicity
criterion asserts 
that \emph{$\hat{C}$ is algebraic if}  \emph{the Arakelov degree
$\hdeg \overline{T_P \hat{C}}$ is positive.} Actually, the converse
implication also holds: when $\hat{C}$ is algebraic,  the canonical
semi-norms  $\|.\|^{\rm can}_{X_{K_v}, 
\hat{C}_{K_v}}$ all vanish, and $\hdeg\overline{T_P \hat{C}}= +
\infty$.

Finally, in Section~\ref{sec.rat},
we derive an extension of the classical theorems 
of Borel, Dwork, P\'olya, and Bertrandias,
which gives a criterion for
the rationality of a formal germ of function~$\varphi$ on some
algebraic curve~$Y$ over a number field. By considering the graph
of~$\varphi$ --- a formal curve~$\widehat C$ in the surface $X:=Y \times \A ^1$ ---
we easily obtain the algebraicity of~$\varphi$ as a corollary of
our previous algebraicity criterion. In this way, we are reduced
to establishing a rationality criterion for an algebraic formal
germ. 
Actually, rationality results for algebraic functions
on the projective line have been investigated by
Harbater~\cite{harbater1988}, and used by Ihara~\cite{ihara1994}
to study the fundamental group of some arithmetic surfaces. Ihara's
results have been extended in~\cite{bost1999} by using Arakelov
geometry on arithmetic surfaces. Our rationality argument in Section~\ref{sec.rat},
based on the Hodge index theorem on arithmetic surfaces of
Faltings-Hriljac, is a variation on the proof of the Lefschetz
theorem on arithmetic surfaces in \emph{loc. cit.}.

It is a pleasure to thank A.~Ducros for his helpful advice on rigid analytic geometry during the preparation of this article.

Some of the results below have been presented, in a preliminary form,  
during the  ``Arithmetic Geometry and Number Theory'' conference  
in honor of N.~Katz, in Princeton, December 2003, 
and  have been announced in~\cite{bost2007}. 

During the preparation of this article, the authors have benefited from the support of the Institut universitaire de France. 

\medskip

It would be difficult to acknowledge fairly the multifaceted influences of Yuri Ivanovich Manin on our work. 
%It is an honour to dedicate this article to him,
%Yuri Ivanovich Manin.
%and we hope that it 
We hope that this article
will appear as a tribute, not only to his multiple contributions
to algebraic geometry and number theory, but also to his global
vision of mathematics, emphasizing geometric insights and analogies.
The presentation of this vision in his 25th-Arbeitstagung report
\emph{New directions in geometry}~\cite{manin1985} has been, since
it was written,  a source of wonder and inspiration to one of the
authors, and we allowed ourselves to borrow the terminology
``$A$-analytic'' from the ``$A$-geometry'' programmatically discussed
in~\cite{manin1985}. It is an honour for us  to dedicate this article
to Yuri Ivanovich Manin.

\section{Preliminary: the geometric case}
\label{sec.geometric}

The theorems we want to prove in this paper are analogues
in arithmetic geometry of classical algebro-geometric results
going back --- at least in an implicit form --- to  Hartshorne,
Hironaka, and Matsumura
(\cite{hartshorne1968}, \cite{hironaka1968}, \cite{hironaka-m1968}).
Conversely, in this section we give  short proofs of
algebraic analogues of our main arithmetic theorems.

\begin{proposition}\label{prop.geom-algebraicity}
Let $\mathscr X$ be a quasi-projective scheme over a field~$k$
and let~$\mathscr P$ be a projective connected subscheme
of dimension~$1$ in~$\mathscr X$.
Let $\widehat {\mathscr C}$ be a formal subscheme of dimension~$2$ in
$\widehat{\mathscr X}_{\mathscr P}$ admitting $\mathscr P$ as a
scheme of definition.
Assume that $\widehat {\mathscr C}$ is (formally) smooth over~$k$,
and that
$\mathscr P$ has no embedded component
(of dimension~$0$), or equivalently, that 
$\mathscr P$  defines a Cartier divisor in~$\widehat{\mathscr C}$,
and 
let $\mathscr N$ be the normal bundle 
%sheaf $\mathscr N_{\mathscrP/\widehat{\mathscr C}}$
 of the immersion 
$\iota:\mathscr P\hra\widehat{\mathscr C}$, that is, the invertible
sheaf $\iota^\ast\cO_{\hat{\cC}}(\cP)$ on~$\cP$.
%it is a line bundle on the scheme~$\mathscr P$.

If the divisor $[\mathscr P]$ on
the formal surface~$\widehat{\mathscr C}$
is \emph{nef} and has \emph{positive self-intersection}, 
%\footnote{Concerning elementary
%intersection theory  on formal
%surfaces in the style of Deligne \cite{deligne1972b},
%see Appendix~\ref{app.formalintersection}.}
then 
the formal surface $\widehat{\mathscr C}$ is algebraic, 
namely
the Zariski-closure of~$\widehat{\mathscr C}$ in~$\mathscr X$ is an
algebraic subvariety of dimension~$2$.
\end{proposition}

Let $(\mathscr P_i)_{i\in I}$ 
be  the family of irreducible components of~$\mathscr P$, 
and $(n_i)_{i\in I}$ their multiplicities in~$\mathscr P$. 
Recall that $[\mathscr P]$ is said to be nef on~$\widehat{\mathscr C}$ when 
 \[ [\mathscr P_i].[\mathscr P]:= \deg_{\mathscr P_i} \mathscr N \geq
0 \;\mbox{ for any $i\in I$,}\] 
and to have positive self-intersection
  if
 \[ [\mathscr P].[\mathscr P]:=\sum_{i\in I} n_i.\deg_{\mathscr P_i}
\mathscr N > 0,\]
or equivalently, when $[\mathscr P]$ is nef,
if one the non-negative integers $\deg_{\mathscr P_i} \mathscr N$ is positive.
Observe that 
%% $[\mathscr P]$ is \emph{nef} and \emph{big} %on~$\widehat{\mathscr C}$ 
these conditions are satisfied if $\mathscr N$ is ample on~$\mathscr P$. 

% Let us indicate that, , the proposition still holds with the smoothness
%assumption on~$\widehat{\mathscr C}$ omitted. 
% we leave this to the
%interested reader. 

More general versions of the  %\emph{algebraicity} 
algebraicity criterion in Proposition~\ref{prop.geom-algebraicity} 
and of its proof below,
without restriction on the dimensions of~$\widehat{\mathscr C}$ and
$\cP$,
can be found in~\cite[\S 3.3]{bost2001} \cite{bogomolov-mq2001},
\cite[Theorem~2.5]{bost2004} (see also~\cite{chambert-loir2002,chen2006}). 
Besides it will  be clear from the proof  that,
 suitably reformulated,
 Proposition~\ref{prop.geom-algebraicity}
still holds with the smoothness
assumption on~$\widehat{\mathscr C}$ omitted;  
 we leave this to the
interested reader. 

Such
algebraicity criteria may also be deduced from the works of Hironaka,
Matsumura, and Hartshorne on the condition~G$_2$
(\cite{hironaka1968}, \cite{hironaka-m1968}, \cite{hartshorne1968}).
We refer the reader to the monographs~\cite{hartshorne1970} and \cite{badescu2004} for extensive
discussions and references about related results concerning formal
functions and projective algebraic varieties.

Note that Proposition~\ref{prop.geom-algebraicity} has consequences for the
study of
algebraic varieties over function fields. Let indeed~$S$
be a smooth, projective, and geometrically connected curve over a
field~$k$
and let $K=k(S)$.  Let $f\colon \mathscr X\ra S$ be a surjective map
of~$k$-schemes and
assume that $\mathscr P$ is the image
of a section of~$f$. Let $X=\mathscr X_K$, $P=\mathscr P_K$
and $\widehat C=\widehat{\mathscr C}_K$ be the generic fibers of
$\mathscr X$, $\mathscr P$, and $\mathscr C$.
Then~$P$ is a~$K$-rational point of~$X$
and $\widehat C$ is a germ of curve in~$X$ at~$P$.
Observe that $\widehat {\mathscr C}$ is algebraic
if and only if $\widehat C$ is algebraic. Consequently, in this
situation, Proposition~\ref{prop.geom-algebraicity} appears as an
algebraicity criterion for a formal germ  of curves $\widehat C$ in
$X$. In particular, it shows that such \emph{a smooth formal 
curve~$\widehat C$ in $X$ is algebraic 
if it extends to a smooth formal scheme $\widehat{\mathscr C}$  
through $\mathscr P$ in~$\mathscr X$ such that the normal bundle 
%the normal bundle %sheaf $\mathscr N_{\mathscrP/\widehat{\mathscr
%C}}$
of  $\mathscr P$ in~$\widehat{\mathscr C}$ has positive degree.}

\begin{proof}[of Proposition~\ref{prop.geom-algebraicity}]
We may assume that $\mathscr X$ is projective and that
$\widehat{\mathscr C}$ is Zariski-dense in~$\mathscr X$.
We let $d=\dim\mathscr X$.
One has obviously $d\geq 2$ and our goal is to prove the equality.

Let $\mathscr O(1)$ be any very ample line bundle on~$\mathscr X$.
The method of ``auxiliary polynomials'', borrowed from transcendence
theory, suggests the introduction of the ``evaluation maps''
\[ \phi_D \colon \Gamma(\mathscr X,\mathscr O(D)) \ra
\Gamma(\widehat{\mathscr C},\mathscr O(D)),
\quad s\mapsto s|_{\widehat{\mathscr C}}, \]
for positive integers~$D$.

%Since $\widehat{\mathscr C}$ is assumed to be Zariski-dense
%in~$\mathscr
%X$, the map $\pi_D$ is injective
%and consequently 

Let us denote $E_D=\Gamma(\mathscr X,\mathscr O(D))$, and
for any integer~$i\geq 0$, let $E^i_D$ be
the set of all $s\in E_D$ such that $\phi_D(s)=s|_{\widehat{\mathscr
C}}$
vanishes at order at least~$i$ along~$\mathscr P$, \ie  such that the
restriction of~$\phi_D(s)$ to $i\mathscr P$ vanishes.
Since $\widehat{\mathscr C}$ is Zariski-dense in~$\mathscr X$,
no non-zero section of  $\mathscr O(D)$ has a restriction to
$\widehat{\mathscr C}$ that vanishes at infinite order
along~$\mathscr P$,  and we have 
\[ \bigcap_{i=0}^\infty E^i_D = 0. \]
Consequently,
\[\rank E_D = \sum_{i=0}^\infty \rank (E^i_D/E^{i+1}_D) .\]

Moreover, there is a canonical injective map of~$k$-vector spaces
\[ E^i_D/E^{i+1}_D \hra \Gamma(\mathscr P, \mathscr O(D)\otimes
\mathscr N^{\vee\otimes i}), \]
which amounts to taking the~$i$th jet along~$\mathscr P$ --- that is,
the restriction to $(i+1)\mathscr P$ ---
of a section which vanishes at order at least~$i$. Indeed the
quotient sheaf 
\[ 
\left( 
  \mathscr O(D) \otimes \mathscr O_{\widehat{\mathscr C}}(-i\mathscr P)
\right) / \left(
  \mathscr O(D) \otimes \mathscr O_{\widehat{\mathscr C}}(-(i+1)\mathscr P)
\right) \] 
over~$\widehat{\mathscr C}$
may be identified with $\mathscr O(D) \otimes \iota_\ast \mathscr
N^{\vee\otimes i}$. 
Observe also that the dimension of the range of this injection
satisfies
an upper bound of the form
\[ \dim \Gamma(\mathscr P, \mathscr O(D)\otimes
{\mathscr N}^{\vee\otimes i}) \leq c(D+i), \]
valid for any non-negative integers~$D$ and~$i$. 

Assume that $E_D^i\neq 0$ and let $s\in E_D^i$ be any nonzero element.
By assumption, $\phi_D(s)$ vanishes at order~$i$ along~$\mathscr P$,
hence $\div\phi_D(s)- i[\mathscr P]$ is an effective divisor on
$\widehat{\mathscr C}$
and its intersection number with~$[\mathscr P]$
is nonnegative,
for $[\mathscr P]$ is nef. Consequently
\[ \div\phi_D(s)\cdot [\mathscr P] \geq i [\mathscr P]\cdot [\mathscr P].\]
Since 
\[ \div\phi_D(s)\cdot[\mathscr P]=\deg_{\mathscr P}(\mathscr O(D))
= D\deg_{\mathscr P}(\mathscr O(1)) \]
and $[\mathscr P]\cdot[\mathscr P]>0$ by the  assumption of positive self-intersection, this
implies
 $i\leq aD$, where $a:=\deg_{\mathscr P}\mathscr O(1)
/ [\mathscr P].[\mathscr P]$. Consequently $E^i_D$ is reduced to~$0$
if $i>aD$.

Finally, we obtain:
\[
 \rank E_D  =  \sum_{i=0}^\infty \rank (E^i_D/E^{i+1}_D) 
 =  \sum_{i=0}^{\lfloor aD\rfloor} \rank (E^i_D/E^{i+1}_D) 
\leq \sum_{i=0}^{\lfloor aD\rfloor } c(D+i). 
\]
This proves that, when~$D$ goes to $+\infty$,
\[ \rank E_D \ll D^2.\]
Besides
%when~$D$ goes to $+\infty$,
\[ \rank E_D = \rank \Gamma(\mathscr X,\mathscr O(D))
 \asymp D^{d}, \]
by 
Hilbert-Samuel's theorem.
This establishes that the integer~$d$, which is at least~$2$,
actually equals~$2$.
%Necessarily $d=2$, as was to be shown.
\end{proof}

\begin{proposition}\label{prop.geom-hodge}
Let $f\colon S'\ra S$ be a dominant morphism between 
two normal projective surfaces over a field~$k$.
Let $D\subset S$ and $D'\subset S'$ be effective divisors
such that $f(D')=D$.

Assume that $f|_{D'}\colon D'\ra D$ is an isomorphism and that
$f$ induces an isomorphism $\widehat f\colon \widehat {S'_{D'}}\ra
\widehat {S_D}$ between formal completions. 
If moreover~$D$ is nef and  $D\cdot D>0$, then~$f$ is birational.
\end{proposition}

Recall that~$D$ is said to be nef if, for any effective divisor~$E$ on~$S$ the (rational) intersection number
$D\cdot E$ is non-negative.

\begin{proof}
By hypothesis, $f$ is \'etale in a neighborhood of~$D'$.
If $\deg(f)>1$, one can therefore write $f^*D=D'+D''$,
where $D''$ is a non zero effective Cartier divisor on~$S'$ which is
disjoint
from~$D'$.
Now, $f^*D$ is a nef divisor on~$S'$
such that $f^*D\cdot f^*D=\deg(f)D\cdot D>0$.
As a classical consequence of the Hodge index theorem
(see~\cite{franchetta1949}, \cite{ramanujam1972},
and also~\cite[Proposition~2.2]{bost1999})
the effective divisor~$f^*D$
is numerically connected, hence connected.
This contradicts the decomposition $f^*D=D'\sqcup D''$.
\end{proof}

\begin{proposition}\label{prop.geom-rationality}
Let $\mathscr S$ be a smooth projective connected surface over
a perfect field~$k$. Let $\mathscr P$ be a smooth projective
connected curve in~$\mathscr S$. 
If the divisor $[\mathscr P]$ on~$\mathscr S$ is big and nef,
then any formal rational function along~$\mathscr P$ 
is defined by a (unique) rational function on~$\mathscr S$.
In other words, one has an isomorphism of fields
\[ k(\mathscr S) \xrightarrow\sim \Gamma(\mathscr P,\operatorname{Frac}
\mathscr O_{\widehat{\mathscr S}_{\mathscr P}}). \]
\end{proposition}
\begin{proof}
Let $\phi$ be any formal rational function along~$\mathscr P$.
We may introduce a sequence of blowing-ups of closed points 
$\nu\colon{\mathscr S}'\ra  \mathscr S$
such that $\phi'=\nu^*\phi$ has no point of indetermination
and may be seen as a map (of formal~$k$-schemes)
$\widehat{\mathscr S}'_{{\mathscr P'}}\ra\P^1_k$,
where ${\mathscr P}'=\nu^*\mathscr P$.

Let us consider the graph $\operatorname{Gr}\phi'$ of~$\phi'$
in~${\mathscr S}'\times\P^1_k$. This is a formally smooth
2-dimensional formal scheme, admitting the graph of  $\phi'_{\mid
{\mathscr P}'}\colon
{\mathscr P'}\ra\P^1_k $
as a scheme of definition, and the morphism $\phi'$ defines an isomorphism of formal schemes 
\[ \psi'\mathrel{\mathord:\mathord=}(\Id, \phi)
\colon \widehat{\mathscr S'}_{\mathscr P'} \ra \operatorname{Gr} \phi'.\]
Like the divisor ${\mathscr P}$ 
in~${\mathscr S}$, its inverse image  ${\mathscr P}'$ 
in~${\mathscr S}'$ is nef and has positive self-intersection.
Proposition~\ref{prop.geom-algebraicity}
therefore  implies that $\text{Gr}\, \phi'$ %the graph of~$\phi'$ 
 is algebraic in~${\mathscr S}'\times\P^1_k$.
In other words, $\phi'$ is an algebraic function.

To establish its rationality, let us introduce  
the Zariski closure $\Gamma$  of
the graph of~$\text{Gr}\,\phi'$ in~${\mathscr S}'\times\P^1_k$, 
the projections $\mathord{\pr_1}\colon \Gamma
\ra {\mathscr S}'$ and 
$\pr_2\colon \Gamma \ra \P^1_k$, and  the normalization
$n\colon\tilde{\Gamma}  \ra \Gamma$ of~$\Gamma$. 
Consider also the Cartier divisor $\mathscr P '_\Gamma$ 
(resp. $\mathscr P '_{\tilde{\Gamma}}$) defined 
as the inverse image $\pr_1^*\mathscr P'$ 
(resp. $n^* \mathscr P '_\Gamma$) 
of~$\mathscr P '$ in~$\Gamma$
(resp. $\tilde{\Gamma}$). 
The morphisms~$n$ and $\pr_1$ define
morphism of formal completions:
\[ \widehat{\tilde{\Gamma}}_{\mathscr P '_{\tilde{\Gamma}}}
\stackrel{\hat{n}}{\longrightarrow}
\widehat{{\Gamma}}_{\mathscr P '_{{\Gamma}}}
\stackrel{\widehat{\pr_1}}{\longrightarrow}
\widehat{\mathscr S'}_{\mathscr P'}.\]
The morphism $\psi'$ may be seen a section of~$\widehat{\pr_1}$;
% \[ \psi\colon \widehat{\mathscr S'_{\mathscr P'}}
% \ra \widehat{{\Gamma}}_{\mathscr P '_{{\Gamma}}}.\]
by normality of~$\widehat{\mathscr S'_{\mathscr P'}}$, 
it admits a factorization through~$\hat{n}$ of the form
$\psi'=\hat{n}\circ \tilde{\psi}$, for some uniquely determined morphism 
of~$k$-formal schemes
$\tilde{\psi}\colon \widehat{\mathscr S'}_{\mathscr P'} \ra \widehat{\tilde{\Gamma}}_{\mathscr P '_{\tilde{\Gamma}}}$.
This morphism~$\tilde\psi$
is a section of~$\mathord{\widehat{\pr_1}} \circ \hat{n}$. 
Therefore the (scheme theoretic) image $\tilde{\psi}(\mathscr P ')$ 
%of~$\mathscr P '$
  defines a (Cartier) divisor in~$\tilde{\Gamma}$ such that 
  \[ (f\colon S' \ra S, D',D)= (\mathord{\pr_1} \circ n \colon \tilde{\Gamma} \ra \mathscr S', \tilde{\mathscr P}, \mathscr P ') \]
  satisfy the hypotheses of Proposition~\ref{prop.geom-hodge}.  Consequently the morphism~$\mathord{\pr_1} \circ n$ 
is birational. Therefore, $\pr_1$ is birational too
and $\phi'$ is the restriction of a rational function
on~$\mathscr S'$, namely $\mathord{\pr_2} \circ \mathord{\pr_1^{-1}}$.
This implies that~$\phi$ is the restriction of 
a rational function on~$\mathscr S$.
The uniqueness of this rational function follows from the Zariski density of the formal neighborhood of~$\mathscr P$ in~$\mathscr S$. 
\end{proof}

\begin{remark}
In the terminology of Hironaka and Matsumura~\cite{hironaka-m1968},
the last proposition asserts that
$\mathscr P$ is~G\textsubscript3 in~$\mathscr S$, and has been established by
Hironaka
in~\cite{hironaka1968}.
%Under such an assumption, 
Hartshorne observes in~\cite[Proposition~4.3, and Remark p.~123]{hartshorne1969}
that Proposition~\ref{prop.geom-hodge} holds more generally under the
assumption that 
$D$ and~$D'$ are~G\textsubscript 3 in~$\mathscr S$
and~$\mathscr S'$.
Our approach to Propositions~\ref{prop.geom-hodge} 
and~\ref{prop.geom-rationality} 
follows an order opposite to the
one in~\cite{hironaka1968} and~\cite{hartshorne1969}, and actually
provides a simple proof of~\cite[Proposition~4.3]{hartshorne1969}.    
\end{remark} 

\section{$A$-analyticity of formal curves}
\label{sec.A-analytic}
\subsection{Size of smooth formal curves over~$p$-adic fields}
\label{subsec.size}

In this Section,
we briefly recall some definitions and results from~\cite{bost2001}.

Let~$K$ be field equipped with some complete ultrametric absolute
value $\abs{.}$ 
and assume that its valuation ring~$R$ is a discrete valuation ring.
Let also $\overline K$ be an algebraic closure of~$K$. We shall still
denote $\abs{.}$ the non-archimedean absolute value on~$\overline K$
that extends the absolute value $\abs{.}$ on~$K$. 

For any positive real number~$r$,
we define the norm $\norm{g}_r$ of a formal power 
series $g=\sum_{I\in\N^N}a_IX^I \in K[[X_1,\dots,X_n]]$
by the formula
\[ \norm{g}_r= \sup_I \abs{a_I} r^{\abs I}\quad; \]
it belongs to $\R_+\cup\{\infty\}$.
The power series~$g$ such that~$\norm{g}_r<\infty$
are precisely those that are convergent and bounded
on the open~$N$-ball of radius~$r$ in~$\bar K^N$.

The group~$G_{\text{for},K}:=\Aut(\widehat{\A^N_{K,0}})$ of
automorphisms
of the formal completion of~$\A^N_K$ at~$0$
may be identified with the set of all~$N$-tuples $f=(f_1,\dots,f_n)$
of power series in~$K[[X_1,\dots,X_N]]$ such that
$f(0)=0$ and $Df(0):=\left(\frac{\partial f_j}{\partial
X_i}(0)\right)$
belongs to $\GL_N(K)$.
We consider its following subgroups:
\begin{itemize}
\item the subgroups $G_{\text{for}}$ 
consisting of all elements $f\in G_{\text{for},K}$
such that $Df(0)\in\GL_N(R)$;
\item the subgroup $G_{\text{an},K}$
consisting of those $f=(f_1,\dots,f_N)$ in $G_{\text{for},K}$
such that, for each~$j$, $f_j$ has a positive radius of convergence;
\item $G_{\text{an}}:=G_{\text{an},K}\cap G_{\text{for}}$;
\item for any positive real number~$r$, the subgroup
$G_{\text{an},r}$ of~$G_{\text{an}}$ consisting of all~$N$-tuples
$f=(f_1,\dots,f_N)$ such that $\norm{f_j}_r\leq r$ for each~$j$.
This subgroup may be identified with the group of all analytic
automorphisms, preserving the origin, of the open~$N$-dimensional
ball of radius~$r$.
\end{itemize}

One has the inclusion 
$  G_{\text{an},r'}\subset G_{\text{an},r}$ for any
$r'>r>0$,  and  the equalities
\[ \bigcup_{r>0} G_{\text{an},r}= G_{\text{an}}\quad\text{and}\quad
        G_{\text{an},1}=\Aut\big( \widehat{\A^N_{R,0}}). \]

\medskip

It is straightforward that a formal subscheme $\widehat V$
of~$\widehat{\A^N_{K,0}}$
is (formally) \emph{smooth} of dimension~$d$ iff there exists $\phi\in
G_{\text{for},K}$
such that $\phi^*\widehat V$ is the formal subscheme
$\widehat{\A^d_{K,0}}\times\{0\}$ of~$\widehat{\A^N_{K,0}}$;
when this holds, one can even find such a $\phi$ in~$G_{\text{for}}$.
Moreover such a smooth formal subscheme~$\widehat V$
is~$K$-\emph{analytic} iff one can find~$\phi$ as above
in~$G_{\text{an},K}$,
or equivalently in~$G_{\text{an}}$.

\medskip

Let $\mathscr X$ be a flat quasi-projective~$R$-scheme, and
 $X=\mathscr X\otimes_R K$  its generic fibre.
Let $\mathscr P\in \mathscr X (R)$ be an~$R$-point, 
and let $P\in X(K)$ be its restriction to $\Spec K$.
In~\cite[\S3.1.1]{bost2001} 
we associated to any smooth formal scheme $\widehat V$
of dimension~$d$  in~$\widehat X_P$, its \emph{size~$S_{\mathscr X}(\widehat
V)$ with respect to the model~${\mathscr X}$} of~$X$. It is a number in~$[0,1]$ whose definition and
basic properties  may be summarized in the following statement:

%Finally, let  be a smooth
%with $\{P\}$ as a scheme of definition.
%The 
%of such a triple $(\mathscr X,\mathscr P,\widehat V)$
%is defined thanks to the following theorem, which
%follows from

\begin{theorem}\label{theo.def-size}
There is a unique way to attach a number $S_{\mathscr X}(\widehat V)$
in~$[0,1]$ 
to any such data $(\mathscr X,\mathscr P,\widehat V)$ so that
the following properties hold:
\begin{enumerate}
\def\theenumi{\alph{enumi}}
\def\labelenumi{\theenumi)}
\item 
if $\mathscr X\ra\mathscr X'$ is an immersion,
then $S_{\mathscr X'}(\widehat V)=S_{\mathscr X}(\widehat V)$
\emph{(invariance under immersions)};
\item for any two triples  $(\mathscr X,\mathscr P, \widehat V)$ and 
$(\mathscr X',\mathscr P',\widehat V')$  as above,
if there exists an~$R$-morphism 
$\phi\colon\mathscr X\ra\mathscr X'$ mapping $\mathscr P$
to~$\mathscr P'$,
\'etale along~$\mathscr P$, and inducing an isomorphism
$\widehat V\simeq\widehat V'$, then $S_{\mathscr X'}(\widehat V')
=S_{\mathscr X}(\widehat V)$
\emph{(invariance by \'etale localization)};
\item 
if $\mathscr X=\A^N_R$ is the affine space over~$R$ and 
$\mathscr P=(0,\dots,0)$, then $S_{\mathscr X}(\widehat V)$
is the supremum in~$[0,1]$
of the real numbers $r\in (0,1]$ for which there exists
$f \in G_{{\rm an},r}$ 
such that
$f^*\widehat V =\widehat{\mathbf A^d_0}\times\{0\}$
\emph{(normalization)}.
\end{enumerate}
\end{theorem}

%\begin{definition}
%The real number $S_{\mathscr X}(\widehat V)$ 
%defined by Theorem~\ref{theo.def-size} is called
%the size of~$\widehat V$ with respect to~$\mathscr X$.
%\end{definition}

As a straightforward consequence of these properties of the size, we
obtain:

\begin{proposition}\label{prop.size-analytic}
A smooth formal subscheme~$\widehat V$ in~$\widehat X_P$ is
$K$-analytic
if and only if its size~$S_{\mathscr X}(\widehat V)$ 
is a positive real number.
\end{proposition}

\begin{proposition}\label{prop.size-smooth}
Let $\mathscr X$, $\mathscr P$, and $\widehat V$ be as above
and assume that there exists a smooth formal~$R$-subscheme $\mathscr
V\subset
\widehat{\mathscr X_{\mathscr P}}$ such that $\widehat V=\mathscr
V_K$.
Then $S_{\mathscr X}(\widehat V)=1$.
\end{proposition}

 The remainder of this section is devoted to further properties of
the size.
%Let us conclude this Section with two remarks on sizes.

\begin{proposition}\label{prop.size-extensions}
The size is invariant under isometric extensions of  valued fields (complete with respect to a discrete valuation).
\end{proposition}
\begin{proof}
It suffices to check this assertion in the case
of a smooth formal subscheme~$\widehat V$ through the origin of the
affine space $\A^N$.
By its very definition, the size cannot decrease under extensions
of the base field. 

To show that it cannot increase either, let us  fix an isomorphism of~$K$-formal schemes
   \[ \xi=(\xi_1,\ldots,\xi_N): \widehat
\A^d_0\stackrel{\sim}{\longrightarrow} \widehat V
   \hookrightarrow \widehat \A^N_0 \]
 given by~$N$ power series
$\xi_i\in K[[T_1,\dots,T_d]]$ such that $\xi_1(0)=\ldots=\xi_N(0)=0$.
We then observe that, for any~$N$-tuple $g=(g_1,\ldots,g_N)$ of series in~$K[[X_1,\dots,X_N]]$, the following two conditions are equivalent:
\begin{enumerate}\def\theenumi{\roman{enumi}}
\item~$g$ belongs to  $G_{\text{for},K}$ and $(g^{-1})^* \widehat{V}=\widehat{\mathbf A^d_0}\times\{0\}$;
\item $g_1(0)=\cdots=g_N(0)=0$, $g_{d+1}(\xi_1,\ldots,\xi_N)=\cdots=g_{N}(\xi_1,\ldots,\xi_N)=0$, and $\left(\frac{\partial g_i}{\partial
X_j}(0)\right)$
belongs to $\GL_N(K)$.
\end{enumerate}

Let $K'$ be a valued field, satisfying the same condition as~$K$,
that contains~$K$ and whose absolute value restricts to the given one on
$K$. Let $R'$ be its valuation ring. We shall denote $G'_{\rm for}, G'_{{\rm an}}, G'_{{\rm an},r},\dots$ the analogues of~$G_{\rm for}, G_{{\rm an}}, G_{{\rm an},r},\ldots$ defined by replacing the valued field~$K$ by $K'$.
Recall that there exists an ``orthogonal projection'' from $K'$ to~$K$, namely a~$K$-linear map $\lambda\colon K'\ra K$ such that 
$\abs{\lambda(a)} \leq \abs{a}$ for any $ a\in K'$
and 
$\lambda(a)=a $ for any~$a\in K$;
see for instance~\cite[p.~58, Corollary~(2.3)]{gruson-vdp1974}.

Let $\widehat V'= \widehat V_{K'}$ be the formal subscheme of
$\widehat{\A^N_{K'}}$ deduced from $\widehat V$ by the extension
of scalars $K \hookrightarrow K'$,  and let~$r$ be  an element  in
$\mathopen]0,S_{ \A^N_{R'}}(\widehat V')\mathclose[$. By the very definition of the
size, there exists some $g'=(g'_1,\ldots,g'_N)$ in~$G'_{{\rm an},r}$
such that $(g'^{-1})^* \widehat{V}=\widehat{\mathbf A^d_0}\times\{0\}$.
Since the tangent space at the origin of  $V'$ is defined over~$K$,
by composing $g'$ with a suitable element in~$GL_N(R')$, we may
even find  $g'$ such that $Dg'(0)$ belongs to $GL_N(R)$.
Then the series $g_i:= \lambda \circ g'_i$, deduced from the series
$g'_i$ by applying the linear map $\lambda$ to their coefficients,
satisfy $g_i(0)=0$, 
$(\partial g_i/\partial X_j)(0) =(\partial g'_i/\partial X_j)(0)$ 
and $\norm{g_i}_r \leq \norm{g'_i}_r$.
Therefore $g:=(g_1,\ldots,g_N)$ is an element of~$G_{{\rm an},r}$.
Moreover, from the equivalence of conditions (i) and (ii) above and
its analogue with $K'$ instead of~$K$, we derive that~$g$ satisfies
$(g^{-1})^* \widehat{V}=\widehat{\mathbf A^d_0}\times\{0\}$. 
This shows that $S_{ \A^N_R}(\widehat V) \geq r$ and establishes the
required inequality $S_{ \A^N_R}(\widehat V) \geq  S_{ \A^N_{R'}}(\widehat
V)$.
\end{proof}

The next proposition relates sizes, radii of convergence, and Newton
polygons.
%\begin{remark}[About sizes, radii of convergence and Newton
%polygons]\label{rema.size-newton}

\begin{proposition}\label{prop.size-radius-newton}
Let~$\phi\in K[[X]]$ be a power series such that $\phi(0)=0$ and
$\phi'(0)\in R$, and let $\widehat C$ be its graph, namely the formal
subscheme of  $\widehat \A^2_0$ defined by the equation
$x_2=\phi(x_1)$.

1) 
The radius of convergence $\rho$ of~$\phi$ satisfies
\[ \rho \geq S_{ \A^2_R}(\widehat C). \]

2)  Suppose that $\rho$ is positive and that $\phi'(0)$ is a unit in
$R$. Then
 \[ S_{ \A^2_R}(\widehat C)=\min(1,\exp \lambda_1),\]
 where~$\lambda_1$ denotes
 the first slope of the Newton polygon of the power
series~$\phi(x)/x$.

\end{proposition}

Recall that, if $\phi= \sum_{i\geq 1} c_i T^i$, under the hypothesis
in  2),
we have:
\[ \lambda_1:=\inf_{i\geq 1}- \frac{\log \abs{c_{i+1}}}{i} \leq \lim
\inf_{i \ra + \infty} 
- \frac{\log \abs{c_{i+1}}}{i} = \log \rho.\]
Moreover $\exp \lambda_1$ is the supremum of the numbers $r \in
]0,\rho[$ such that, for any~$t$ in~$\overline K$ satisfying $\abs{t}
< r$, we have $\abs{\phi(t)} =\abs{t}$.

% when $\rho$ is positive,  $\exp \lambda_1$ is the supremum of
%numbers~$r$ in~$]0,\rho[$ such that $%\abs{\phi(t)} = \abs{t}$ for
%any $t\in \overline K$ satisfying $\abs{t} < r$.
%Let $\phi\in K[[x]]$ be a power series with positive radius of
%convergence such that $\phi(0)=0$ and $\phi'(0)\in R$.
%Let $\widehat C$ be its formal graph in~$\A^2_0$, defined
%by the equation $x_2=\phi(x_1)$; 

\begin{proof}

Let~$r$ be a positive real number
such that $r<S_{\A^2_0}(\widehat C)$.
By assumption, there are power series $f_1$ and
$f_2\in K[[X_1,X_2]]$ such that $f=(f_1,f_2)$ belongs
to~$G_{\text{an},r}$
and such that $f^*\widehat C=\widehat{\A^1}\times\{0\}$. This last
condition implies (actually is equivalent to) the identity in 
$K[[T]]$:
\[ f_2(T,0)=\phi(f_1(T,0)).\]
Let us write $f_1(T,0)=\sum_{i \geq 1} a_{i}T^i $,
$f_2(T,0)=\sum_{i \geq 1} b_{i}T^i$,
and $\phi(X)=\sum_{i\geq 1} c_i X^i$. 

One has $b_1=c_1a_1$ and $c_1=\phi'(0)$ belongs to~$R$ by hypothesis.
Moreover, the first column of the matrix~$Df(0)$ is
$\left(\begin{smallmatrix} a_1\\b_1\end{smallmatrix}\right)
=a_1\left(\begin{smallmatrix} 1 \\ c_1 \end{smallmatrix}\right)$.
Since $Df(0)$ belongs to $\GL_2(R)$ and $c_1$ to~$R$,
this imply that $a_1$ is a unit in~$R$.
Then, looking at the expansion of~$f_1(T,0)$ (which satisfies
$\norm{f_1(T,0)}_r \leq r$), we see
that $\abs {f_1(t,0)}=\abs t$ for any  $t\in \overline K$ such that
$\abs t<r$. Consequently, if~$g\in K[[T]]$ denotes the reciprocal
power series of~$f_1(T,0)$, then~$g$
converges on the open disc of radius~$r$ and
 satisfies $\abs{g(t)}=\abs t$ for any $t \in \overline K$ such that
$\abs t<r$.
 
 The identity in~$K[[T]]$
\[ \phi(T) = \phi( f_1(g(T),0)) = f_2(g(T),0) \]
then shows that
the radius of convergence of~$\phi$ is at least~$r$. This establishes~1).

Let us now assume that $\rho$ is positive and that  $\phi'(0) (=c_1)$
is a unit of~$R$. Then $b_1= a_1 c_1$
also is a unit and similarly, we have $\abs{f_2(t,0)}=\abs{t}$
for any $t \in \overline K$ such that $\abs t<r$.  This implies that
$\abs{\phi(t)}=\abs{t}$ for any such~$t$. This shows that $\exp
\lambda_1 \geq S_{ \A^2_R}(\widehat C)$.

To complete the proof of 2),  observe that the element~$f$ of~$G_{\rm
an}$ defined as $f(T_1,T_2)=(T_1+T_2, \phi(T_1))$ satisfies
$f^*\widehat C=\widehat{\A^1}\times\{0\}$ and belongs to $G_{{\rm
an},r}$ for any~$r$ in~$\mathopen]0, \min(1, \exp \lambda_1)\mathclose[$.
%In other words, \emph{the Newton polygon of~$\phi(x)/x$ is
%above the line of slope~$\log(r)$.} 
%We conclude: \emph{Let $\widehat\C$ be the formal graph of
%a power series~$\phi\in K[[x]]$ be a power series
%such that $\phi(0)=0$ and $\phi'(0)\in R^*$.}
\end{proof}

%One has $S_{\widehat \A^2}(\widehat C)=\min(1,r,\exp(\lambda_1))$,
%where~$r$ is the radius of convergence of~$\phi$ and~$\lambda_1$
%is the first slope of the Newton polygon of the power
%series~$\phi(x)/x$.}

%\begin{proposition}\label{prop.newton}
%Let~$\phi\in K[[x]]$ be a power series such that $\phi(0)=0$ and
%$\phi'(0)\in R^*$, and let $\widehat C$ be its graph, namely the
%formal subscheme of  $\widehat \A^2_0$ defined by the equation
%$x_2=\phi(x_1)$. If $\rho$ denotes the radius of convergence
%of~$\phi$ and~$\lambda_1$
% the first slope of the Newton polygon of the power
%series~$\phi(x)/x$, then 
% \[ S_{ \A^2_R}(\widehat C)=\min(1,\rho,\exp \lambda_1).\]
%\end{proposition}

Observe that for any non zero $a\in R$, 
the series $\phi(T)=T/(a-T)$ has radius
of convergence $\rho=\abs{a}$ while the size of its graph $\widehat
C$ is 1 (observe that 
$f(T_1,T_2):=(aT_1+T_2,T_1/(1-T_1))$ satisfies $f^*\widehat
C=\widehat{\A^1}\times\{0\}$).
Taking $\abs a<1$, this shows
that  the size of the graph of
a power series~$\phi$ can be larger than its radius of convergence 
when the assumption $\phi'(0)\in R$ is omitted.

As an application of  the second assertion in
Proposition~\ref{prop.size-radius-newton},  we obtain that, when~$K$ is a
$p$-adic field,  the size of the graph of~$\log(1+x)$
is equal to~$\abs p^{1/(p-1)}$.
%where~$p$ is the characteristic of the residue field of~$R$.
Considering this graph as the graph of the exponential
power series with axes exchanged, this also follows from the first
assertions of 
Proposition~\ref{prop.size-radius-newton} and %, combined with
 Proposition~\ref{prop.size-foliation} below.
%\end{remark}

Finally, let us indicate that, by analyzing the construction
\emph{\`a la} Cauchy
of local solutions of analytic ordinary differential equations, one
may establish
the following lower bounds
on the size of a formal curve obtained by integrating
an algebraic one-dimensional foliation over a~$p$-adic field
(\cf \cite[Proposition~4.1]{bost2004}):
\begin{proposition}\label{prop.size-foliation} Assume that~$K$ is a field of
characteristic~$0$, and that its residue field~$k$ has positive
characteristic~$p$.
Assume also that $\mathscr X$ is smooth over~$R$ in a neighborhood
of~$\mathscr P$.
Let $\mathscr F\subset T_{\mathscr X/R}$ be a rank~$1$ subbundle
and let  $\widehat C$ be the formal integral curve through~$P$ of the
one-dimensional foliation $F=\mathscr F_K$.
Then
 \[ S_{\mathscr X}(\widehat C)\geq \abs p^{1/(p-1)}.\]
If moreover~$K$ is absolutely unramified %over~$\Q_p$
 \emph{(that is, if the maximal ideal of~$R$ is $pR$)} and   
if the one-dimensional subbundle $\mathscr F_{k}\subset T_{\mathscr
X_k}$ 
is closed  under~$p$-th power, then
 \[ S_{\mathscr X}(\widehat C)\geq \abs p^{1/p(p-1)}. \] 
\end{proposition}

\subsection{$A$-analyticity of formal curves in algebraic varieties
over number fields}
\label{sec.A-analyticity}

Let~$K$ be a number field and let~$R$
denote its ring of integers. For any maximal ideal
$\mathfrak p$ of~$R$, 
let $\abs{\cdot}_{\mathfrak p}$ denote the $\mathfrak p$-adic
absolute value, normalized by the
condition $\abs{\pi}_{\mathfrak p}=(\#(R/\mathfrak p))^{-1}$
for any uniformizing element~$\pi$ at~$\mathfrak p$.
Let $K_{\mathfrak p}$ and~$R_{\mathfrak p}$
be the $\mathfrak p$-adic completions of~$K$ and~$R$, and $\F_{\mathfrak
p}:=R/{\mathfrak p}$ the residue field of~${\mathfrak p}$.

In this Section, we consider a quasi-projective algebraic variety
$X$ over~$K$, a rational point~$P$ in~$X(K)$,
and   a smooth formal curve $\widehat C$ in~$\widehat X_P$.
% through~$P$ in~$X$..

It is straightforward that, if~$N$ denotes a sufficiently divisible
positive integer, there exists a model $\mathscr X$ of~$X$,
quasi-projective over~$R[1/N]$, such that~$P$ extends to a point
$\mathscr P$ in~$\mathscr X(R[1/N])$. 
Then, for any maximal ideal $\mathfrak p$ not dividing~$N$,
the size $S_{\mathscr X_{R_{\mathfrak p}}}(\widehat C_{K_{\mathfrak
p}})$
is a well-defined real number in~$[0,1]$.

\begin{definition}
We will say that the formal curve $\widehat C$ in~$X$
is \emph{$A$-analytic} if the following conditions are satisfied:
\begin{enumerate}\def\theenumi{\roman{enumi}}\def\labelenumi{(\theenumi)}
\item for any place~$v$ of~$K$, the formal curve $\widehat C_{K_v}$
is $K_v$-analytic;
\item the infinite product $\prod_{\mathfrak p\nmid N}
 S_{\mathscr X_{R_{\mathfrak p}}}(\widehat C_{K_{\mathfrak p}})$
converges to a positive real number.
\end{enumerate}
\end{definition}

Condition (ii) asserts precisely that the series with non-negative
terms
\[ \sum_{\mathfrak p\nmid N}
 \log S_{\mathscr X_{R_{\mathfrak p}}}(\widehat C_{K_{\mathfrak
p}})^{-1}\]
 is convergent.

Observe that the above definition does not depend
on the choices required to formulate it.
Indeed, condition~(i) does not involve any choice.
Moreover, if condition~(i) holds and if $N'$ is any positive multiple
of~$N$, condition~(ii) holds for 
 $(N,\mathscr X,\mathscr P)$ if and only if it
holds for  $(N',\mathscr X_{R[1/N']},\mathscr P_{R[1/N']})$.
Moreover, for any two such triples $(N_1,\mathscr X_1,\mathscr P_1)$
and $(N_2,\mathscr X_2,\mathscr
P_2)$, there is a positive integer~$M$, multiple of both $N_1$ and
$N_2$, 
such that the models $(\mathscr X_1,\mathscr P_1)$ and $(\mathscr
X_2,\mathscr P_2)$ of~$(X,P)$
become isomorphic over $R[1/M]$. This shows that, when (i) is
satisfied,  conditions~(ii)
for any two triples $(N,\mathscr X,\mathscr P)$ are indeed equivalent.

It  follows from the properties of the size recalled
in Proposition~\ref{theo.def-size}
that  $A$-analyticity is  invariant under  immersions
and  compatible to \'etale localization.

As a consequence of Proposition~\ref{prop.size-analytic} and
\ref{prop.size-smooth}, we also have:

\begin{proposition}\label{preEisenstein}
Let $\widehat C$ be a smooth formal curve which is $K_v$-analytic
for any place~$v$ of~$K$. Assume that 
$\widehat C$
extends to a smooth formal curve $\widehat {\mathscr C}\hra\mathscr X$
over $R[1/N]$, for some $N\geq 1$.
Then $\widehat C$ is~$A$-analytic.
\end{proposition}
%\begin{proof}
\removelastskip
Indeed, these conditions imply that the size of~$\widehat C$
at almost every finite place of~$K$ is equal to~$1$,
while being positive at every place.
%\end{proof}

\medskip

As observed in substance by Eisenstein \cite{eisenstein1852},
any algebraic smooth formal
curve satisfies the hypothesis of Proposition~\ref{preEisenstein}.
Therefore: 

\begin{corollary}%[Eisenstein]
If the smooth formal curve $\widehat C$  is algebraic,
then it is~$A$-analytic.
\end{corollary}

The invariance of size under extensions of valued fields established in Proposition~\ref{prop.size-extensions} easily
implies that, \emph{for any number field $K'$ containing~$K$, the
smooth formal curve $\widehat{C}':=\widehat{C}_{K'}$ in~$X_{K'}$
deduced from $\widehat C$ by the extension of scalars $K
\hookrightarrow K'$ is~$A$-analytic iff $\widehat C$ is~$A$-analytic.}

Let 
$\phi\in K[[X]]$ be any formal power series, and let $P:=(0,\phi(0))$.
From the inequality in Proposition~\ref{prop.size-radius-newton}, 1),  between the 
convergence radius of a power series and the size of its graph, it follows that
\emph{the~$A$-analyticity of the graph~$\widehat C$ of~$\phi$ in~$\widehat{\A^2_P}$
implies that the convergence
radii $R_v$ of~$\phi$ at the places~$v$ of~$K$  satisfy the so-called Bombieri's condition
 \[ \prod_v \min(1,R_v)>0, \]
or equivalently
\[\sum_v \log^+ R_v^{-1}< + \infty.\]}
 %%%%%%%%%%%%%%%%%%%%%%%%%%%%%%%%%%%%%%%%%%%%%%%%%%%%%%%%%%%%%%%%%%%%%%%%%%%
%                                                                         %
%    , the so-called Bombieri condition.                                  %
%   Indeed, for almost all places~$v$,                                    %
%   $\min(1,R_v)$ is at least equal to the size of~$\widehat C$ at~$v$.   %
%                                                                         %
%%%%%%%%%%%%%%%%%%%%%%%%%%%%%%%%%%%%%%%%%%%%%%%%%%%%%%%%%%%%%%%%%%%%%%%%%%%
However the converse does not hold, as can be seen by considering
the power series
$\phi(X)=\log(1+X)$, that satisfies that Bombieri's condition (since all the
$R_{v}$ equal 1), but is not~$A$-analytic (its~$p$-adic size is
$\abs{p}^{1/(p-1)}$ and 
 the infinite series $\sum \frac{1}{p-1}\log p$ diverges).

Let us  conclude this section by a brief discussion
of the relevance of~$A$-analyticity in the arithmetic theory of
differential equations (we refer to \cite{bost2001,chambert-loir2002, bost2004} for more details).
%case of differential equations. % and Grothendieck-Katz conjecture.

Assume that~$X$ is smooth over~$K$, that~$F$ is sub-vector bundle
of rank one in the tangent bundle $T_{X}$ (defined over~$K$), and that
$\widehat C$ is the formal leaf
at~$P$ of
the one-dimensional algebraic foliation on~$X$ defined by~$F$.
By a model of~$(X,F)$ over~$R[1/N]$, we mean the  data of 
a scheme $\mathscr X$ quasi-projective and  smooth over $\Spec R$, of  a coherent subsheaf
$\mathscr F$ 
of~$T_{\mathscr X/R}$, and of an isomorphism $X\simeq \mathscr X\otimes K$
inducing an isomorphism $F\simeq \mathscr F\otimes K$. Such models clearly
exist if~$N$ is sufficiently divisible. Let us choose one of them
 $(\mathscr X,\mathscr F)$.
We say that the foliation~$F$ satisfies the \emph{Grothendieck-Katz
condition}
if for almost every maximal ideal $\mathfrak p\subset R$,
the subsheaf $\mathscr F_{\F_{\mathfrak
p}}$ of~$T_{\mathscr X_{\F_{\mathfrak
p}/\F_{\mathfrak
p}}}$ is closed under~$p$-th
powers, where~$p$ denotes the characteristic of~$\F_{\mathfrak
p}$. 
As above, this condition does not depend on the  choice
of the model 
$(\mathscr X,\mathscr F)$. 

\begin{proposition} With the above notation, if  $F$ satisfies the Grothendieck-Katz
condition, then its formal integral curve
 $\widehat C$ through any rational point $P$ in $X(K)$ is~$A$-analytic.
\end{proposition}
\begin{proof}
It follows from Cauchy's theory of  analytic ordinary
differential equations over local fields that
the formal curve $\widehat C$ is $K_v$-analytic for any 
place~$v$ of~$K$.

After possibly increasing~$N$, we may assume that~$P$ extends to a
section $\mathscr P$ in~$\mathscr X(R[1/N])$.
For any maximal ideal $\mathfrak p\subset R$ that
is unramified over a prime number~$p$,
and such that $\mathscr F_{\F_{\mathfrak p}}$ is closed under~$p$-th power,
Proposition~\ref{prop.size-foliation} show that
the $\mathfrak p$-adic size of~$\widehat C$ is at least
$\abs{p}^{1/p(p-1)}$.  
When  $F$ satisfies the Grothendieck-Katz
condition, this inequality holds for almost
all maximal ideals of~$R$.
Since the series over primes $\sum_p \frac1{p(p-1)}\log p$
converges, this implies the convergence of 
the series
$\sum_{\mathfrak p\nmid N}
 \log S_{\mathscr X_{R_{\mathfrak p}}}(\widehat C_{K_{\mathfrak
p}})^{-1}$
and consequently the~$A$-analyticity of~$\widehat C$.
\end{proof}

\section{Analytic curves in algebraic varieties over local fields and
canonical semi-norms}
\label{sec.can}

\subsection{Consistent sequences of norms}\label{subsec:consistent}

Let~$K$ be a local field, $X$ a projective scheme over~$K$, and~$L$ 
a line bundle over~$X$.

We may consider the following natural constructions of sequences of
norms on the spaces of sections $\Gamma(X, L^{\otimes n})$:
\begin{enumerate}
\item When $K=\C$ and~$X$ is reduced, we may choose an arbitrary
continuous norm
$\left\|.\right\|_{L}$ on the $\C$-analytic line bundle~$L_{\an}$
defined by~$L$ on the compact and reduced complex analytic space
$X(\C)$. Then, for any integer~$n$, the space of algebraic regular
sections $\Gamma(X, L^{\otimes n})$ may be identified with a subspace
of the space of continuous sections of~$L^{\otimes n}_{\an}$ over
$X(\C)$. It may therefore be equipped with the restriction of the
 ${\rm L}^{\infty}$-norm, defined by:
 \begin{equation}\label{supnorm}
  \left\| s\right\|_{{\rm L}^{\infty},n}:=\sup_{x \in
  X(\C)}\left\| s(x)\right\|_{L^{\otimes n}} \,\,\mbox{ for any
  }s\in \Gamma(X, L^{\otimes n}),  
 \end{equation}
where $\left\| .\right\|_{L^{\otimes n}}$ denotes the continuous norm
on~$L_{\an}^{\otimes n}$ deduced from $\left\| .\right\|_{L}$ by
taking the~$n$-th
tensor power.

This construction admits a variant where, instead of the sup-norms
(\ref{supnorm}), one considers the
${\rm L}^p$-norms defined by using some ``Lebesgue measure''
(\cf \cite[4.1.3]{bost2001}, and 
\cite[Th\'eor\`eme~3.10]{randriam2006}).
%This construction admits the following variant. Let us assume for
%simplicity that~$X$ is of pure dimension~$d$, and let $\mu$ be a
%positive Lebesgue measure\footnote{Namely, a Borel measure on
%$X(\C)$ 
%such that, for any complex analytic embedding~$j$ of some open
%subset 
%$U$ of~$X(\C)$ in some numerical space $\C^N$, the direct image
%measure $j_{\star}\mu$ and the Lebesgue measure $\lambda_{j(U)}$ ---
%defined by the current $\frac{1}{d!}(\frac{1}{2i} \sum_{1\leq k \leq
%N}
%dz_{k}\wedge d\overline{z}_{k})^d \delta_{j(U)}$ --- on~$j(U)$ are
%comparable, in the sense that $j_{\star}\mu =e^{\phi}
%\lambda_{j(U)}$ 
%for some locally bounded Borel real function $\phi$ on~$j(U)$.
%}
%on~$X(\C)$.
%Then, instead of the sup-norm (\ref{supnorm}), we may consider the
%${\rm L}^2$-norm $\left\|.\right\|_{{\rm L}^2,n}$ defined by 
%\begin{equation}\label{L2}
%  \left\|s\right\|_{{\rm L}^2,n}^2:=\int_{X(\C)}
%  \left\|s(x)\right\|_{L^{\otimes n}} \, d\mu(x).
% \end{equation}

\item When $K=\R$ and~$X$ is reduced, we may choose
a continuous norm on~$L_\C$ that is invariant under
complex conjugation. The previous
constructions define
complex norms on the complex vector spaces
 \[ \Gamma(X,L^{\otimes n})\otimes_{\R}\C \simeq
\Gamma(X_{\C},L_{\C}^{\otimes
 n})\]
 which are invariant under complex conjugation, and  by 
 restriction, real norms on the real vector spaces
$\Gamma(X,L^{\otimes n})$. 

\item When~$K$ is a~$p$-adic field, with ring of integers
$\cO$, we may
choose a pair $(\cX, \cL)$, where $\cX$ is a projective flat model of 
$X$ over $\cO$, and $\cL$ a line bundle over $\cX$ extending~$L$.
Then, for any integer~$n$, the $\cO$-module $\Gamma(\cX, \cL^{\otimes n})$ 
is 
%% (torsion-)
free of finite rank and may be identified with an
$\cO$-lattice in the~$K$-vector space $\Gamma(X, L^{\otimes n})$, and 
consequently defines a norm on the latter 
% (\cite[prop.~6, p.~28]{weil1967})
--- namely, the norm
$\left\|\cdot\right\|_{n}$ such that a 
section $s\in \Gamma(X, L^{\otimes n})$
satisfies $\left\|s\right\|_{n}\leq 1$ iff~$s$ extends to a section of
$\cL^{\otimes n}$ over $\cX$.

\item 
A variant of Construction~(1) can be used when~$K$ is a~$p$-adic field and $X$ is reduced.
Let~$\norm{\cdot}$ be a metric on~$L$ (see Appendix~\ref{app.metrics}
for basic definitions concerning metrics in the~$p$-adic setting).
For any integer~$n$, the space $\Gamma(X,L^{\otimes n})$
admits a $\mathrm L^\infty$-norm, defined for any $s\in\Gamma(X,L^{\otimes n})$
by 
$ \norm{s}_{\mathrm L^\infty,n} \mathrm{:=}
   \sup_{x\in X(C)} \norm{s(x)}$, where $C$ denotes the completion of an algebraic closure of $K$.
When the metric of~$L$ is defined by a model~$\mathscr L$
of~$L$ on a \emph{normal} projective model~$\mathscr X$ of~$X$ on~$R$,
then this norm coincides with that defined by construction~(3)
(see, \eg, \cite[Proposition~1.2]{rumely-l-v2000}).
 \end{enumerate}
 \smallskip

For any given~$K$, $X$, and~$L$ as above, we
shall say that two sequences $(\left\|.\right\|_{n})_{n \in \N}$
and $(\left\|.\right\|'_{n})_{n \in \N}$ of norms on the finite
dimensional~$K$-vector spaces $(\Gamma(X, L^{\otimes n}))_{n \in \N}$
are \emph{equivalent} when, for some positive constant~$C$ and any
positive integer~$n$,
\[ C^{-n} \left\|.\right\|_{n}' \leq \left\|.\right\|_{n}
\leq C^n \left\|.\right\|_{n}'.\]

One easily checks that, for any given~$K$, $X$ and~$L$, the above
 constructions provide sequences of
norms  $(\left\|.\right\|_{n})_{n \in \N}$ on the sequence of spaces 
$(\Gamma(X, L^{\otimes n}))_{n \in \N}$ that are all equivalent. In particular, their equivalence class does not depend on
 the auxiliary data (models, norms on~$L$, \dots) involved.
(For the comparison of the ${\rm L}^2$ and ${\rm L}^\infty$
norms in the archimedean case, 
see notably~\cite[Th\'eor\`eme~3.10]{randriam2006}.) 

A sequence of norms on the spaces $\Gamma(X, L^{\otimes n})$
that is equivalent to one 
(or, equivalently, to any) 
of the sequences thus constructed will be called
\emph{consistent}. This notion immediately extends to sequences
$(\left\|.\right\|_{n})_{n\geq n_{0}}$ of norms on the spaces
$\Gamma(X, L^{\otimes n})$, defined only for~$n$ large enough.

When the line bundle~$L$ is ample, consistent sequences of
norms are also provided by additional constructions. Indeed we have:
\begin{proposition}\label{ampleconsistent}
    Let~$K$ be a local field, $X$ a projective scheme over~$K$, and
    $L$ an \emph{ample} line bundle over~$X$. Let moreover~$Y$ be a
    closed subscheme of~$X$, and assume~$X$ and~$Y$ reduced when~$K$
    is archimedean.
    
    For any consistent sequence of norms $(\left\|.\right\|_{n})_{n
\in \N}$ on 
$(\Gamma(X, L^{\otimes n}))_{n \in \N}$, the quotient norms 
$(\left\|.\right\|'_{n})_{n\geq n_{0}}$ on the spaces 
$(\Gamma(Y, L_{\mid Y}^{\otimes n}))_{n\geq n_{0}}$, deduced from the 
norms $\left\|.\right\|_{n}$ by means of the restriction maps
$\Gamma(X, L^{\otimes n}) \longrightarrow \Gamma(Y, L_{\mid
Y}^{\otimes n})$
--- \emph{which are surjective for % $n$
$n\geq n_0$ 
large enough since $L$ is ample} --- constitute a consistent sequence.
\end{proposition}

When~$K$ is archimedean, this is proved in~\cite[Appendix]{bost2004},
by
introducing a positive metric on~$L$, as 
a consequence of Grauert's finiteness theorem for pseudo-convex
domains applied to the unit disk bundle of~${L}^\vee$ (see
also~\cite{randriam2006}). 

When~$K$ is a 
$p$-adic field with ring of integers $\mathcal O$, Proposition
\ref{ampleconsistent} follows from the
basic properties of ample line bundles over projective $\mathcal
O$-schemes.
Let indeed $\cX$ be a projective flat model of~$X$ over~$\cO$, 
$\cL$ an ample line bundle on~$\cX$,  $\cY$ the closure of~$Y$ in~$\cX$, 
and $\mathcal{I_Y}$ the ideal sheaf of~$\cY$. If the 
positive integer~$n$ is large enough, then the cohomology group 
$H^1(\cY,\mathcal{I_Y}\cdot \cL^{\otimes n})$ vanishes, and 
the restriction morphism  
$\Gamma(\cX, \cL^{\otimes n}) \ra \Gamma(\cY, \cL_{\mid
\cY}^{\otimes n})$ is therefore surjective. Consequently, the norm on 
$\Gamma(Y, L_{\mid Y}^{\otimes n})$ attached to the lattice
$\Gamma(\cY, \cL_{\mid \cY}^{\otimes n})$ is the quotient of the norm 
on 
$\Gamma(X, L^{\otimes n})$ attached to $\Gamma(\cX, \cL^{\otimes n})$.

Let~$E$ be a finite dimensional vector space over the local field~$K$,
equipped with some norm, supposed to be euclidean or hermitian in the 
archimedean case. This norm induces similar norms on the tensor powers
$E^{\otimes n}$, $n\in \N$, hence --- by taking the quotient norms ---
on the symmetric powers $\Sym^n E$.
If~$X$ is the projective space $\P(E):=\operatorname{\mathbf {Proj}}
\Sym^\cdot (E)$ and~$L$ the line bundle $\cO (1)$ over $\P(E)$, then the canonical
isomorphisms $\Sym^n E \simeq \Gamma (X, L^{\otimes n})$ allow one to
see these norms as a sequence of norms on 
$(\Gamma(X, L^{\otimes n}))_{n \in \N}$.
One easily checks that this sequence is consistent.  (This is 
straightforward in the~$p$-adic case. When~$K$ is archimedean, this
follows for instance from~\cite[Lemma~4.3.6]{bost-g-s94}.) 

For any closed subvariety~$Y$ in~$\P(E)$ and any $n \in \N$, 
we may consider the commutative diagram of~$K$-linear maps:
%$$\begin{array}{ccccc}
%    S^n E & \stackrel{\sim}{\longrightarrow}  & S^n \Gamma ({\mathbb
%P}(E), \cO(1)) &\stackrel{\sim}{\longrightarrow}   & 
%    \Gamma ({\mathbb P}(E), \cO(n))
%  \\
%     &  & \Big \downarrow &  & \Big \downarrow \\
%     &  & S^n \Gamma (Y, \cO(1)) &
%\stackrel{\beta_{n}}{\longrightarrow}  & \Gamma (Y,
%     \cO(n))
%\end{array}$$
%%  
%%  \begin{array}{ccccc}
%%      S^n E & \stackrel{\sim}{\longrightarrow}  & S^n \Gamma ({\mathbb
%%  P}(E), \cO(1)) &\stackrel{\sim}{\longrightarrow}   & 
%%      \Gamma ({\mathbb P}(E), \cO(n))
%%    \\
%%       &  & \Big \downarrow &  & \Big \downarrow\rlap{$\alpha_n$} \\
%%       &  & S^n \Gamma (Y, \cO(1)) &
%%  \stackrel{\beta_{n}}{\longrightarrow}  & \Gamma (Y,
%%       \cO(n))
%%  \end{array}\]
\[
\xymatrix{
   {\Sym^n E} \ar[r]^-{\sim}
& {\Sym ^n \Gamma(\P(E),\cO(1))} \ar[d] \ar[r]^-{\sim} 
& \Gamma(\P(E),\cO(n))  \ar[d]^{\alpha_n}  \\
& {\Sym^n \Gamma(Y,\cO(1))} \ar[r]^-{\beta_n} 
& \Gamma(Y,\cO(n)) } \]
where the vertical maps are the obvious restriction morphisms. The
maps $\alpha_{n}$, and consequently $\beta_{n}$, are surjective if
$n$ is large
enough.

Together with Proposition~\ref{ampleconsistent}, these observations
yield the following corollary:

\begin{corollary}
    Let~$K$, $E$ and~$Y$ a closed subscheme of~$\P(E)$ be as above.
Assume that~$Y$ is reduced if~$K$ is archimedean. 
Let us choose a
norm 
   on~$E$ (resp. on~$\Gamma(Y, \cO(1))$) and let us equip $\Sym^n E$
(resp.
   $\Sym^n \Gamma(Y, \cO(1))$) with the induced norm, for any $n\in \N$.
   
   Then the sequence of quotient norms on~$\Gamma(Y, \cO(n))$ defined
   for~$n$ large enough by means of the surjective morphisms
   $\alpha_{n}\colon \Sym^n E \rightarrow \Gamma(Y, \cO(n))$ (resp. by
   means of
   $\beta_{n}\colon \Sym^n \Gamma (Y, \cO(1)) \rightarrow \Gamma(Y,
   \cO(n))$) is consistent. 
\end{corollary}

\subsection{Canonical semi-norms}
\label{subsec.canonical}

Let~$K$ be a local field.
Let~$X$ be a \emph{projective} variety over~$K$, $P$ a rational point in~$X(K)$,
and $\widehat C$ be a smooth~$K$-\emph{analytic} formal curve
in~$\widehat X_P$.
To these data, we are going to attach a canonical semi-norm
 $\|\cdot\|^{\can}_{X, \hat{C}}$
on the tangent line $T_P\widehat C$ of~$\widehat C$ at~$P$. It will be defined by
considering an avatar of the evaluation map 
\[ 
 E^i_D/E^{i+1}_D \hra \Gamma(\mathscr P, \mathscr O(D)\otimes
\mathscr N^{\vee\otimes i}) \]
which played a prominent role in our
proof of Proposition~\ref{prop.geom-algebraicity}.

The construction of~$\|\cdot\|^{\can}_{X, \hat{C}}$ will require auxiliary data, 
of which it will
eventually not depend.

Let us choose a line bundle~$L$ on~$X$ and a consistent
sequence of norms on the 
$K$-vector spaces~$E_D=\Gamma(X,L^{\otimes D})$, for $D\in\N$.
% and a continuous
% metric on~$L$ (see Appendix~\ref{app.metrics} for definitions).
% Let us denote $E_D=\Gamma(X,L^{\otimes D})$; the line
% bundle~$L^{\otimes
% D}$ being endowed with the tensor product metric, the formula
% \[ \norm{s} = \sup_{x\in X(\bar K)} \norm{s(x)}, \]
% for $s\in E_D$, defines a norm on~$E_D$.
Let us also fix norms~$\norm{\cdot}_0$ 
on the~$K$-lines~$T_P\widehat C$ and $L_{\mid P}$.

Let us denote by~$C_i$ the~$i$th neighborhood of~$P$
in~$\widehat C$. Thus we have $C_{-1}=\emptyset$,
$C_0=\{P\}$, and $C_i$ is a~$K$-scheme isomorphic to $\Spec K[t]/(t^{i+1})$;
moreover, $\widehat C=\varinjlim C_i$.
Let us denote by $E^i_D$ the~$K$-vector subspace  of the $s\in E_D$
such that $s_{\mid C_{i-1}}=0$.
The restriction map $E_D\ra\Gamma(C_i,L^{\otimes D})$
induces a linear map of finite dimensional~$K$-vector spaces
\[ \phi^i_D \colon E^i_D \ra \Gamma(C_{i},
\mathcal{I}_{C_{i-1}}\otimes L^{\otimes D})\simeq (T^\vee_P\widehat C)^{\otimes i}
\otimes L_{\mid P}^{\otimes D}. \]

We may consider the $\norm{\phi^i_D}$ of this map, computed by using
the chosen norms on~$E_{D}$, $T_P\widehat C$, and $L_{\mid P}$, and
the ones they induce by restiction, duality and tensor product on
$E^i_D$ and on~$(T^\vee_P\widehat C)^{\otimes i}
\otimes L_{\mid P}^{\otimes D}$.

Let us now define $\rho(L)$ by the following formula:
\[ \rho (L) = \limsup_{i/D\ra +\infty} \frac1i \log \norm{\phi^i_D}
.\]

The analyticity of~$\widehat C$ implies that  $\rho(L)$ belongs 
to $[-\infty,+\infty\mathclose[$. Indeed,
when~$K$ is $\C$ or $\R$,  as observed in \cite[\S3.1]{bost2004},
from Cauchy inequality we easily derive the existence of positive real numbers $r$ and $C$ such that
\begin{equation}
\label{archCauchy}
\norm{\phi^i_D} \leq C^{D+1} r^{-i}.
\end{equation}
 When~$K$ is ultrametric, we may actually bound $\rho(L)$ in terms of
the size of~$\widehat C$:

\begin{lemma}\label{lemm.can/size}
Assume that~$K$ is ultrametric and let~$R$ be its ring of integers.
Let $\mathscr X$ be a projective flat~$R$-model of~$X$
and let $\mathscr P\colon \Spec R\ra\mathscr X$  the section
extending~$P$. Assume moreover that the metric of~$L$
is given by a line bundle~$\mathscr L$ on~$\mathscr X$ extending~$L$
and the consistent sequence of norms on~$(E_D)$ by the construction~(3)
in  Section~\ref{subsec:consistent},
and fix the norm~$\norm{\cdot}_0$ on~$T_P\widehat C$ 
so that its unit ball is equal to $N_{\mathscr P}\mathscr X
\cap T_P\widehat C$.

Then, one has 
\[ \rho(L)\leq -\log S_{\mathscr X,\mathscr P}(\widehat C).\]
\end{lemma}
\begin{proof} 
Let~$r$ be an element of~$\mathopen]0,  S_{\mathscr X}(\widehat C)\mathclose[$. 
We claim that, with the notation above, we have:
\[ \norm{\phi_i^D} \leq r^{-i}.\] 
This will establish that 
$\rho(L)=\limsup_{i/D \ra + \infty}\frac1i \log \norm{\phi_i^D} \leq - \log r$, 
hence the required inequality
by letting~$r$ go to~$ S_{\mathscr X}(\widehat C)$.

To establish the above estimate on~$\norm{\phi_i^D}$, let us choose
an affine open neighbourhood~$U$ of~$\mathscr P$ in~$\mathscr X$
such that~$\mathscr L_{U}$ admits a non-vanishing section~$l$, and a
closed embedding $i\colon U \hookrightarrow \A^N_{R}$ such that 
$i(\mathscr P)=(0,\ldots,0)$.
Let $\widehat{C}'$ denote the image of~$\widehat{C}$
by the embedding of formal schemes $\widehat{i_K}_P: \hat{X}_P \hookrightarrow \widehat{\A^N_{K,0}}$. 
%%%%%%%%%%%%%%%%%%%%%%%%%%%%%%%%%%%%%%%%%%%%%%%%%%%%%%%%%%%%%%%%%%%%%%%%%%%
%                                                                         %
%   The invariance of the size under closed embeddings and                %
%   \'etale localization allows us                                        %
%   to assume that $\mathscr X$ is the affine space $\A^N_R$,             %
%   that $P=(0,\dots,0)$ and that $\mathscr L=\mathscr O_{\mathscr X}$.   %
%                                                                         %
%%%%%%%%%%%%%%%%%%%%%%%%%%%%%%%%%%%%%%%%%%%%%%%%%%%%%%%%%%%%%%%%%%%%%%%%%%%
% Consider $r<S_{\mathscr X}(\widehat C)$ and a formal automorphism
By the very definition of the size, we may find
$\Phi$ in~$G_{{\rm an},r}$ such that 
%   \in \Aut_0(\B_r^N)$ with
$\Phi^*\widehat{C}'=\widehat{\A^1_0}\times\{0\}^{N-1}$.
Let~$s$ be an element of
$\Gamma(\mathscr X,\mathscr L^{\otimes D})$. 
% vanishing at $\mathscr P$. 
We may write
$s_{\mid U}=i^*Q\cdot l^{\otimes D}$
for some~$Q$ in~$R[X_1,\ldots,X_N]$. % such that $P(0,\ldots,0)=0$. 
Then, $\Phi^* Q$ is given by a
formal series $g=\sum b_I X^I$ which satisfies $\norm{g}_r\leq 1$, 
or equivalently, $\abs{b_I}r^{\abs I}\leq 1$
for any multiindex~$I$.
If~$s$ belongs to $E^i_D$,
with the chosen normalizations of norms, we have: $\norm{\phi_i^D(s)}
= \abs{b_{i,0,\dots,0}}\leq r^{-i}$. 
\end{proof}

The exponential of~$\rho(L)$ is a well defined element in~$[0,+\infty\mathclose[$,
and we may introduce the following

\begin{definition} 
The \emph{canonical semi-norm} on~$T_P\widehat C$ attached 
to~$(X,\widehat C,L)$ is 
\[ \norm{\cdot}^{\can}_{X,\widehat C,L}:=e^{\rho(L)} \norm{\cdot}_0.\]
\end{definition}

Observe that, if $\widehat C$ is algebraic, then there exists a real number $\lambda$ such that the filtration $(E^i_D)_{i\in \N}$ becomes stationary --- or equivalently $\phi^i_D$ vanishes --- for $i/D > \lambda$ (for instance, we may take the degree of the Zariski closure of $\widehat C$ for $\lambda$). Consequently, in this case, $\rho(L)= - \infty$ and the canonical semi-norm $\norm{\cdot}^{\can}_{X,\widehat C,L}$ vanishes.

The notation $\norm{\cdot}^{\can}_{X,\widehat C,L}$ for the canonical semi-norm --- which makes reference to~$X$, $\widehat C$, and~$L$ only --- 
is justified by the first part in the next Proposition:

\begin{proposition}
\begin{enumerate}\def\theenumi{\alph{enumi}}\def\labelenumi{\theenumi)}
\item 
The semi-norm  $\norm{\cdot}^{\can}_{X,\widehat C,L}$ is independent of the choices
of  norms on~$T_P\widehat C$ and $L_{\mid P}$, and of the consistent sequence
of norms on the spaces~$E_D:=\Gamma(X,L^{\otimes D})$.

\item
For any positive integer~$k$, the semi-norm  $\norm{\cdot}^{\can}_{X,\widehat C,L}$ is unchanged if~$L$ is replaced by $L^{\otimes k}$.

\item
Let $L_1$ and~$L_2$ be two line bundles
and assume that $L_2\otimes L_1^{-1}$ has a regular section~$\sigma$
over~$X$ 
that does not vanish at~$P$.
Then
\[ \norm{\cdot}^{\can}_{X,\widehat C,L_1} \leq
\norm{\cdot}^{\can}_{X,\widehat C,L_2}. \]
\end{enumerate}
\end{proposition}
\begin{proof}
\emph a)
Let us denote with primes another family of norms on
the spaces~$T_P\widehat C$, $L_{\mid P}$, and $E_D$, and by $\rho'(L)$ and $(\norm{\cdot}^{\can}_{X,L,\widehat C})' $ the attached ``rho-invariant" and canonical semi-norm.
There are positive real numbers~$a$, $b$, $c$ such that
$ \norm{t}'_0= a \norm{t}_0 $ for any $t\in T_P\widehat C$,
$\norm{s(P)}'=b\norm{s(P)}$ for any local section~$s$ of~$L$ at~$P$,
and 
\[ c^{-D} \norm{s} \leq \norm{s}' \leq c^D \norm{s}  \]
for any positive integer~$D$ and any global section~$s\in E_D$.
Consequently, for $(i,D)\in\N^2$ and $s\in E^i_D$, 
\[ \norm{\phi^i_D(s)}' = a^{-i}b^D \norm{\phi^i_D(s)}
\leq a^{-i}b^D \norm{\phi^i_D} \norm{s} 
\leq a^{-i}b^D \norm{\phi^i_D} c^D \norm{s}', \]
hence
\[  \norm{\phi^i_D}' \leq  a^{-i}c^D b^D \norm{\phi^i_D} \]
and
\[  \frac1i \log\norm{\phi^i_D}' \leq - \log a +  \frac Di \log (bc)
 + \frac1i \log\norm{\phi^i_D}. \]
When $i/D$ goes to infinity, this implies:
\[ \rho'(L) \leq -\log  a + \rho(L), \]
from which follows:
\[  (\norm{\cdot}^{\can}_{X,L,\widehat C})' \leq
\norm{\cdot}^{\can}_{X,L,\widehat C}, \]
by definition of the canonical semi-norm.
The opposite inequality also holds by symmetry,
hence the desired equality.

\emph b) To define $\rho(L)$ and $\rho(L^{\otimes k})$, let us use the same norm $\norm{\cdot}_0$ 
on~$T_P\widehat C$, and assume that the consistent sequence of norms chosen on $(\Gamma(X,L^{\otimes D})$ is defined by one of the constructions (1-4) in the above subsection~\ref{subsec:consistent}, and finally that the one on   
$(\Gamma(X,(L^{\otimes k})^{\otimes D}))=(\Gamma(X,L^{\otimes kD}))$ is extracted from the one on $(\Gamma(X,L^{\otimes D})$.

Specifying the line bundle with a supplementary index, one has
\[ \phi^i_{D,L^{\otimes k}}=\phi^i_{kD,L}.\]
The definition of an upper limit therefore
implies that $\rho(L^k)\leq\rho(L)$. 

To establish the opposite inequality, observe that, for any section $s$ in $E^i_{D,L}$ and any positive integer $k$, the $k$-th tensor power $s^{\otimes k}$ belongs to $E^{ki}_{D,L^{\otimes k}}$ and 
\[ \phi^{ki}_{D, L^{\otimes k}}(s^{\otimes k})=(\phi^i_{D,L}(s))^{\otimes k}.\]
 Let $\rho$ be
any real number such that $\rho<\rho(L)$,
and choose~$i$, $D$, and $s\in E^i_{D,L}$ such that
$\norm{\phi^i_{D,L}(s)}\geq e^{\rho i} \norm{s}$.
Then, for any positive integer $k$, we have
 \[ \norm{\phi^{ki}_{D, L^{\otimes k}}(s^{\otimes k})}=\norm{\phi^i_{D,L}(s)}^k
\geq e^{\rho ki} \norm{s}^k = e^{\rho ki} \norm{s^{\otimes k}}, \]
so that $\norm{\phi^{ki}_{D, L^{\otimes k}}}^{1/ki}\geq e^\rho$. Consequently,
$\rho(L^k)\geq \rho$.

\emph c) Here again, we may use
 the same norm $\norm{\cdot}_0$ 
on~$T_P\widehat C$ to define $\rho(L_1)$ and $\rho(L_2)$, and 
 assume that the consistent sequence of norms chosen on $(\Gamma(X,L_1^{\otimes D})$ and
  $(\Gamma(X,L_2^{\otimes D})$ are defined by one of the constructions (1-4) above.

If~$s$ is a global section of~$L_1^{\otimes D}$,
then $s\otimes\sigma^{\otimes D}$ is a global section of~$L_2^{\otimes D}$;
if~$s$ vanishes at order~$i$ along~$\widehat C$, so does
$s\otimes\sigma^{\otimes D}$ and
\[ \phi^i_{D,L_2}(s\otimes\sigma^{\otimes D}) =  \phi^i_{D,L_1}(s)\otimes \sigma(P)^{\otimes D}.
\]
Consequently, 
\[  \norm{\phi^i_{D,L_1}(s)} \leq \norm{\phi^i_{D,L_2}}.
 \norm{s\otimes\sigma^{\otimes D}}.\norm{\sigma(P)}^{-D}  \leq
(\norm{\sigma(P)}^{-1}\norm{\sigma})^D.
\norm{\phi^i_{D,L_2}} .\norm{s}, \]
and $\rho(L_1)\leq \rho(L_2)$, as was to be shown.
\end{proof}

\begin{corollary} 
The set of semi-norms on~$T_P\widehat C$ described 
by~$\norm{\cdot}^\can_{X,\widehat C,L}$ when~$L$ varies in the class of
line bundles on~$X$ possesses a maximal element, namely the canonical
semi-norm $\norm{\cdot}^\can_{X,\widehat C,L}$ attached to any ample
line bundle~$L$ on~$X$. 
\end{corollary}

We shall denote~$\norm{\cdot}^\can_{X,\widehat C}$ this maximal element.
The formation of~$\norm{\cdot}^\can_{X,\widehat C}$ satisfies the
following compatibility properties with respect to rational morphisms.

\begin{proposition} 
Let $X'$ be another projective algebraic variety over~$K$,
and let $f\colon X\dashrightarrow X'$ be a rational
map that is defined near~$P$. Let $P':=f(P)$, and assume that~$f$
defines an (analytic, or equivalently, formal) isomorphism from~$\widehat C$
onto a smooth~$K$-analytic formal curve~$\widehat{C'}$ in~$\widehat{X'}_{P'}$.

Then for any $v\in T_P\widehat C$, 
\[ \norm{Df(P)v}_{X',f(\widehat C)}\leq \norm{v}_{X,\widehat C}.\]

If moreover~$f$ is an immersion in a neighborhood of~$P$, then the
equality holds. 
\end{proposition}

When~$K$ is archimedean, this summarizes the results established
in~\cite[Sections~3.2 and~3.3]{bost2004}. The arguments in \emph{loc. cit.}
may be immediately transposed to the ultrametric case, by using consistent
norms as defined above instead of~$\mathrm{L}^\infty$ norms on the spaces
of sections~$E_D$. We leave the details to the reader.

Observe finally that this Proposition allows us to define the canonical
semi-norm~$\norm{\cdot}^\can_{X,\widehat C}$ when the algebraic
variety~$X$ over~$K$ is only supposed to be quasi-projective. Indeed,
if $\overline{X}$ denotes some projective variety containing~$X$
as an open subvariety,  the semi-norm~$\norm{\cdot}^\can_{\overline
X,\widehat C}$ is independent of the choice of~$\overline X$, and we
let 
\[ \norm{\cdot}^\can_{X,\widehat C}:=\norm{\cdot}^\can_{\overline
X,\widehat C}.\]

\section{Capacitary metrics on~$p$-adic curves}
\label{sec.capa}

\subsection{Review of the complex case}

Let~$M$ be a compact Riemann surface and let~$\Omega$
be an open subset of~$M$. 
We assume that the compact subset complementary to~$\Omega$
in any connected component of~$M$ is not polar.
Let~$D$ be an effective divisor on~$M$ whose support is contained in~$\Omega$. 
Potential theory on Riemann surfaces (see~\cite[3.1.3--4]{bost1999})
shows the existence of a unique subharmonic function $g_{D,\Omega}$ on~$M$
satisfying the following assumptions:
\begin{enumerate}
\item $g_{D,\Omega}$ is harmonic on~$\Omega\setminus\abs D$;
\item the set of points~$z\in M\setminus\Omega$ such that $g_{D,\Omega}(z)\neq0$
is a polar subset of~$\partial\Omega$;
\item for any open subset~$V$ of~$\Omega$ and any holomorphic
function~$f$ on~$V$ such that $\div(f)=D$, 
the function $g_{D,\Omega}-\log\abs f^{-2}$ on $V\setminus \abs D$ is the restriction of a harmonic function
 on~$V$.
\end{enumerate}
Moreover, $g_{D,\Omega}$ takes non-negative values, is locally integrable on~$M$
and defines  a $\mathrm L_1^2$-Green current for~$D$ in the
sense of~\cite{bost1999}. It is the so-called \emph{equilibrium potential} attached to the divisor $D$ in $\Omega$.

If~$E$ is another effective divisor on~$M$ supported in $\Omega$, 
one has $g_{D+E,\Omega}=g_{D,\Omega}+g_{E,\Omega}$. 
We can therefore extend by linearity
the definition of the equilibrium potential~$g_{D,\Omega}$
to arbitrary divisors~$D$ on~$M$ that are supported on~$\Omega$.
Recall also that, if $\Omega_0$  denotes the union of the connected components
of~$\Omega$ which meet~$\abs D$, then $g_{D,\Omega_0}=g_{D,\Omega}$
(\emph{loc. cit.}, p.~258).

The function $g_{D,\Omega}$ allows one to define
a generalized metric on the line bundle $\mathscr O_M(D)$, by the formula
\[ \norm{\mathbf 1_D}^2 (z)  = \exp(-g_{D,\Omega}(z)),\]
where $\mathbf 1_D$ denotes the canonical global section of~$\mathscr O_M(D)$.
We will call this metric the \emph{capacitary metric}\footnote{Our terminology
differs slightly from that in~\cite{bost1999}. In the present article, the term \emph{capacitary metric} 
will be used for two distinct notions: for the
metrics on line  bundles defined using equilibrium potentials just defined, 
and for some metrics on the tangent line to~$M$ at a point,
see Subsection~\ref{subsec:capacitary-norm}. In~\cite{bost1999}, it was used for the latter notion only.} on~$\mathscr O_M(D)$
 attached to~$\Omega$
and denote by $\norm{f}_\Omega^{\CAP}$ the norm of a local
section~$f$ of~$\mathscr O_M(D)$.

\subsection{Equilibrium potential and capacity on~$p$-adic curves}
\label{sec.capacity}

Let~$R$ be a complete discrete valuation ring, 
and let~$K$ be its field of fractions and $k$ its residue field.
Let~$X$ be a smooth projective curve over~$K$ 
and let~$U$ be an affinoid subspace of
the associated rigid~$K$-analytic curve~$X^\an$.
We shall always require that~$U$ \emph{meets every connected component
of~$X^\an$} --- this hypothesis is analogous to the non-polarity
assumption in the complex case. 
We also let $\Omega=X^\an\setminus U$, which we view
as a (non quasi-compact) rigid~$K$-analytic curve; its affinoid subspaces
are just affinoid subspaces of~$X^\an$ disjoint of~$U$
---
see Appendix~\ref{app.rigid} for a detailed proof that
this endowes~$\Omega$ with the structure of a rigid~$K$-analytic space in the sense of Tate.

The aim of this subsection is to endow  the line
bundle $\mathscr O(D)$, where~$D$ is a divisor which does not meet~$U$,
with a metric (in the sense of~Appendix~\ref{app.metrics}) canonically attached attached to~$\Omega$, in a way that 
parallels the construction over Riemann surfaces
recalled in the previous subsection.

Related constructions of equilibrium potentials over $p$-adic curves have been developed by various authors, notably Rumely~\cite{rumely1989} and Thuillier~\cite{thuillier2005} (see also~\cite{kani1989}). Our approach will be self-contained, and formulated in the framework of classical rigid analytic geometry. Our main tool will be intersection theory on a model $\mathscr X$ of $X$ over $R$. This point of view  will allow us to combine potential theory on $p$-adic curves and Arakelov intersection theory on arithmetic surfaces in a straightforward way. 

We want to indicate that, by using 
an adequate potential theory on analytic curves in
the sense of Berkovich~\cite{berkovich1990} such as the one developped by Thuillier~\cite{thuillier2005}, one could give a treatment of equilibrium potential on $p$-adic curves and their relations to canonical semi-norms that would more closely parallel the one in the complex case. For instance,    in the Berkovich 
setting, the affinoid subspace~$U$ is a compact subset of
the analytic curve attached to~$X$, and $\Omega$
is an open subset. We leave the transposition and the extension of our results in the framework of Berkovich and Thuillier to the interested reader.

\medskip

By Raynaud's general results on formal/rigid geometry,
see for instance~\cite{bosch-l1993,bosch-l1993b},
there exists a normal projective flat model~$\mathscr X$ of~$X$
over~$R$ such that~$U$ is the set of rig points of~$X^\an$
reducing to some open subset~$\mathsf U$ of the special fibre~$\mathsf X$.
We shall write $U=\tube{\mathsf U}_{\mathscr X}$ and say that~$U$ is the
tube of~$\mathsf U$ in~$\mathscr X$;
similarly, we write $\Omega=\tube{\mathscr X\setminus\mathsf U}_{\mathscr X}$.
(We remove the index~$\mathscr X$ from the notation when it is clear
from the context.)
The reduction map identifies the connected components of~$U$ with
those of~$\mathsf U$, and the connected components of~$\Omega$
with those of~$\mathsf X\setminus\mathsf U$.
Since we assumed that $U$ meets every connected component of~$X$,
this shows that $\mathsf U$ meets every connected component of~$\mathsf X$.

Recall that  to  any two Weil divisors $Z_1$ and $Z_2$ on $\mathscr X$ such that $Z_{1,K}$ and $Z_{2,K}$ have disjoint supports is attached their
the \emph{intersection number} $(Z_1,Z_2)$. It is a rational number, which depends linearly on $Z_1$ and $Z_2$. It may be defined \emph{\`a la} Mumford (see~\cite[II.(b)]{mumford1961}), and coincide with the degree over the residue field $k$ of the intersection class $Z_1.Z_2$ in $\operatorname{CH}_0(\mathsf X)$ when $Z_1$ or $Z_2$ is Cartier. Actually, when the residue field $k$ is an algebraic extension of a finite field --- for instance when $K$ is a $p$-adic field, the case we are interested in in the sequel ---  any Weil divisor on $\mathscr X$ has a multiple which is Cartier (see~\cite[Th\'eor\`eme 2.8]{moret-bailly1989}), and this last property, together with their bilinearity, completely determines the intersection numbers.

The definition of  intersection numbers immediately extends by bilinearity to pairs of Weil divisors with coefficients in $\Q$ (shortly, \emph{$\Q$-divisors}) in $\mathscr X$ whose supports do not meet in $X$. 

\begin{proposition}\label{prop.metric}
For any divisor~$D$ on~$X$,
there is a unique $\Q$-divisor~$\mathscr D$ on~$\mathscr X$
extending~$D$ and satisfying the following two conditions:
\begin{enumerate}
\item For any irreducible component~$v$ of codimension~$1$ of~$\mathsf X\setminus \mathsf U$,
$\mathscr D\cdot v=0$.
\item The vertical components of~$\mathscr D$ do not meet~$\mathsf U$.
\end{enumerate}
Moreover, the map $D\mapsto\mathscr D$ so defined is linear
and sends effective divisors to effective divisors.
\end{proposition}
\begin{proof}
Let~$S$ denote  the set of irreducible components of~$\mathsf X$
and let~$T\subset S$ be the subset consisting 
of components which do not meet~$\mathsf U$.
Let~$\mathscr D_0$ be the schematic closure of~$D$ in~$\mathscr X$.
Since $\mathsf U$ meets every connected component of~$\mathsf X$,
$T$ does not contain all of the irreducible components 
of some connected component of~$\mathsf X$,
so that the restriction of the intersection
pairing of~$\Div_\Q(\mathscr X)$ to the subspace generated by 
the components of~$\mathsf X$ which belong to~$T$ is negative  definite
(see for instance~\cite[Corollaire~1.8]{deligne1972b} when $\mathscr X$ is regular; one reduces to this case by considering a resolution of $\mathscr X$, as in~\cite[II.(b)]{mumford1961}).
Therefore, there is a unique vertical divisor~$V$, linear combination
of components in~$T$, such that~$(\mathscr D_0+V,s)=0$
for any~$s\in T$. (In the analogy with the theory of electric networks,
the linear system one has to solve corresponds to that of a Dirichlet
problem on a graph, with at least one electric source per connected
component.)
Set $\mathscr D=\mathscr D_0+V$; it satisfies assumptions~1) and~2).
The linearity of the map $D\mapsto\mathscr D$ follows immediately
from the uniqueness of~$V$.

Let us assume that $D$ is effective and show that so is~$V$.
(In the graph theoretic language, this is a consequence of the maximum
principle for the discrete Laplacian.)
Denote by~$m_s$ the multiplicity of the component~$s$ in the special fibre
of~$\mathscr X$, so that $\sum_{s\in S} m_s s$ belongs to the kernel
of the intersection pairing. Write $V=\sum_{s\in S} c_s s$,
where $c_s=0$ if $s\not\in T$. 

Let $S'$ be the set of elements~$s\in S$ where $c_s/m_s$ achieves
its minimal value. 
Then, for any element~$\tau$ of~$S'\cap T$,
\begin{align*}
0
&  =    (c_\tau/m_\tau) \big(\sum_{s\in S} m_s s, \tau\big)  
 =    c_\tau (\tau,\tau)+ \sum_{s\neq\tau} (c_\tau/m_\tau) m_s (s,\tau) \\
&  \leq c_\tau (\tau,\tau)+\sum_{s\neq\tau} c_s (s,\tau) 
=   \sum_{s\in S} (c_s s,\tau)  \\
&\leq  (\mathscr D,\tau)-(\mathscr D_0,\tau)
=    - (\mathscr D_0,\tau).
\end{align*}
Since $\mathscr D_0$ is effective and horizontal, 
$(\mathscr D_0,\tau)\geq 0$, 
hence all previous inequalities
are in fact equalities. In particular, $(\mathscr D_0,\tau)=0$ and
$c_s/m_s=c_\tau/m_\tau$ for any~$s\in S$ such that $(s,\tau)\neq 0$.

Assume by contradiction that $V$ is not effective, \ie, that there is some~$s$
with $c_s$ negative. Then $S'$ is contained in~$T$ (for $c_s=0$
is $s\not\in T$) and the preceding argument implies that
$S'$ is a union of connected components of~$\mathsf X$. 
(In the graph theoretical analogue, all neighbours of a vertex in~$S'$ 
belong to~$S'$.)
This contradicts the assumption that $U$ meets every connected
component of~$X^\an$ and concludes the proof that $V$ is effective.
\end{proof}

In order to describe the functoriality properties of the assignment $D \mapsto \mathscr D$ constructed in Proposition \ref{prop.metric}, we consider two  smooth projective curves ~$X$ and~$X'$  over~$K$, some
 normal projective 
flat models $\mathscr X$ and~$\mathscr X'$ over~$R$ of these curves, 
and  $\pi\colon\mathscr X' \ra\mathscr X$ 
 an $R$-morphism such that the $K$-morphism $\pi_K\colon X'\ra X$ is finite.
 
 Recall that the direct image of  $1$-dimensional cycles defines a $\Q$-linear map between spaces of $\Q$-divisors:
\[ \pi_*\colon \Div_\Q(\mathscr X')  \longrightarrow \Div_\Q(\mathscr X) ,\]
and that the inverse image of Cartier divisors defines a $\Q$-linear map between spaces of $\Q$-Cartier divisors:
\[ \pi^*\colon \Div^{\rm Cartier}_\Q(\mathscr X) \longrightarrow \Div^{\rm Cartier}_\Q(\mathscr X').\]
These two maps satisfy the following adjunction formula, valid for any $Z$ in $\Div^{\rm Cartier}_\Q(\mathscr X)$  and any  $Z'$ in $\Div_\Q(\mathscr X')$:
\begin{equation}\label{adjunction}
(\pi^* Z, Z')=(Z,\pi_*Z').
\end{equation}

When $k$ is an algebraic extension of a finite field, as recalled above, $\Q$-divisors and $\Q$-Cartier divisors on $\mathscr X$ or $\mathscr X'$ coincide, and $ \pi^*$ may be seen as a linear map from $\Div_\Q(\mathscr X)$ to  $\Div_\Q(\mathscr X')$ adjoint to $\pi_*.$

 In general, the map  $\pi^*$ above admits a unique extension to a $\Q$-linear map 
\[ \pi^*\colon \Div_\Q(\mathscr X)  \longrightarrow \Div_\Q(\mathscr X'),\]
compatible with the pull-back of divisors on the generic fiber
\[ \pi_K^*\colon \Div_\Q(X)  \longrightarrow \Div_\Q( X'),\]
such that the adjunction formula (\ref{adjunction}) holds for any $(Z,Z')$
in $\Div_\Q(\mathscr X)\times\Div_\Q(\mathscr X')$. The unicity of such
a map map follows from the non-degeneracy properties of the intersection
pairing, which show that if a divisor $Z_1'$ supported by the closed fiber
$\mathsf X$ of $\mathscr X$ satisfies $Z'_1\cdot Z'_2=0$ for every $Z'_2$ in
$\Div_\Q(\mathscr X')$, then $Z'_1=0$. The existence of $\pi^*$ is known
when $\mathscr X'$ is regular (then $\Div_\Q(\mathscr X)$ and $\Div^{\rm
Cartier}_\Q(\mathscr X)$ coincide), and when  $\pi$ is birational ---
\emph{i.e.,} when $\pi_K$ is an isomorphism --- and $\mathscr X$ is
regular, according to Mumford's construction in~\cite[II.(b)]{mumford1961}. 
To deal with the general case, observe that there exist two
projective flat regular curves $\tilde{\mathscr X}$ and $\tilde{\mathscr
X'}$ equipped with birational $R$-morphisms $\nu\colon \tilde{\mathscr X}
\ra \mathscr X$ and $\nu'\colon \tilde{\mathscr X'} \ra \mathscr X'$, and an
$R$-morphism $\tilde{\pi}\colon \tilde{\mathscr X'} \ra \tilde{\mathscr X}$
such that $\pi \circ \tilde{\nu} = \tilde{\pi} \circ \nu$. Then it is
straightforward that $\pi^*:= \tilde{\nu}_*\tilde{\pi}^*\nu^*$ satisfies
the required properties.

Observe also that the assignment $\pi \mapsto \pi^*$ so defined is
functorial, as follows easily from its definition.

\begin{proposition}\label{prop.capa-functorial}
Let $\mathsf U$ be a Zariski open subset of the special
fibre~$\mathsf X$ and let $\mathsf U'=\pi^{-1}(\mathsf U)$.
Assume that $\tube{\mathsf U}_{\mathscr X}$ meets
every connected component of~$X^\an$; then $\tube{\mathsf U'}_{\mathscr X'}$
meets every connected component of~$(X')^\an$.

Let~$D$ and~$D'$ be divisors on~$X$ and~$X'$ respectively,
let $\mathscr D$ and~$\mathscr D'$ be the extensions
to~$\mathscr X$ and~$\mathscr X'$, relative to the open subsets
$\mathsf U$ and~$\mathsf U'$ respectively, given
by Proposition~\ref{prop.metric}.

\begin{enumerate}\def\theenumi{\alph{enumi}}\def\labelenumi{\theenumi)}
\item
Assume that $D'=\pi^*D$.
If $\abs D$ does not meet~$\tube{\mathsf U}$, then $\abs {D'}$ is disjoint
from~$\tube{\mathsf U'}$ and  $\mathscr D'=\pi^*\mathscr D$.

\item
Assume that $D=\pi_*D'$.
If $\abs {D'}$ does not meet~$\tube{\mathsf U'}$, then $\abs{D}\cap\tube{\mathsf U}=\emptyset $ and $\mathscr D=\pi_*\mathscr D'$.
\end{enumerate}
\end{proposition}

\begin{proof} % [Proof of Proposition~\ref{prop.capa-functorial}]
Let us denote by~$S$ the set of irreducible components
of the closed fibre~$\mathsf X$ of~$\mathscr X$,
and by~$T$ its subset 
of the components which do not meet~$\mathsf U$.
Define similarly $S'$  and~$T'$ to be the set of irreducible
components of~$\mathsf X'$ and its subset
corresponding to the components that do not meet~$\mathsf U'$.
Let also~$N$ denote the set of all irreducible components
of~$\mathsf X'$ which are contracted to a point by~$\pi$.

By construction of $\pi^*$, the divisor~$\pi^*(\mathscr D)$ satisfies
$(\pi^*(\mathscr D),n)=0$ for any~$n\in N$ and has no
multiplicity along the components of~$N$ that are not
contained in~$\pi^{-1}(\abs{\mathscr D})$.

Since $\mathsf U'=\pi^{-1}(\mathsf U)$,
$T'$ is the union of all components of~$\mathsf X'$
that are mapped by~$\pi$, either to a point outside~$\mathsf U$,
or to a component in~$T$.

\emph a)
Let $t'\in T'$.
One has $(\pi^*\mathscr D,t')=(\mathscr D,\pi_*t')=0$
since $t'$ maps to a component in~$T$,  or to a point.
Moreover, by the construction
of~$\pi^*$, the vertical components of~$\pi^*\mathscr D$
are elements $s'\in S'$ such $\pi(s')$ meets the support of~$\mathscr D$.
By assumption, the Zariski closure of~$D$ in~$\mathsf X$ is disjoint
from~$\mathsf U$; in other words, the vertical components of~$\pi^*\mathscr D$
all belong to~$T'$.
This shows that the divisor~$\pi^*\mathscr D$ on~$\mathscr X'$
satisfies the conditions of Proposition~\ref{prop.metric};
since it extends~$D'=\pi^*D$, one has
$\pi^*\mathscr D=\mathscr D'$.

\emph b)
Let~$s$ be a vertical component appearing in~$\pi_*(\mathscr D')$;
necessarily, there is a vertical component $s'$ of~$\mathscr D'$ 
such that $s=\pi(s')$. This implies that $s'\in T'$, hence $s\in T$.
For any $t\in T$, $\pi^*(t)$ is a linear combination of vertical
components of~$\mathsf X'$ contained in~$\pi^{-1}(t)$.
Consequently, they all belong to~$T'$ and
one has
$
(\pi_*(\mathscr D'),t) = (\mathscr D',\pi^*(t))= 0
$.
By uniqueness, $\pi_*(\mathscr D)=\mathscr E$.
\end{proof}

\begin{corollary}\label{lemm.indep}
Let~$X$ be a projective smooth algebraic curve over~$K$,
let~$U$ be an affinoid subspace of~$X^{\rm an}$ which meets 
any connected component of~$X^{\rm an}$.
Let~$D$ be a divisor on~$X$ whose support is disjoint from~$U$.

Then the metrics on the line bundle~$\mathscr O_X(D)$ induced
by the line bundle~$\mathscr O_{\mathscr X}(\mathscr D)$
defined by Prop~\ref{prop.metric}
does not depend on the choice of the projective flat model~$\mathscr X$
of~$X$ such that~$U$ is the tube of a Zariski open subset of
the special fibre of~$\mathscr X$.
\end{corollary}
\begin{proof}
For $i=1,2$, let $(\mathscr X_i,\mathsf U_i)$ 
be a pair as above, consisting of a normal flat, 
projective model~$\mathscr X_i$
of~$X$ over~$R$ and an open subset~$\mathsf U_i$ 
of its special fibre~$\mathsf X_i$ 
such that $\tube{\mathsf U_i}_{\mathscr X_i}=U$.
Let $\mathscr D_i$ denote the extension of~$D$
on~$\mathscr X_i$ relative to~$\mathsf U_i$.

There exists a third model~$(\mathscr X',\mathsf U')$
which admits maps $\pi_i\colon\mathscr X'\ra\mathscr X_i$, for $i=1,2$,
extending the identity on the generic fibre. Let $\mathscr D'$ denote
the extension of~$D$ on~$\mathscr X'$.
For $i=1,2$, one has $\pi_i^{-1}(\mathsf U_i)=\mathsf U'$.
By Proposition~\ref{prop.capa-functorial},
one thus has the equalities $\pi^*\mathscr D_i=\mathscr D'$
hence the line bundles $\mathscr O_{\mathscr X'}(\mathscr D')$
on~$\mathscr X'$
and $\mathscr O_{\mathscr X_i}(\mathscr D_i)$ on~$\mathscr X$
induce the same metric on~$\mathscr O_X(D)$.
\end{proof}

We shall call this
metric the \emph{capacitary metric}
and denote as $\norm{f}_{\Omega}^{\CAP}$ the norm
of a local section~$f$ of~$\mathscr O_X(D)$
for this metric.

\begin{proposition}\label{prop.capa.connected}
Let~$X$ be a projective smooth algebraic curve over~$K$,
let~$U$ be an affinoid subspace of~$X^\an$ which meets 
any connected component of~$X^\an$.
Let~$D$ be a divisor on~$X$ whose support is disjoint from~$U$
and let $\Omega=X^\an\setminus U$.

If $\Omega'$ denotes the union of the connected components
of~$\Omega$ which meet~$\abs D$,
then the capacitary metrics of~$\mathscr O(D)$ relative
to~$\Omega$ and to~$\Omega'$ coincide.
\end{proposition}
\begin{proof}
Let us fix a normal projective flat model $\mathscr X$ of~$X$ over~$R$
and a Zariski open subset~$\mathsf U$ of its special fibre~$\mathsf X$
such that $U=\tube{\mathsf U}_{\mathscr X}$.
Let $\mathsf Z=\mathsf X\setminus\mathsf U$ and
let $\mathsf Z'$ denote the union of those connected components
of~$\mathsf Z$ which meet the specialization of~$\abs D$.
Then $\Omega'=\tube{\mathsf Z'}$ is the complementary subset
to the affinoid $\tube{\mathsf U'}$, where $\mathsf U'=\mathsf X\setminus \mathsf Z'$; in particular, $\tube{\mathsf U'}$ meets every
connected component of~$X^\an$.

Let $\mathscr D_0$ denote the horizontal divisor on~$\mathscr X$
which extends~$D$.
The divisor~$\mathscr D' \mathrel{:=} \mathscr D_{\Omega'}$
is the unique $\Q$-divisor of the form $\mathscr D_0+V$ on~$\mathscr X$
where~$V$ is a vertical divisor supported by~$\mathsf Z'$ 
such that $(\mathscr D',t)=0$ for any irreducible component
of~$\mathsf Z'$.
By the definition of~$\mathsf Z'$, an irreducible component of~$\mathsf Z$ 
which is not contained in~$\mathsf Z'$ doesn't meet
neither~$\mathsf Z'$, nor~$\mathscr D_0$.
It follows that for any such component~$t$, $(\mathscr D',t)=(\mathscr D_0,t)+(V,t)=0$.
By uniqueness, $\mathscr D'$ is the extension of~$D$ on~$\mathscr X$
relative to~$\mathsf U$, so that $\mathscr D_{\Omega'}=\mathscr D_{\Omega}$.
This implies the proposition.
\end{proof}

As an application of the capacitary metric, in the next proposition
we establish a variant of  a classical
theorem by Fresnel and Matignon
(\cite[th\'eor\`eme~1]{fresnel-m1986})
asserting that affinoids of a curve can be defined by one equation.
(While these authors make no hypothesis
on the residue field of~$k$, or on
the complementary subset of the affinoid~$U$, 
we are able to impose the  polar divisor of~$f$.)
Using the terminology of Rumely \cite[\S4.2, p.~220]{rumely1989}, 
this proposition means that affinoid subsets of a curve are RL-domains
(``rational lemniscates''),
and that RL-domains with connected complement are PL-domains
(``polynomial lemniscates'').
It is thus essentially equivalent to Rumely's theorem 
\cite[Theorem~4.2.12, p.~244]{rumely1989} asserting
that \emph{island domains} coincide with PL-domains.
Rumely's proof relies on his non-archimedean potential theory,
which we replace here by Proposition~\ref{prop.metric}.

This proposition will also be used to derive further properties
of the capacitary metric.

\begin{proposition}\label{prop.equations}
Assume that the residue field~$k$ of~$K$ is algebraic over
a finite field.
Let~$(X,U,\Omega)$ be as above, and let~$D$
be an effective divisor which does not meet~$U$
but meets every connected component of~$\Omega$.
There is a rational function~$f\in K(X)$ with polar divisor
a multiple of~$D$
such that
$U=\{x\in X\sozat \abs{f(x)}\leq 1\}$.
\end{proposition}
\begin{proof}
Keep notations as in the proof of Proposition~\ref{prop.metric};
in particular, $S$ denotes the set of irreducible
components of~$\mathsf X$. 
The closed subset~$\mathsf X\setminus\mathsf  U$ has only finitely
many connected components, say $\mathsf V_1,\dots,\mathsf V_r$. Moreover,
we may assume that for each~$i$, $\mathsf V_i$ is the union
of a family $T_i\subset S$ of components of~$\mathsf X$.
For any~$i$, the tube $\tube{\mathsf V_i}$ in~$X^\an$
consisting of the rig. points of~$X$
which reduce to points of~$\mathsf V_i$  is a connected  analytic  subset
of~$X^\an$, albeit not quasi-compact, and $X^\an$ is the disjoint
union of~$U=\tube{\mathsf U}$ and of the~$\tube{\mathsf V_i}$.
(See~\cite{raynaud94} for more details.)
We let~$m_s$ denote the multiplicity of the component~$s$
in the special fibre, and $F=\sum_{s\in S} m_s s$.

Let $\mathscr D=\mathscr D_0+V$ be the extension of~$D$ to
a $\Q$-divisor of~$\mathscr X$ given
by Proposition~\ref{prop.metric},
where $\mathscr D_0$ is horizontal and~$V=\sum_{s\in S} c_s s$ 
is a vertical divisor supported by the special fibre~$\mathsf X$.
One has $c_s=0$ for $s\not\in T$ and $c_s\geq 0$ if $s\in T$.
For any $s\not\in T$, we define $a_s=(V,s)$.
This is a nonnegative
rational number and we have
\[ \sum_{s\in S\setminus T} a_s m_s
= (V,F) - \sum_{t\in T} m_t(V,t)
= \sum_{t\in T} m_t (\mathscr D_0,t)
= \sum_{s\in S} m_s (\mathscr D_0,s)
\]
since~$D$ does not meet~$U$,
hence 
\begin{equation}\label{eq.sum.asms}
 \sum_{s\in S\setminus T} a_s m_s = (\mathscr D_0,F) = \deg(D).
\end{equation}

\smallskip

For any $s\in S\setminus T$,
let us fix
 a point~$x_s$ of~$\mathsf X$ which is contained on the component~$s$
as well as on the smooth locus of~$\mathscr X$.
Using either a theorem of Rumely \cite[Theorem~1.3.1, p.~48]{rumely1989}, 
or van der Put's description of the Picard group of any one-dimensional 
$K$-affinoid, \cf\cite[Proposition~3.1]{vanderput1980},\footnote
  {The proofs in both references
  are similar and rely on the Abel-Jacobi map, together
  with the fact that~$K$ is the union of its locally compact
  subfields.}
there is a rational function~$f_s\in K(X)$
with polar divisor a multiple of~$D$
and of which all zeroes specialize to~$x_s$.
We may write its divisor as a sum
\[ \div(f_s)=-n_s \mathscr D + E_s + W_s, \]
where $n_s$ is a positive integer, $E_s$ is a horizontal
effective divisor having no common component with~$\mathscr D$,
and $W_s$ is a vertical divisor.
Since $E_s$ is the closure of the divisor of zeroes of~$f_s$, it only
meets the component labelled~$s$. One thus
has $(E_s,s')=0$ for $s'\in S\setminus\{s\}$, while
\[ (E_s,s)=\frac1{m_s} (E_s,m_s s)= \frac1{m_s}(E_s,F)=\frac{n_s}{m_s}\deg D.\]

Let $t\in T$. One has $(\div(f_s),t)=0$, hence
\[ (W_s,t) = n_s(\mathscr D,t) - (E_s,t) = 0. \]
Similarly, if $s'\in S\setminus T$, $(\div(f_s),s')=0$ and
\begin{align*}
 (W_s,  s') & = n_s (\mathscr D,  s')-(E_s, s') \\
& = n_s (\mathscr D_0, s') + n_s(V, s') - (E_s,s') \\
& =  0 + n_s a_{s'}   - (E_s,s'). 
\end{align*}
If $s'\neq s$, it follows that
\[ (W_s,s')= n_s a_{s'}  , \]
while
\[ (W_s,s) = n_s a_s - \frac{n_s}{m_s} \deg(D). \]

We now define a vertical divisor 
\[ W = \sum_{s\not\in T} \frac{a_s m_s}{n_s} W_s. \]
For any~$t\in T$, $(W,t)=0$. Moreover, for any $s'\in S\setminus T$,
\begin{align*}
 (W, s') &
=\sum_{s\not\in T} \frac{ a_s m_s}{n_s} (W_s,s') \\
&= \sum_{s\not\in T} a_s m_s a_{s'} 
   - \frac{a_{s'}}{n_{s'}} n_{s'} \deg(D) \\
& = a_{s'}  \left(\sum_{s\not\in T} a_sm_s\right) 
   - a_{s'} \deg(D) ,
\end{align*}
hence $(W,s')=0$ by~\eqref{eq.sum.asms}.
Therefore, the vertical $\Q$-divisor~$W$ is a multiple of the special fibre
and there is $\lambda\in\Q$ such that $W=\lambda F$.
Finally,
\[ \sum_{s\not\in T} \frac{a_s m_s}{n_s} \div(f_s)
 -\lambda F = - \deg(D) \mathscr D
 + \sum_{s\not\in T} \frac{a_s m_s}{n_s} E_s \]
is a principal $\Q$-divisor.
It follows that there are positive integers $\mu$
and~$\lambda_s$, for $s\not\in T$,
such that
\[ \mathscr P = \sum_{s\not\in T} \lambda_s E_s - \mu \mathscr D \]
is the divisor of a rational function $f\in K(X)$.

By construction, the polar divisor of~$f$ on~$X$
is a multiple of~$D$.
Moreover, the reduction of any $x\not\in U$ belongs
to a component labelled by~$T$ at which the multiplicity
of~$\mathscr P$ is positive. Consequently, $\abs{f(x)}>1$.
On the contrary, if $x\in U$, it reduces to a component
outside~$T$ and $\abs{f(x)}\leq 1$. More precisely, 
$\abs{f(x)}<1$ if and only if~$x$ reduces to one of the points~$x_s$,
$s\not\in T$.
\end{proof}

The definition of an algebraic metric now implies the
following explicit description of the capacitary metric.
\begin{corollary}
Let $(X,U,\Omega)$ be as above
let~$D$ be any divisor which does not meet~$U$,
and let~$f$ be a rational function defining~$U$, as in 
preceding proposition, and whose polar divisor
is equal to~$m D$, for some positive integer~$m$.
Then, the capacitary metric on~$\mathscr O_{X}(D)$
can be computed as
\[ - \log \norm{1_{\mathsf D}}^{\CAP}_{\Omega}(x) = \frac1m \log^+ \abs{f(x)}
= \max (0,\log\abs{f(x)}^{1/m}). \]
\end{corollary}
\begin{proposition}\label{prop.capa.isom}
Let $(X,U,\Omega)$ and $(X',U',\Omega')$ be as above
% , denote $\Omega'=(X')^\an\setminus U'$
% and $\Omega=X^\an\setminus U$ 
and let $\phi\colon \Omega'\ra\Omega$
be any rigid analytic isomorphism.
Let~$D'$ be any divisor in~$X'$ 
whose support does not meet~$U'$ and let~$D=\phi(D')$.

Then, for any $x\in \Omega'$,
\[ \norm{1_{D'}}^\CAP_{\Omega'}(x) = \norm{1_{D}}^\CAP_{\Omega}(\phi(x)). \]
\end{proposition}
\begin{proof}
By linearity, we may assume that~$D$ is effective.
Let $f\in K(X)$ and $f'\in K(X')$ be rational
functions as in Proposition~\ref{prop.equations}.
Let~$m$ and $m'$ be positive integers such that the polar divisor
of~$f$ and~$f'$ are $mD$ and $m'D'$ respectively.
The function $f\circ\phi$ is a meromorphic function on~$\Omega'$
whose divisor is $mD'$.
Consequently, the meromorphic function on~$\Omega'$ 
\[ g= (f\circ \phi)^{m'}/(f')^m \]
is in fact invertible. We have to prove that $\abs g(x)=1$
for any $x\in \Omega'$.

Let $(\eps_n)$ be any decreasing sequence of elements
of~$\sqrt{\abs{K^*}}$
converging to~$1$. The sets $V'_n=\{x\in X'\sozat
\abs{f'(x)}\geq \eps_n\}$ are affinoid subspaces of~$\Omega'$
and exhaust it.
By the maximum principle (see Proposition~\ref{prop.max-principle} below),
one has
\[ \sup_{x\in V'_n} \abs{g(x)} = \sup_{\abs{f'(x)}=\eps_n}
           \abs{g(x)} \leq  1/ (\eps_n)^m \leq 1. \]
Consequently,
$ \sup_{x\in\Omega'} \abs{g(x)} \leq 1 $.
The opposite inequality is shown similarly by
considering the isomorphism $\phi^{-1}\colon\Omega\ra\Omega'$.
This proves the proposition.
\end{proof}

\subsection{Capacitary norms on tangent spaces}
\label{subsec:capacitary-norm}
\begin{definition}
Let $(X,U,\Omega)$ be as above
and let $P\in X(K)$ be a rational point lying in~$\Omega$.
Let us endow the line bundle $\mathscr O_X(P)$  
with its capacitary metric relative to~$\Omega$.
The capacitary norm 
$\norm{\cdot}_{P,\Omega}^\CAP$
on the~$K$-line $T_P X$ is then defined
as the restriction of~$(\mathscr O_X(P),\norm{\cdot}_{\Omega}^\CAP)$
to the point~$P$, composed with 
the adjunction isomorphism $\mathscr O_X(P)|_P\simeq T_PX$.
\end{definition}

\begin{example}
Let us fix a normal projective flat model~$\mathscr X$,
let~$\mathscr P$ be the divisor extending~$P$,
meeting the special fibre~$\mathsf X$ in a smooth point~$\mathsf P$.
Let $\mathsf U=\mathsf X\setminus\{\mathsf P\}$
and define $U=\tube{\mathsf U}$, $\Omega=\tube{\mathsf P}$.
In other words, $\Omega$ is the set of rig-points of~$X^\an$
which have the same reduction~$\mathsf P$ as~$P$.
Then $\Omega=\tube {\mathsf P}$ is isomorphic to an open unit ball,
the divisor~$\mathscr P$ is simply the image of the section
which extends the point~$P$,
and the capacitary metric on~$T_PS$ is simply the metric induced
by the integral model.
\end{example}

\begin{example}[Comparison with other definitions]
Let us show that how this norm 
fits with Rumely's definition in~\cite{rumely1989}
of the capacity of~$U$ with respect to the point~$P$.
Let~$f$ be a rational function on~$X$, without pole
except~$P$,  such that
$U=\{x\in X\sozat \abs{f(x)}\leq 1\}$.
Let~$m$ be the order of~$f$ at~$P$ and let us define $c_P\in K^*$ so that
$f(x)=c_P t(x)^{-m}+\dots$ around~$P$, where~$t$ is a fixed local
parameter at~$P$.
By definition of the adjunction map, the local section $\frac1t 1_P$
of~$\mathscr O_X(P)$  maps to the tangent vector~$\frac\partial{\partial t}$.
Consequently,
\begin{equation}
 \norm{\dfrac\partial{\partial t}}_{P,\Omega}^\CAP
  = \norm{ \dfrac1t 1_P} (P) 
  = \lim_{x\ra P} \abs{t(x)}^{-1} \min (1, \abs{f(x)}^{-1/m})
 = \abs{c_P}^{-1/m}. \end{equation}

As an example, and to make explicit the relation of our rationality
criterion below with the classical theorem of Borel-Dwork later on,
let us consider the classical case in 
which $X=\mathbf P^1$ (containing the affine line with~$t$ coordinate),
and~$U$ is the affinoid subspace of~$\mathbf P^1$
defined by the inequality $\abs t \geq r$ (to which we add the point
at infinity), where $r\in\sqrt{\abs {K^*}}$. Let us note
$\Omega=\complement U$ and choose for
the point~$P\in\Omega$ the point with coordinate $t=0$.
Let~$m$ be a positive integer and $a\in K^*$ such that $r^m=\abs a$;
let $f=a/t^m$; this is a rational function on~$\P^1$ with a single
pole at~$P$ and~$U$ is defined by the inequality $\abs f\leq 1$.
It follows that
\[ \norm{\dfrac \partial{\partial t}}_{P,\Omega}^\CAP
 = \abs a^{-1/m}= 1/r. \]
Similarly, assume that~$U$ is an affinoid subset of~$\mathbf P^1$
which does not contain the point~$P=\infty$. Then~$U$
is bounded and 
$ \norm{t^2\frac\partial {\partial t}}_{P,\Omega} $.
is nothing but its transfinite diameter in the sense of Fekete.
(See~\cite{amice75}, the equivalence of both notions
follows from~\cite[Theorem~4.1.19, p.~204]{rumely1989};
see also \cite[Theorem~3.1.18, p.~151]{rumely1989} for its archimedean
counterpart.)
\end{example}

\begin{remarks}\label{rema.capa.inequalities}
\begin{enumerate}\def\theenumi{\alph{enumi}}\def\labelenumi{\theenumi)}
\item
Let~$(X,U)$ be as above, let $P\in X(K)$ be a rational point
such that $P\not\in U$. Let $\Omega=X^\an\setminus U$ and
define~$\Omega_0$ to be the connected component of~$\Omega$ which contains~$P$.
It follows from Proposition~\ref{prop.capa.connected} that 
the norms $\norm{\cdot}_{P,\Omega_0}^\CAP$
and $\norm{\cdot}_{P,\Omega}^\CAP$ on~$T_PX$ coincide.

\item
Let $U'$ be another affinoid subspace of~$X^\an$ such that $U'\subset U$;
the complementary subset~$\Omega'$ to~$U'$ satisfies $\Omega\subset\Omega'$.
If moreover $\Omega$ and~$\Omega'$ are connected, then
for any $P\in \Omega$ and any vector~$v\in T_PX$, one has
\[ \norm{v}_{P,\Omega'}^\CAP \leq \norm{v}_{P,\Omega}^\CAP. \]

Indeed, since $\Omega$ and~$\Omega'$ are connected
and contain~$P$,  Proposition~\ref{prop.equations}
implies that there exist rational functions~$f$ and~$f'$ on~$X$,
without pole except~$P$, such that
the affinoids~$U$ and~$U'$ are defined by the inequalities
$\abs f\leq 1$ and $\abs {f'}\leq 1$ respectively.
Replacing~$f$ and~$f'$ by some positive powers, 
we may also assume that $\ord_P(f)=\ord_P(f')$; let us denote it by~$-d$.
Let~$t$ be a local parameter at~$P$; it is enough
to prove the desired inequality for $v=\frac{\partial}{\partial t}$.

We may expand~$f$ and~$f'$ around~$P$ as Laurent series in~$t-t(P)$,
writing
\[ f = \frac{c}{(t-t(P))^d} + \dots, \quad f'= \frac{c'}{(t-t(P))^d}+\dots \quad.
\]
The rational function
$g=f/f'$ on~$X$ defines a holomorphic  function on the affinoid
subspace defined by the inequality $\{\abs {f'}\geq 1\}$, since
the poles at~$P$ at the numerator and at the denominator cancel
each other; moreover, $g(P)=c/c'$.
Using twice the maximum  principle (Proposition~\ref{prop.max-principle}),
we have
\begin{align*} \abs{g(P)} & \leq \sup_{\abs{f'(x)}\geq 1} \abs{g(x)}
 = \sup_{\abs{f'(x)}=1} \abs{g(x)}
 = \sup_{\abs{f'(x)}=1} \abs{f(x)} \\
&  = \sup_{\abs{f'(x)}\leq 1} \abs{f(x)} \leq 1 \end{align*}
since $\Omega\subset\Omega'$. 
This implies that $\abs{g(P)}\leq 1$, so that $\abs{c}\leq\abs {c'}$.
Therefore,
\[ \norm{\frac{\partial}{\partial t}}_{P,\Omega'}^\CAP
= \abs{c'}^{-1/d}  \leq \abs{c}^{-1/d} =  \norm{\frac{\partial}{\partial t}}_{P,\Omega}^\CAP,
\]
as was to be shown.
\end{enumerate}
\end{remarks}

\subsection{Canonical semi-norms and capacities}
\label{sec.schwarz}

Let~$K$ be a local field.

In the case~$K$ is archimedean, we assume moreover that $K=\C$;
let~$M$ be a connected Riemann surface, $\Omega$ be an open
subset in~$M$, relatively compact. 
In the case~$K$ is ultrametric, let~$M$ be a smooth projective
curve over~$K$, let~$U$ be an affinoid in~$M^\an$, let
us denote $\Omega=M^\an\setminus U$. 

In both cases, let~$O$ be a point in~$\Omega$.

We endow the~$K$-line $T_OM$ with its capacitary semi-norm,
as defined by the first author in~\cite{bost2004} when $K=\C$,
or in the previous section in the~$p$-adic case.

\smallskip

Let~$X$ be a  projective variety over~$K$, 
let $P\in X(K)$ be a rational point and let $\widehat C$ be a 
smooth formal curve in~$\widehat X_P$. 
Assume that $\widehat C$ is $K$-analytic
and let $\phi\colon \Omega\ra X^\an$ be an analytic map
such that $\phi(O)=P$ which maps the germ of~$\Omega$ at~$O$
to~$\widehat C$. (Consequently, if $D\phi(O)\neq 0$,
then $\phi$ defines an analytic isomorphism
from the formal germ of~$\Omega$ at~$O$ to~$\widehat C$.)
We endow $T_PX$ with its canonical semi-norm
$\norm{\cdot}^\can_{X,\widehat C}$.

\begin{proposition}
For any $v\in T_O\Omega$, one has
\[ \norm{D\phi(O) (v)}_{X,\widehat C}^\can \leq \norm{v}_{P,\Omega}^\CAP. \]
\end{proposition}
\begin{proof}
The case $K=\C$ is treated in~\cite[Proposition~3.6]{bost2004}.
It therefore remains to treat the ultrametric case.

In view of remark~\ref{rema.capa.inequalities}, \emph a),
we may assume that $\Omega$ is connected.
By Proposition~\ref{prop.equations},
there exists a rational function $f\in K(M)$ without pole
except~$O$ such that $U=\{x\in M\sozat \abs{f(x)}\leq 1\}$.
Let $m>0$ denote the order of the pole of~$f$ at the point~$O$.
For any real number $r>1$ belonging to $\sqrt{\abs{K}^*}$, let us denote
by $U_r$ and~$\partial U_r$ the affinoids $\{\abs{f(x)}\geq r\}$ and
$\{\abs{f(x)}=r\}$ in~$M$.
One has $\bigcup_{r>1}U_r=\Omega$. We shall denote by
$\phi_r$ the restriction of~$\phi$ to the affinoid~$U_r$.
Let us also fix a local parameter~$t$ at~$O$ and let us define
$c_P=\lim\limits_{x\ra O}t(x)^m f(x)$.
One has $\norm{\frac\partial{\partial t}}_{O,\Omega}^\CAP=\abs{c_P}^{-1/m}$.

Let~$L$ be an ample line bundle on~$X$
For the proof of the proposition,
we may assume that $D\phi(O)$ is non zero;
then $\phi$ is a formal isomorphism and we may consider
the formal parameter $\tau=t\circ\phi^{-1}$ on~$\widehat C$ at~$P$.
We have $dt=\phi^*d\tau$, hence
$D\phi(O)(\frac\partial{\partial t})=\frac\partial{\partial\tau}$.
Let us also fix a norm~$\norm{\cdot}_0$ on the~$K$-line~$T_P\widehat C$, 
and let us still denote by~$\norm{\cdot}_0$
the associated norm on on its dual $T^\vee_P\widehat C$.

Let us choose a real number $r>1$ such that $r\in\sqrt{\abs{K^*}}$,
fixed for the moment.
Since the residue field of~$K$ is finite,
the line bundle $\phi_r^*L$ on~$U_r$ is torsion 
(see~\cite[Proposition~3.1]{vanderput1980}); we may therefore
consider a positive integer~$n$ and a nonvanishing section~$\eps$
of~$\phi_r^*L^{\otimes n}$.
For any integer~$D$ and any section 
$s\in\Gamma(X,L^{\otimes nD})$, let us write
$\phi_r^* s = \sigma \eps^{\otimes nD}$, where
$\sigma$ is an analytic function on~$U_r$.
Since we assumed that $D\phi(O)\neq 0$,
the condition that~$s$ vanishes at order~$i$ along~$\widehat C$
means exactly that $\sigma$ vanishes at order~$i$ at~$O$.
Consequently, the~$i$-th jet of~$\phi_r^*s$ at~$O$ is given by
\[ \jet^i_O(\phi_r^*s) = (\sigma t^{-i})(O) \eps^{nD}(O) \otimes
d\tau^{\otimes i}. \]
Writing $(\sigma t^{-i})^m = (\sigma^m f^i)(ft^m)^{-i}$,
it  follows that
\[ \norm{\jet^i_{O}(\phi_r^*s)}^m
         = \abs{\sigma^m f^i}(O) \abs{c_P}^{-i} \norm{\eps(O)}^{nmD}
\norm{d\tau}_0^{im} .\]
Notice that $\sigma^mf^{i} $ is an analytic function on~$U_r$.
By the maximum
principle (Proposition~\ref{prop.max-principle}),
\[ \abs{\sigma^m f^i}(O) \leq \sup_{U_r} \abs{\sigma^m f^i}
           = \sup_{x\in \partial U_r} \abs{\sigma^m f^i(x)}
        = \norm{\sigma}_{\partial U_r}^m   r^i. \]
Consequently,
\begin{align*}
 \norm{\jet^i_{O}(s)}
&\leq  \norm{\sigma}_{\partial U_r}
 \abs{c_P}^{-i/m} r^i \norm{\eps(O)}^{nD} \norm{d\tau}_0^i  \\
&\leq \norm{s}_{\partial U_r} \abs{c_P}^{-i/m} r^i \left(\frac{\norm{\eps(O)}}
{\inf_{x\in \partial U_r} \norm{\eps(x)}} \right)^{nD} \norm{d\tau}_0^i.
\end{align*}
With the notations of Section~\ref{subsec.canonical},
it follows that the norm of
the evaluation morphism 
\[ \phi^i_{nD}\colon E^i_{nD}
\ra L^{\otimes nD}_{\mid P}\otimes (T^\vee_P\widehat C)^{\otimes i} \]
satisfies the inequality 
\[ \norm{\phi^i_{nD}}^{1/i}
 \leq  r^{1/m} \abs{c_P}^{-1/m} \left ( \norm{\eps(O)}/\inf_{\partial
U_r}\norm{\eps} \right)^{nD/i} 
\norm{d\tau}_0 , \]
hence 
\[ 
 = \limsup_{i/D\ra\infty} \frac1i\log\norm{\phi^i_{nD}}
\leq \frac1m \log\frac{r}{\abs{c_P}} = \log\norm{d\tau}_0 .\]
Using the notations introduced for defining the canonical semi-norm,
we thus have $ \rho(L)=\rho(L^{\otimes n})\leq \log\norm{d\tau}_0$ and
\begin{align*}
 \norm{D\phi(O)(\frac\partial{\partial t})}^\can_{X,\widehat C,P}
&= \norm{\frac\partial{\partial \tau}}^\can_{X,\widehat C,P}
= e^{\rho(L)}\norm{\frac\partial{\partial \tau}}_0  \\
& \leq
\left(\frac{r}{\abs{c_P}}\right)^{1/m} 
= r^{1/m} \norm{\frac\partial{\partial t}}^\CAP_{\Omega,P}. \end{align*}
Letting~$r$ go to~$1$, we obtain the desired inequality.
\end{proof}

\subsection{Global capacities}

Let~$K$ be a number field, let~$R$ denote the
ring of integers in~$K$.
Let~$X$ be a projective smooth algebraic curve over~$K$.
For any ultrametric place~$v$ of~$R$, let us denote
by~$\F_v$ the residue field of~$R$ at~$v$,
by~$K_v$ the completion of~$K$ at~$v$, and by~$X_v$ the rigid
$K_v$-analytic variety attached to~$X_{K_v}$.
For any archimedean place~$v$ of~$X$, corresponding
to an embedding $\sigma\colon K\hra\C$, we
let $X_v$ be the compact Riemann surface $X_\sigma(\C)$. When $v$ is real, by an open subset of $X_v$, we shall mean an open subset of $X_\sigma(\C)$ invariant under complex conjugation.

Our goal in this Section is to show how capacitary metrics
at all places fit within the framework of the Arakelov
intersection theory (with $\mathrm L_1^2$-regularity) 
introduced in~\cite{bost1999}.
Let us briefly recall here the main notations and properties
of this arithmetic intersection theory, referring to
this article for more details.

For any normal projective flat model~$\mathscr X$ of~$X$ over~$R$,
the Arakelov Chow group $\hCH^1_\R(\mathscr X)$ 
consists of equivalence classes of pairs 
$(\mathscr D,g)\in\hZ^1_\R(\mathscr X)$
where $\mathscr D$ is a $\R$-divisor on~$\mathscr X$ 
and~$g$ is a Green current with $\mathrm L_1^2$-regularity
on~$\mathscr X(\C)$ for the real divisor~$\mathscr D_K$,
stable under complex conjugation. For any
class~$\alpha$ of an Arakelov divisor $(\mathscr D,g)$,
we shall denote, as usual,  $\omega(\alpha)=dd^c g+\delta_{\mathscr D_\C}$.

Arithmetic intersection theory endowes the space~$\hCH^1_\R(\mathscr X)$
with a symmetric $\R$-valued bilinear form.
Any morphism $\pi\colon\mathscr X'\ra\mathscr X$ between normal projective
flat models of curves~$X'$ and~$X$ induces morphisms of abelian groups
$\pi_*\colon \hCH^1_\R(\mathscr X')\ra\hCH^1_\R(\mathscr X)$
and $\pi^*\colon\hCH^1_\R(\mathscr X)\ra\hCH^1_\R(\mathscr X')$.
For any classes $\alpha$ and~$\beta\in\hCH^1_\R(\mathscr X)$,
$\gamma\in\hCH^1_\R(\mathscr X')$, one has
$\pi^*\alpha\cdot\pi^*\beta=\alpha\cdot\beta$
and a projection formula
$\pi_*(\pi^*\alpha\cdot\gamma)=\deg(\pi)\alpha\cdot\pi_*(\gamma)$,
when $\pi$ has constant generic degree~$\deg(\pi)$.

Any class $\alpha\in\hCH^1_\R(\mathscr X)$,
defines a \emph{height} function~$h_\alpha$ which is a linear
function on the subspace of~$Z^1_\R(\mathscr X)$ 
consisting of real~$1$-cycles~$Z$ on~$\mathscr X$
such that $\omega(\alpha)$ is locally~$\mathrm L^\infty$ on a neighbourhood
of~$\abs Z(\C)$. If $D$ is a  real divisor on~$X$
such that $\omega(\alpha)$ is locally~$\mathrm L^\infty$ in
a neighbourhood of~$\abs D(\C)$,
we shall still denote $h_\alpha(D)$ the height of the
unique horizontal~$1$-cycle on~$\mathscr X$ which extends~$D$.
Moreover, for any effective divisor~$D$ on~$X$
such that~$\omega(\alpha)$ is locally~$\mathrm L^\infty$ in
a neighbourhood of~$\abs{\pi_*(D)}(\C)$, then $\omega(\pi^*\alpha)$
is locally~$\mathrm L^\infty$ in a neighbourhood of~$\abs D(\C)$
and one has the equality $ h_{\pi^*\alpha}(D)=h_\alpha(\pi_*(D))$.

\begin{definition}
Let~$D$ be a divisor on~$X$.
For each place~$v$ of~$K$, let $\Omega_v$ be an open subset
of~$X_v$ (stable under complex conjugation if $v$ is archimedean). 
One says that the collection~$(\Omega_v)$
is an \emph{adelic tube adapted to~$D$} if the following conditions
are satisfied:
\begin{enumerate}
\item for any ultrametric place~$v$,
the complement of~$\Omega_v$ in any connected component of~$X_v$
is a nonempty affinoid subset ;
\item for any archimedean place~$v$,
the complement of~$\Omega_v$ in any connected component of~$X_v$
is non-polar;
\item there exists an effective reduced divisor~$E$ containing~$\abs D$, 
a finite set of places~$F$ of~$K$,
a normal projective flat model~$\mathscr X$ of~$X$ over~$R$
such that for any ultrametric place~$v$ of~$K$
such that $v\not\in F$, $\Omega_v=\tube{\overline E}_v$
is the tube in~$X_v$ around the specialization of~$E$ in
the special fibre~$\mathscr X_{\F_v}$.
\end{enumerate}
% If moreover the divisor~$E$ in~(3) can be taken to be~$D$ itself, we
% shall say that $\Omega$ is \emph{strictly adapted to~$D$}.
\end{definition}

Let $\Omega=(\Omega_v)$ be a family, where, for each place~$v$ of~$K$,
$\Omega_v$ is an open subset of the analytic curve~$X_v$
satisfying the conditions~(1)  and~(2).
Let~$D$ be a divisor on~$X$ whose support is contained in~$\Omega_v$
for any place~$v$ of~$K$.
By the considerations of this section, the line bundle
$\mathscr O_X(D)$ is then endowed, for each place~$v$ of~$K$,
of a~$v$-adic metric~$\norm{\cdot}^\CAP_{\Omega_v}$.
If $\Omega$ is an adelic tube adapted to~$D$, then,
for almost all places of~$K$, this metric
is in fact induced by the horizontal extension of the divisor~$D$
in an adequate model~$\mathscr X$ of~$X$.
Actually, one has the following proposition:
\begin{proposition}
Assume that~$\Omega$ is an adelic tube adapted to~$\abs D$.
There is a normal, flat, projective model~$\mathscr X$
of~$X$ over~$R$ and a (unique) Arakelov $\Q$-divisor extending~$D$,
inducing at any place~$v$ of~$K$
the~$v$-adic capacitary metric on~$\mathscr O_X(D)$.
\end{proposition}
Such an arithmetic surface~$\mathscr X$
will be said \emph{adapted to~$\Omega$}.
Then, the Arakelov $\Q$-divisor on~$\mathscr X$
whose existence is asserted   by
the proposition will be denoted~$\widehat D_{\Omega}$.
Observe moreover that the current $\omega(\widehat D_\Omega)$
is locall~$\mathrm L^\infty$ on~$\Omega$ since it vanishes there.
Consequently, the height $h_{\widehat D_\Omega}(E)$ is defined
when $E$ is any $0$-cycle on~$X$ which is supported
by~$\Omega$.

\begin{proof}
It has already been recalled that archimedean Green
functions defined by potential theory have the
required $\mathrm L_1^2$-regularity.
It thus remains to show that the metrics at finite places
can be defined using a single model~$(\mathscr X,\mathscr D)$ of~$(X,D)$ 
over~$R$.

\begin{lemma}
There exists 
a normal, flat projective model $\mathscr X$ of~$X$ over~$R$,
and, for any ultrametric place~$v$ of~$K$,
a Zariski closed subset~$\mathsf Z_v$ of 
the special fibre~$\mathsf X_{\F_v}$ at~$v$
such that $\Omega_v=\tube{\mathsf Z_v}$.
We may moreover assume that for almost all ultrametric places~$v$ of~$K$,
$\mathsf Z_v=\mathscr E\cap\mathsf X_{\F_v}$,
where $\mathscr E$ is an effective reduced horizontal divisor on~$\mathscr X$.
\end{lemma}
\begin{proof}
Let~$\mathscr X_1$ be a projective flat model of~$X$ over~$R$,
$E$ an effective reduced divisor on~$\mathscr X$,
and~$F$ be a finite set of places 
satisfying the condition~3) of the definition of an adelic tube.
Up to enlarging~$F$, we may assume that the fibre product
$\mathscr X_1\otimes_R R_1$ is normal,
where $R_1$ denotes the subring of~$K$ obtained from~$R$
by localizing outside places in~$F$.

By Raynaud's formal/rigid geometry comparison theorem,
there is, for each finite place~$v\in F$,
a normal projective and flat model~$\mathscr X_v$
of~$X$ over the completion $\widehat{R_v}$,
a Zariski closed subset~$\mathsf Z_v$ the special fibre of~$\mathscr X_v$,
such that $\Omega_v=\tube{\mathsf Z_v}$.

By a general descent theorem of Moret-Bailly
(\cite[Theorem~1.1]{moret-bailly1996}; 
see also~\cite[6.2, Lemma~D]{bosch-l-r90}), 
there exists a projective and flat
$R$-scheme~$\mathscr X$ which coincides with~$\mathscr X_1$
over~$\Spec R_1$ and such that
its completion at any finite place~$v\in F$ is isomorphic
to~$\mathscr X_v$. By faithfully
flat descent, such a scheme is normal
(see~\cite[21.E, Corollary]{matsumura1980}).

For any ultrametric place~$v$ over~$\Spec R_1$, we just let~$\mathsf Z_v$
be the specialization of~$E$ in~$\mathscr X_{\F_v}=(\mathscr X_0)_{\F_v}$;
one has $\Omega_v=\tube{\mathsf Z_v}$ by assumption since $v$ 
does not belong to the finite set~$F$ of excluded places.
For any ultrametric place~$v\in F$, $\mathsf Z_v$ identifies with
a Zariski closed subset of the special fibre~$\mathscr X_{\F_v}$
and its tube is equal to~$\Omega_v$ by construction.
This concludes the proof of the lemma.
\end{proof}

Fix such a model~$\mathscr X$ and
let $\mathscr D_0$ be the Zariski closure of~$D$ in~$\mathscr X$.
For any ultrametric place~$v$ of~$F$, let~$V_v$
be the unique divisor on the special fibre~$\mathscr X_{\F_v}$
such that $\mathscr D_0+V_v$ satisfies the assumptions
of Proposition~\ref{prop.metric}.
One has $V_v=0$ for any ultrametric place~$v$ such that~$\mathsf Z_v$ has 
no component of dimension~$1$, hence for all but finitely places~$v$.
We thus may consider the $\Q$-divisor $\mathscr D=\mathscr D_0+\sum_v V_v$
on~$\mathscr X$ and observe that it induces the capacitary metric
at all ultrametric places.
\end{proof}

\begin{proposition}
Let~$D$ be a divisor on~$X$ and let~$\Omega$
be an adelic tube adapted to~$\abs D$.
One has the equality
\[ \widehat D_\Omega\cdot\widehat D_\Omega=
        h_{\widehat D_\Omega}(D).\]
\end{proposition}
%%  More generally, for two divisors~$D_1$, $D_2$
%%  with adelic tubes~$\Omega_1$, $\Omega_2$,
%%  one has 
%%  \[ \widehat D_1_{\Omega_1}\cdot\widehat D_2_{\Omega_2}
%%  =h_{\widehat D_1_{\Omega_2})(D_2)\]
%%  if ``$\partial\Omega_1$ does not meet~$\Omega_2$''.
%%  (The conditions $\Omega_{2,v}\subset\Omega_{1,v}$ of~\cite{bost1999}
%%  at almost finite places imply $D_2<D_1$, 
%%  hence would be inadequate here.)
\begin{proof}
Let us consider a model~$\mathscr X$ of~$X$
and an Arakelov~$\Q$-divisor~$\mathscr D$ on~$\mathscr X$
defining the capacitary metric~$\norm{\cdot}^\CAP_{\Omega_v}$
at all ultrametric places~$v$ of~$K$.

Let~$\mathscr D_0$ denote the Zariski closure of~$D$
in~$\mathscr X$. For any ultrametric place~$v$ of~$K$,
and let~$V_v$ be the vertical part of~$\mathscr D$ lying
above~$v$  so that $\mathscr D=\mathscr D_0+\sum_v V_v$.
% For any complex embedding~$\sigma\colon K\hra\C$,
% let $g_\sigma$ denote the Green function for~$D$ at the place
% defined by~$\sigma$.
By~\cite[Corollary~5.4]{bost1999}, one has
\[ \widehat D_\Omega\cdot \widehat D_\Omega
 = h_{\widehat D_\Omega}(\mathscr D). \]
By the definition of the capacitary metric at ultrametric places,
the geometric intersection number of~$\mathscr D$ with
any vertical component of~$\mathscr D$ is zero. 
Consequently, 
\[ \widehat D_\Omega\cdot \widehat D_\Omega
 = h_{\widehat D_\Omega}(\mathscr D_0)
+  \sum_v h_{\widehat D_\Omega}(V_v)
= h_{\widehat D_\Omega}(\mathscr D_0), \]
as was to be shown.
\end{proof}

\begin{corollary}\label{coro.height.capa}
Let $P\in X(K) $ be a rational point of~$X$
and let~$\Omega$ be an adelic tube adapted to~$P$.
One has 
\[ \widehat P_\Omega\cdot \widehat P_\Omega=\hdeg 
(T_PX,\norm{\cdot}^\CAP_\Omega).\]
% In particular, $\widehat P_\Omega$ is a numerically
% effective Arakelov $\Q$-divisor if and only if
% $ \hdeg (T_PX,(\norm{\cdot}^\CAP_{\Omega_v})_v)\geq 0$.
\end{corollary}

\section{An algebraicity criterion for~$A$-analytic curves}
\label{sec.algebraic}

Let~$K$ be a number field, $R$ its ring of integers, 
$X$ a quasi-projective algebraic variety over~$K$
and let~$P$ a point in~$X(K)$. 
Let $\widehat C\hra \widehat X_P$
be a smooth formal curve that is~$A$-analytic.

For any place~$v$ of~$K$, 
the formal curve $\widehat C$ is $K_v$-analytic,
and we may equip the~$K$-line $T_P\widehat C$
with the canonical~$v$-adic semi-norm 
$\norm{\cdot}_v^\can=\norm{\cdot}^\can_{X,\widehat C,P,v}$ constructed in Section~\ref{subsec.canonical}.
We claim that, equipped with these semi-norms,  $T_P\widehat C$ defines a semi-normed~$K$-line
$(T_P\widehat C,\norm{\cdot}^\can)$ with a well-defined Arakelov degree 
in~$\mathopen]-\infty, +\infty]$, in the sense of~\cite[4.2]{bost2004}. 
Recall that it means that, for any (or equivalently, for some)
 non-zero element in~$T_P\widehat C$, the series $\sum_v \log^+\norm{t}^\can_v$ is convergent. To see this, consider 
 a quasi-projective
flat~$R$-scheme~$\mathscr X$ with generic fibre~$X$, together with
a section $\mathscr P\colon \Spec R\ra\mathscr X$ which extends~$P$.
According to Lemma~\ref{lemm.can/size}
(applied to projective compactifications of~$X$ and $\cX$, and an ample line bundle $\cL$), 
the inequality
\[  \log \norm{t}^\can_v \leq - \log  S_{\mathscr X,v}(\widehat C), \]
holds for almost all finite places~$v$,
where $S_{\mathscr X,v}$ denotes the size of~$\widehat C$
with respect to the $R_v$ model $\mathscr X\otimes R_v$.
Since by definition of~$A$-analyticity
the series with non-negative terms $\sum_v\log S_{\mathscr X,v}(\widehat C)^{-1}$ has a finite sum, this establishes the required convergence.

The Arakelov degree of  $(T_P\widehat C, \norm{\cdot}^\can)$ is defined as the sum:
\[ \hdeg (T_P\widehat C, \norm{\cdot}^\can):= \sum_v (-\log\norm{t}^\can_v).\]
It is a well defined element in~$\mathopen]-\infty,+\infty]$, independent of the choice of~$t$ by the product formula (we follow the usual convention $-\log 0= + \infty$.)

The following criterion extends  Theorem~4.2 of~\cite{bost2004},
where instead of canonical semi-norms, larger norms constructed by means of the sizes were used at finite places.

\begin{theorem}\label{theo.algebraic} Let $\widehat C$ be, as above, an~$A$-analytic curve through a rational point~$P$ in some algebraic variety~$X$ over~$K$.

If $\hdeg (T_P\widehat C,\norm{\cdot}^\can)>0$, then $\widehat C$
is algebraic.
\end{theorem}

\begin{proof} We keep the above notation, and we assume, as we may, $X$ (resp. $\cX$) to be projective over~$K$ (resp. over~$R$). We choose an ample line bundle $\cL$ over $\cX$ and we let $L:=\cL_K$.

We let $\cE_D:=\Gamma(\cX, \cL^{\otimes D})$ and, 
for any embedding $\sigma:K\hookrightarrow \C$, we choose 
a consistent sequence of hermitian norms~$(\norm{\cdot}_{D,\sigma})$ 
on the $\C$-vector spaces 
$\cE_{D,\sigma} \simeq \Gamma (X_\sigma, L_\sigma^{\otimes D})$, 
in a way compatible with complex conjugation. 
Using these norms, we define hermitian vector bundles
$\overline{\cE}_D:= (\cE_D,(\norm{\cdot}_{D,\sigma})_\sigma)$ over~$\Spec R$.

We also choose an hermitian structure on~$\cP^*\cL$, and we denote
$\overline{\cP^*\cL}$ the so-defined hermitian line bundle over $\Spec
R$. Finally, we equip $T_P\widehat C$ with the~$R$-structure defined
by $N_{\mathscr P} \mathscr X \cap T_P\widehat C$ and with an arbitrary
hermitian structure, and in this way we define an hermitian line bundle
$\overline T_0$ over $\Spec R$ such that $(T_0)_K=T_P\widehat C$.

We define the~$K$-vector spaces $E_D:=\cE_{D,K} \simeq \Gamma( X,
L^{\otimes D})$, their subspaces~$E^i_D$, and the evaluation maps
\[ \phi^i_D \colon E^i_D \ra (T^\vee_P\widehat C)^{\otimes i}
           \otimes L_{\mid P}^{\otimes D}. \]
as in the ``local'' situation considered in Section~\ref{subsec.canonical}.
According to the basic algebraicity criteria in~\cite[2.2]{bost2004},
to prove that $\widehat C$ is algebraic, it suffices to prove that the
ratio
\begin{equation} \label{ratio}
\frac{ \sum\limits_{i\geq 0} (i/D)
\rank (E^i_D/E^{i+1}_D)}
         {\sum\limits_{i\geq 0}  \rank (E^i_D/E^{i+1}_D)}
\end{equation}
stays bounded when~$D$ goes to~$+\infty$.

For any place~$v$ of~$K$, the morphism $\phi^i_D$
has a~$v$-adic norm, defined by means of the integral and hermitian 
structures introduced above. If $\phi^i_D\neq 0$, the height of~$\phi^i_D$
is the real number defined as the (finite) sum:
\[ h(\phi^i_D) = \sum_v \log \norm{\phi^i_D}_v. \]
When $\phi^i_D$ vanishes, we define $h(\phi^i_D)=-\infty$; observe that,
in this case, $E^{i+1}_D=E^i_D$.

As established in the proof of Lemma~\ref{lemm.can/size} above
(see also \cite[Lemma~3.3]{bost2001}), the following inequality holds for any finite
place~$v$, and any two non-negative integers~$i$ and~$D$:
\begin{equation}
\label{nonarchCauchy}
\log \norm{\phi^i_D}_v\leq - i\log S_{\mathscr X,v}(\widehat C). 
\end{equation}
Since $\widehat C$ is $A$-analytic,
the upper bounds~(\ref{archCauchy}) and~(\ref{nonarchCauchy}) 
show the existence of some positive real number $c$ such that
\begin{equation}
\label{linearbound}
h(\phi^i_D) \leq c (i+D).
\end{equation}

For any place~$v$ of~$K$, we let
\[ \rho_v(L)=\limsup_{i/D\ra\infty} \frac1i \log\norm{\phi^i_D}_v.\]
This is an element in~$[-\infty, +\infty\mathclose[$, which, according to (\ref{nonarchCauchy}), satisfies:
\[ \rho_v(L) \leq - \log  S_{\mathscr X,v}(\widehat C)\] 
for any finite place~$v$. Moreover, by its very definition, the Arakelov degree of $(T_P\widehat C,\norm{\cdot}^\can)$ is given by:
%\begin{equation}
\begin{align*}
 \hdeg (T_P\widehat C,\norm{\cdot}^\can) &
= \sum_{v} (-\rho_{v}(L)) + \hdeg \overline{T}_{0} \\
& =\sum_{v\, \text{finite}} \left(-\rho_{v}(L) - \log  S_{\mathscr X,v}(\widehat C)\right) \\
& \qquad\qquad {} + \sum_{v\, \text{finite}}  \log  S_{\mathscr X,v}(\widehat C) + \sum_{v \mid \infty} (-\rho_{v}(L))
 + \hdeg \overline{T}_{0}.
\end{align*}
%\end{equation}
In the last expression, the terms of the first sum belong to
$[0,+\infty]$ --- and the sum itself is therefore well-defined in 
$[0,+\infty]$ --- and the second sum is convergent by~$A$-analyticity of 
$\widehat{C}$.

Observe also that, since the sums
\[ \sum_{v\,{\rm finite}} \left(-\frac1i \log \norm{\phi^i_D}_v + \log  S_{\mathscr X,v}(\widehat
C)\right)\]
have non-negative terms, we get, as a special instance of Fatou's 
Lemma:
\begin{multline*}
 \sum_{v\, {\rm finite}} \liminf_{i/D\ra\infty}\left(-\frac1i \log \norm{\phi^i_D}_v + \log  S_{\mathscr X,v}(\widehat
C)\right)
\\ \leq 
\liminf_{i/D\ra\infty}
\sum_{v\, {\rm finite}} \left(-\frac1i \log \norm{\phi^i_D}_v + \log  S_{\mathscr X,v}(\widehat
C)\right).\end{multline*}
Consequently
\[  \limsup_{i/D\ra\infty} \frac1i h(\phi^i_D)
\leq \sum_v \rho_v(L),\]
and
\begin{equation}\label{etdeun}
\hdeg(T_P\widehat C,\norm{\cdot}^\can)\leq  
-\limsup_{i/D\ra\infty} \frac1i h(\phi^i_D) + \hdeg \overline{T}_{0}.
\end{equation}
When $ \hdeg (T_P\widehat C, \norm{\cdot}^\can)$ is positive, the
inequality (\ref{etdeun}) implies the existence of positive real
numbers $\epsilon$ and $\lambda$ such that, 
for any two positive integers~$i$ and~$D$, 
\begin{equation}\label{etdequatre} 
    \hdeg \overline{T}_{0} - \frac1i h(\phi^i_D) \geq
    \epsilon \quad\text{if $i\geq\lambda D$.} 
\end{equation}

Let $\cE^i_D:=\cE_D \cap E^i_D$ and let $\overline{\cE^i_D/\cE^{i+1}_D}$ 
be the hermitian vector bundle on 
$\Spec R$ defined by the quotient  $\cE^i_D/\cE^{i+1}_D$ equipped 
with the hermitian structure induced by the one of~$\overline{\cE}_D$. 
The evaluation map~$\phi^i_D$ induces 
an injection $E^i_D/E^{i+1}_D \hookrightarrow  (T^\vee_P\widehat C)^{\otimes i}
\otimes L_{\mid P}^{\otimes D}$. Actually, either $\phi^i_{D}=0$ and
then $E^i_D=E^{i+1}_D$, or $\phi^i_{D}\neq 0$, and this inclusion is
an isomorphism of~$K$-lines. In either case, we have:
\[ \hdeg \overline{\cE^i_D/\cE^{i+1}_D}
= \rank (E^i_D/E^{i+1}_D) \left(\hdeg (\overline{\cP^*\cL}^{\otimes D}
\otimes \overline{T}_{0}^{\vee\otimes i}) + h(\phi^i_{D})\right).\]  
Indeed, if $\phi^i_{D}=0$, both sides vanish (we follow the usual
convention $0\cdot (-\infty)=0$). If $\phi^i_{D}\neq 0$, the equality is a
straightforward consequence of the definitions of the Arakelov degree
of an hermitian line bundle over $\Spec R$ and of the heights
$h(\phi^i_{D})$.

The above equality may also be written:
\begin{equation}\label{etdedeux}\hdeg \overline{\cE^i_D/\cE^{i+1}_D}
= \rank (E^i_D/E^{i+1}_D) \left(D\hdeg \overline{\cP^*\cL}
-i\hdeg \overline{T}_{0} + h(\phi^i_{D})\right).
\end{equation}

Moreover, 
by~\cite[Proposition~4.4]{bost2001}, there is a constant~$c'$,
such that for any $D\geq 0$ and any saturated submodule $\cF$ of
$\cE_{D}$,
\[ \hdeg \overline{\cE_{D}/\cF} \geq - c' D\rank(\cE_{D}/\cF).\]
(This is an easy consequence of the fact that the~$K$-algebra 
$\bigoplus_{D\geq 0}\cE_{D,K}$ is finitely generated.)
Applied to $\cF:=\bigcap_{i\geq 0}\cE^i_{D}$, this estimate becomes:
\begin{equation}\label{etdetrois}
\sum_{i\geq 0} \hdeg \overline{\cE^i_{D}/\cE^{i+1}_{D}}
\geq -c'D\sum_{i\geq 0} \rank (E^i_D/E^{i+1}_D).
\end{equation}
Using (\ref{etdedeux}) and (\ref{etdetrois}), we derive the inequality:
\begin{multline}
\label{fourstars}
-(c'+\hdeg \overline{\cP^*\cL})D\sum_{i\geq 0} \rank (E^i_D/E^{i+1}_D) \\
\leq
\sum_{i\geq 0} \rank (E^i_D/E^{i+1}_D) (-i\hdeg
\overline{T}_{0} + h(\phi^i_{D})).
\end{multline}

Finally, using (\ref{fourstars}), (\ref{linearbound}), and (\ref{etdequatre}),
we obtain
\begin{multline*}
\sum_{i< \lambda D} \rank (E^i_D/E^{i+1}_D)(\frac{i}{D} \hdeg \overline{T}_{0} -c \frac{i+D}{D}) 
+
\sum_{i\geq \lambda D}  \rank (E^i_D/E^{i+1}_D) \epsilon\frac{i}{D} \\
\leq
(c'+\hdeg \overline{\cP^*\cL})
\sum_{i\geq 0}  \rank (E^i_D/E^{i+1}_D).
\end{multline*}
This implies that the ratio (\ref{ratio}) is bounded by
\[ \lambda + \frac1\epsilon\left( c'+ \hdeg \overline{\cP^*\cL}+ c + \lambda \max(0,c-\hdeg \overline{T}_{0} )\right),\]
and completes the proof.
\end{proof}

\section{Rationality criteria}
\label{sec.rat}
\subsection{Numerical equivalence and numerical effectivity
on arithmetic surfaces}

The following results are variations on a
classical theme in Arakelov geometry of arithmetic surfaces.
The first theorem characterizes numerically
trivial Arakelov divisors with real coefficients.
It is used in the next proposition to describe
effective Arakelov divisors whose sum is numerically effective.
We allow ourselves to use freely the notations of~\cite{bost1999}.

\begin{theorem}[{Compare~\cite[Theorem~5.5]{bost1999}}]
\label{thm5.5}
Let $\mathscr X$ be a normal flat projective scheme over
the ring of integers of a number field~$K$ whose
generic fibre is a smooth and geomerically connected curve.
Let $(D,g)$ be any element in~$\hZ^1_\R(\mathscr X)$
which is numerically trivial. Then there exist
an integer~$n$, real numbers $\lambda_i$ and rational functions
$f_i\in K(\mathscr X)^*$, for $1\leq i\leq n$,
and a family $(c_\sigma)_{\sigma\colon K\hra\C}$ of real numbers
such that $c_{\bar\sigma}=c_{\sigma}$, $\sum c_\sigma=0$, and
$(D,g)=(0,(c_\sigma))+\sum_{i=1}^n\lambda_i \hdiv(f_i)$.
\end{theorem}
\begin{proof}
There are real numbers $\lambda_i$
and Arakelov divisors $(D_i,g_i)\in\hZ^1(\mathscr X)$ 
such that $(D,g)=\sum \lambda_i (D_i,g_i)$.
We may assume that the~$\lambda_i$ are linearly independent over~$\Q$.
By assumption, the degree of~$D$ on any vertical component
of~$\mathscr X$ is zero; the linear independence of the~$\lambda_i$
implies that the same holds
for any~$D_i$. Let us then denote by~$g'_i$ any Green current
for~$D_i$ such that $\omega(D_i,g'_i)=0$. One has 
\[ 0=\omega(D,g)=\sum \lambda_i\omega(D_i,g_i) = \sum
\lambda_i\omega(D_i,g'_i), \]
so that the difference $g-\sum \lambda_i g'_i$ is harmonic,
and therefore constant on any connected component of~$\mathscr X(\C)$. 
By adding a locally constant function to some~$g'_i$, we may assume
that $g=\sum\lambda_i g'_i$. Then, $(D,g)=\sum\lambda_i (D_i,g'_i)$.
This shows that we may assume that one has $\omega(D_i,g_i)=0$ for any~$i$.
By Faltings-Hriljac's formula, 
the N\'eron-Tate quadratic form on~$\Pic^0(\mathscr X_K)\otimes\R$
takes the value~$0$ on the class of the real 
divisor~$\sum \lambda_i(D_i)_K$. Since
this quadratic form is
positive definite (see~\cite[3.8, p.~42]{serre1997}),
this class is zero.
Using that the~$\lambda_i$ are linearly independent over~$\Q$, we deduce that
the class of each divisor~$(D_i)_K$ in~$\Pic^0(\mathscr X_K)$
is torsion.
Since $D_i$ has degree zero 
on any vertical component of~$\mathscr X$ and the Picard
group of the ring of integers of~$K$ is finite, the class
in~$\Pic(\mathscr X)$ of the divisor~$D_i$ is torsion too.
Let then choose positive integers~$n_i$ and rational
functions~$f_i$ on~$\mathscr X$ such that $\div(f_i)=n_i D_i$.
The Arakelov divisors $\hdiv(f_i)-n_i(D_i,g_i)$
are of the form $(0,c_i)$, where $c_i=(c_{i,\sigma})_{\sigma\colon
K\hra\C}$ is a family of real numbers such that 
$c_{i,\bar\sigma}=c_{i,\sigma}$ and $\sum_\sigma c_{i,\sigma}=0$.
Then, letting $c_\sigma=\sum_i (\lambda_i/n_i)c_{i,\sigma}$, one has 
\[ (D,g)= (0,(c_\sigma))+ \sum \frac{\lambda_i}{n_i} \hdiv(f_i)
\]
as requested.
\end{proof}

Let $f_1,\dots,f_n$ be meromorphic functions on some Rieman surface~$M$,
and $\lambda_1,\dots,\lambda_n$  real numbers, and
let $f\in \C(M)^*\otimes_\Z\R$ be defined as 
$f=\sum_{i=1}^n f_i\otimes \lambda_i$. We shall denote 
$\abs f$ the real function on~$M$ given by $\prod \abs{f_i}^{\lambda_i}$,
and by~$\div f$ the $\R$-divisor $\sum\lambda_i \div(f_i)$;
they don't depend on the decomposition of~$f$ as a sum
of tensors. 
One has $\ddc  \log\abs{f}^{-2}+\delta_{\div(f)}=0$.

We shall say that a pair $(D,g)$ formed of a divisor~$D$ on~$M$
and of a Green current~$g$ with $L_1^2$ regularity for~$D$
is \emph{effective}\footnote{in the terminology of~\cite{bost1999},
\emph{nonnegative}}
if the divisor~$D$ is effective,
and if the Green current~$g$ of degree~$0$ for~$D$
may be represented by a nonnegative summable function
(see~\cite[Definition~6.1]{bost1999}).
% we say that $(D,g)$ is supereffective if 
% moreover the current $\omega(g)=\ddc g+\delta_D$ is a positive measure.

%%  Let $\mathscr X$ be an arithmetic surface, and let us denote
%%  by $\pi\colon\mathscr X\ra\Spec\mathfrak o_K$ the Stein
%%  factorization of the canonical morphism $\mathscr X\ra\Spec\Z$.
Similarly,
we say that an Arakelov divisor $(D,g)\in\hZ^1_\R(\mathscr X)$ on
the arithmetic surface~$\mathscr X$ is effective
if~$D$ is effective on~$\mathscr X$ and if
$(D_\C,g)$ is effective on~$\mathscr X(\C)$.

We say that an Arakelov divisor, or the class~$\alpha$
of an Arakelov divisor,
is \emph{numerically effective} (or shortly, \emph{nef})
if $[(D,g)]\cdot\alpha\geq 0$ for any effective Arakelov divisor 
$(D,g)\in\hZ^1_\R(\mathscr X)$ 
(according to~\cite[Lemma~6.6]{bost1999},
it is sufficient to consider Arakelov divisors $(D,g)$
with $\mathscr C^\infty$-regularity).
If $(D,g)$ is an effective and numerically effective Arakelov divisor,
then the current $\omega(g)\mathrel{:=}\ddc g+\delta_D$ 
is a positive measure
(see~\cite[proof of Proposition~6.9]{bost1999}).

\begin{proposition}\label{prop.lefschetz}
Let $\mathscr X$ be a normal, flat projective scheme
over the ring of integers of a number field~$K$
whose generic fibre is a smooth geometrically connected algebraic curve.

Let $(D,g)$ and $(E,h)$ be non-zero elements of~$\hZ^1_\R(\mathscr X)$;
let $\alpha$ and $\beta$ denote their classes in~$\hCH^1_\R(\mathscr X)$.
Let us assume that the  following conditions are satisfied:
\begin{enumerate}
\item the Arakelov divisors $(D,g)$ and $(E,h)$ are effective;
\item the supports of~$D$ and~$E$ do not meet 
and $\int_{\mathscr X(\C)}g \ast h=0$.
% \item the form $\omega(D,g)=\ddc g+\delta_D$ (\resp the form $\omega(E,h)$)
% is locally bounded near the support of~$D$ (\resp of~$E$).
\end{enumerate}
If the class $\alpha+\beta$ is numerically effective,
then
there exist a positive real number~$\lambda$, 
an element $f\in K(\mathscr X)^*\otimes_\Z\R$ 
and a family $(c_\sigma)_{\sigma\colon K\hra\C}$
of real numbers which is invariant by conjugation and
satisfies $\sum_\sigma c_\sigma=0$, such that for any embedding
$\sigma\colon K\hra\C$, 
\[ g_\sigma = (c_\sigma+\log \abs f^{-2})^+, \quad\text{and}\quad
  h_\sigma = \lambda (c_\sigma+\log \abs f^{-2})^-, \]
where, for any real valued function~$\phi$,
we denoted $\phi^+=\max(0,\phi)$ and $\phi^-=\max(0,-\phi)$, so that
$\phi^+-\phi^-=\phi$.

Moreover, $\alpha^2=\alpha\beta=\beta^2=0$.
\end{proposition}

\begin{proof}
Since $(D,g)$ and~$(E,h)$ are effective and non-zero,
the classes $\alpha$ and~$\beta$ are not equal to zero
\cite[Proposition~6.10]{bost1999}.
Moreover, the assumptions of the proposition imply that
\[ \alpha\cdot\beta = \deg\pi_* (D,E) + \frac12 \int_{\mathscr X(\C)} g\ast
h = 0. \]
Since $\alpha+\beta$ is numerically effective,
it follows from Lemma~6.11 of~\cite{bost1999}
(which in turn is an application of the Hodge index theorem
in Arakelov geometry)
that there exists $\lambda\in\R_+^*$
such that $\beta=\lambda\alpha$ in~$\hCH^1_\R(\mathscr X)$.
In particular, $\alpha$ and~$\beta$ are nef, and
$\alpha^2=\beta^2=\alpha\cdot\beta=0$.

Replacing $(E,h)$ by~$(\lambda E,\lambda h)$,
we may assume that $\lambda=1$. 
Then, $(D-E,g-h)$ belongs to the kernel of
the canonical map 
$\rho\colon \hZ^1_\R(\mathscr X)\ra\hCH^1_\R(\mathscr X)$, 
so is numerically trivial.
By Theorem~\ref{thm5.5}, there exist real numbers $\lambda_i$,
rational functions~$f_i\in K(\mathscr X)^*$ 
and a
family $c=(c_\sigma)_{\sigma\colon K\hra\C}$ of real numbers,
invariant by conjugation, 
such that $\sum_{\sigma} c_\sigma=0$
and $(D-E,g-h)=(0,c)+\sum\lambda_i \hdiv(f_i)$ in~$\hZ^1_\R(\mathscr X)$.
Let us denote by~$f$ the element $\sum f_i\otimes \lambda_i$
of~$K(\mathscr X)^*\otimes_\Z\R$. 
The proposition now follows by applying
Lemma~\ref{lemm.g*h=0} below
to the connected Riemann surface $\mathscr X_\sigma(\C)$, 
the pairs $(D,g_\sigma)$, $(E,h_\sigma)$ and
the ``meromorphic function'' $e^{-2c_\sigma}f_{\mathscr X_\sigma(\C)}$,
for each embedding $\sigma\colon K\hra\C$.
\end{proof}

\begin{lemma}\label{lemm.g*h=0}
Let~$M$ be a compact connected Riemann surface, let~$D$ and $D'$  two
nonzero $\R$-divisors on~$M$, and let~$g$ and $g'$
be two Green functions with $L_1^2$ regularity for~$D$
and $D'$.
We make the following assumptions: $\abs D\cap \abs{D'}=\emptyset$,
the pairs $(D,g)$ and $(D',g')$
are effective, the currents $\omega(g)=\ddc g+\delta_D$
and $\omega(g')=\ddc g'+\delta_{D'}$ are positive
measures, $\int_{M} g*g'=0$.
If there exists an element $f\in \C(M)^*\otimes\R$
such that $g-g'=\log\abs{f}^{-2}$,
then $g=\max(0,\log\abs{f}^{-2})$ and
$g'=\max(0,\log\abs f^2)$.
\end{lemma}
\begin{proof}
First observe that 
\[ \omega(g)-\omega(g')=\ddc (g-g')+\delta_D-\delta_{D'} 
= \ddc \log\abs{f}^{-2}+\delta_D-\delta_{D'}=\delta_{D-D'-\div(f)}, \]
by the Poincar\'e-Lelong formula. By assumption, the current
$\omega(g)-\omega(g')$ belongs to the Sobolev space~$L_{-1}^2$; it is therefore 
non-atomic (see~\cite[Appendix, A.3.1]{bost1999}),
so that $D-D'=\div(f)$ and $\omega(g)=\omega(g')$.

Observe also that $g_{\mid M \setminus \abs D}$ (resp. $g'_{\mid M \setminus \abs D'}$) is a subharmonic current. In the sequel, we denote by $g$ (resp. $g'$) the unique subharmonic function on $M \subset \abs D$ (resp. on $M' \subset \abs D'$) which represents this current.

Let~$F$ be the set of points $x\in M$ where $\abs{f(x)}=1$
and let $\Omega=M\setminus F$ be its complementary subset.
The functions
$h=\max(0,\log\abs f^{-2})$ and $h'=\max(0,\log\abs f^2)$
are continuous Green functions with $L_1^2$ regularity
for~$D$ and~$D'$ respectively. The currents
$\ddc  h+\delta_D$, $\ddc h'+\delta_{D'}$ are equal
to a common positive measure, which we denote by~$\nu$.
Since $h$ (resp. $h'$) is harmonic on $M\setminus(\abs D \cup F)$ (resp. on
$M\setminus(\abs D' \cup F)$), this measure is supported by $F.$
%The classical formula  
%$\ddc \max(0,\log\abs z^2)=\frac1{2\pi}d\Arg z\wedge \delta_{\abs z=1}$ 
%implies that one has 
%$\nu=\frac1{2\pi}d\Arg f\wedge\delta_{\abs f=1}$. In particular
%the measure~$\nu$ is positive and supported by~$F$.

Let~$S$ be the support of the positive measure~$\omega(g)$.
It follows from~\cite[Remark~6.5]{bost1999} that
$g$ and $g'$ vanish $\omega(g)$-almost everywhere on~$M$.
Consequently, the equality $\log\abs{f}^{-2}=g-g'=0$ holds
$\omega(g)$-almost everywhere;
in particular, $S\subset F$.

Let us pose $u=h-g=h'-g'$; this is a current with $L_1^2$
regularity on~$M$ and $\ddc u=\ddc  h -\ddc  g = \nu-\omega(g)$.
In particular, $\ddc (u|_\Omega)=0$: $u$ is harmonic on~$\Omega$.
Since~$g$ is nonnegative, one has $u\leq 0$ on~$F=\complement\Omega$.
By the maximum principle, this implies that $u\leq 0$ on~$\Omega$ (\emph{cf.}~\cite[Theorem A.6.1]{bost1999}; observe that $u$ is finely continuous on $M$).

Finally, one has 
\[ 0=\int_M g*g'=\int_M h*h' - \int_M u\nu - \int_M u\omega(g)
\geq \int_M h*h' . \]
% \[ 0=\int_M g*g'=\int_M (h+\alpha)*g'
% = \int_M h*g' + \int_M \alpha *g'
% = \int_M h*(h'+\alpha)+ \int_M \alpha\omega  \\
% =\int_M h*h'+\int_M \alpha\omega + \int_M \alpha\nu. \]
By~\cite[Corollary~6.4]{bost1999}, this last term is nonnegative,
so that all terms of the formula vanish.
In particular, $\int u\nu=0$, hence $u=0$ ($\nu$-a.e.).
Using again that~$u$ is harmonic on~$\Omega$, it follows that its Dirichlet norm vanishes, and finally that $u\equiv 0$.
\end{proof}

\begin{remark}\label{rema.lefschetz.green}
The Green currents~$g$ and~$h$ appearing in
the conclusion of Proposition~\ref{prop.lefschetz} are very particular.
Assume for example that  the Arakelov divisors $\widehat D$ 
and~$\widehat E$
are defined using capacity theory at the place~$\sigma$,
with respect to an open subset $\Omega_\sigma$ of~$X_\sigma$.
Then, $g_\sigma$ and~$h_\sigma$ vanish nearly everywhere 
on~$\complement\Omega_\sigma$. In other words,
$\complement \Omega_\sigma$ is contained in the set of~$x\in X_\sigma$
such that $\abs{f(x)}^{2}=\exp(-c_\sigma)$,
which is a real semi-algebraic curve in~$X_\sigma$,
viewed as a real algebraic surface.
In particular, it contradicts any
of the following hypothesis on~$\Omega_\sigma$, 
respectively denoted~(4.2)$_{\mathscr X,\Omega_\sigma}$ 
and (4.3)$_{\mathscr X,\Omega_\sigma}$
in~\cite{bost1999}:
\begin{enumerate}
\item 
the interior of~$\mathscr X_\sigma(\C)\setminus\Omega_\sigma$
is not empty;
\item 
there exists an open subset~$U$ of
\mbox{$\mathscr X_\sigma(\C)\setminus\abs{D}(\C)$}
not contained in~$\Omega$ such that
any harmonic function on~$U$ which vanishes nearly everywhere
on~\mbox{$U\setminus\Omega$} vanishes on~$U$.
\end{enumerate}
\end{remark}

\subsection{Rationality criteria for algebraic and analytic functions on curves over number fields}
\label{sec.rational}

Let~$K$ be a number field and~$X$ be a smooth projective geometrically connected curve
over~$K$.
For any place~$v$ of~$K$, we denote by~$X_v$ the associated
rigid analytic curve over~$K_v$ if~$v$ is ultrametric,
resp.~the corresponding Riemann surface $X_\sigma(\C)$ if~$v$ is induced by
an embedding of~$K$ in~$\C$.

Let~$D$ be an effective divisor in~$X$
and $\Omega=(\Omega_v)_v$ an adelic tube adapted to~$\abs D$.
We choose a normal projective flat model of~$X$
over the ring of integers~$\mathcal O_K$ of $K$, say $\mathscr X$, 
and an Arakelov $\Q$-divisor $\widehat D_\Omega$ 
on~$\mathscr X$ inducing the capacitary metrics~$\norm{\cdot}^\CAP_{\Omega_v}$
at all places~$v$ of~$K$.
In particular, we assume that for any ultrametric place~$v$,
$\Omega_v$ is the tube~$\tube{\mathsf Z_v}$
around a closed Zariski subset~$\mathsf Z_v$
of its special fibre~$\mathscr X_{\F_v}$,
and $\mathsf Z_v=\overline D\cap \mathscr X_{\F_v}$
for almost all places~$v$.

Our first statement in this section is the following
arithmetic analogue of Proposition~\ref{prop.geom-hodge}.
%%  We begin with a definition.
%%  
%%  \begin{definition}
%%  Let~$K$ be a  number field and let~$R$ be the ring of integers of~$X$.
%%  Let~$X$ be a smooth projective
%%  curve over~$K$ and let~$D$ be an effective divisor on~$X$.
%%  An \emph{adelic tube} of~$D$ is the following data :
%%  \begin{enumerate}
%%  \item a model~$\mathscr X$ over~$R$ of~$X$ 
%%  which is a projective normal arithmetic surface ;
%%  \item an effective divisor~$\mathscr D$ on~$\mathscr X$ such
%%  that $\mathscr D_K=D$;
%%  \item an open neighbourhood~$\Omega_\infty$ of~$\abs{D}$ 
%%  in the Riemann surface~$X(\C)$ whose complement is not
%%  polar in any connected component of~$X(\C)$.
%%  \end{enumerate}
%%  \end{definition}
%%  Given such an adelic tube, we can define as follows
%%  an open subset~$\Omega_v$ of the analytic curve $X_v$, 
%%  for each place~$v$ of~$K$.
%%  At archimedean places, observe that $X(\C)$ is the disjoint union the
%%  Riemann surfaces~$X_v$; we let $\Omega_v=\Omega_\infty\cap X_v$.
%%  At a finite place~$v$, we let~$\Omega_v$ be the preimage
%%  of~$\mathscr D_{\F_v}$ by the reduction map; it is the complementary
%%  subset of an affinoid subspace of the rigid $K_v$-analytic curve~$X_v$.

\begin{proposition}\label{prop.arith-hodge}
Let $X'$ be another geometrically connected smooth projective curve over~$K$
and $f\colon X'\ra X$ be a non constant morphism. 
Let~$D'$ be an effective divisor in~$X'$.
We make the following assumptions:
\begin{enumerate}
\item by restriction, $f$ defines an isomorphism 
from the subscheme~$D'$ of~$X'$ to the subscheme~$D$ of~$X$
and is \'etale in a neighbourhood of~$\abs{D'}$;
\item for any place~$v$ of~$K$, the morphism~$f$ admits an analytic section 
$\phi_v\colon \Omega_v\ra X'_v$ defined over~$\Omega_v$ whose
formal germ is equal to $\widehat{f_{D}}_{K_v}^{-1}$;
\item the class of the Arakelov $\Q$-divisor
$\widehat D_\Omega$  is numerically effective.
\end{enumerate}
Assume moreover 
\begin{enumerate}
\item [\textup{4$'$)}] either that $\widehat D_\Omega\cdot\widehat D_\Omega >0$;
\item [\textup{4$''$)}] or that there is an archimedean place~$v$ 
such that the complementary subset to~$\Omega_v$ in~$X_v$ 
is not contained in a real semi-algebraic curve of~$X_v$.
\end{enumerate}
Then~$f$ is an isomorphism.
\end{proposition}

\begin{proof}
Let us denote by~$E$ the divisor~$f^*D$ on~$X'$; we will
prove that $E=D'$. Observe that, according to assumption~1), this divisor may be written
\[ E:=f^*D=D' + R,\]
where $R$ denotes an effective or zero divisor on $X$, whose support is disjoint from the one of $D'$.  

Let~$\mathscr X'$ denote the normalization of~$\mathscr X$
in the function field of~$X'$ and let us still denote by~$f$
the natural map from~$\mathscr X'$ to~$\mathscr X$ which extends~$f$.
Then~$\mathscr X'$ is a normal projective flat model of~$X'$ over~$\mathcal O_K$.
For any place~$v$ of~$K$, let $\Omega'_v$ denote
the preimage~$f^{-1}(\Omega_v)$ of~$\Omega_v$ by~$f$.
The complementary subset of~$\Omega_v$ is
a nonempty affinoid subspace of~$X'_v$ if~$v$ is ultrametric,
and a non-polar compact subset of~$X'_v$ if~$v$ is archimedean.
Moreover, for almost all ultrametric places~$v$,
$\Omega'_v$ is the tube around the specialization
in~$\mathscr X'_{\F_v}$ of~$f^{-1}(D)$.
In particular, the collection~$\Omega'=(\Omega'_v)$ is an adelic
tube adapted to~$\abs E$.

We thus may assume that the capacitary metrics on
% the line bundles 
$\mathscr O_{X'}(D')$
and~$\mathscr O_{X'}(E)$ relative
to the open subsets~$\Omega'_v$ are induced by
Arakelov $\Q$-divisors on~$\mathscr X'$.
Let us denote them by~$\widehat {D'}_{\Omega'}$ and~$\widehat E_{\Omega'}$ 
respectively.

Since~$X$ and~$X'$ are normal, and the associated rigid analytic spaces as well, the image $\phi_v(\Omega_v)$
of~$\Omega_v$ by the analytic section~$\phi_v$ 
is a closed and open subset~$\Omega_v^1$
of~$\Omega'_v$ containing~$\abs{D'}$, and the collection~$\Omega^1=(\Omega_v^1)$
is an adelic tube adapted to~$\abs{D'}$. Consequently, by Proposition~\ref{prop.capa.connected},
one has $\widehat{D'}_{\Omega'}=\widehat{D'}_{\Omega^1}$.
Similarly, denoting $\Omega^2_v=\Omega'_v\setminus\Omega_v^1$,
the collection $\Omega^2=(\Omega^2_v)$ is an adelic tube
adapted to~$\abs R$ and 
$\widehat R_{\Omega'}=\widehat R_{\Omega^2}$.
One has 
$\widehat {E}_{\Omega'}=f^*\widehat D_{\Omega}=\widehat {D'}_{\Omega^1}+\widehat R_{\Omega^2}$.
Since $\Omega^1_v\cap\Omega^2_v=\emptyset$
for any place~$v$, Lemma~\ref{lemma.5.1} below implies that
$[\widehat R_{\Omega^2}]\cdot[\widehat {D'}_{\Omega^1}]=0$.

Since~$\widehat D$ is non-zero and its class is numerically effective,
the class 
in $\hCH^1_{\Q}(\mathscr X')$
of the Arakelov divisor~$f^*\widehat D=\widehat D+\widehat R$
is numerically effective too.
Proposition~\ref{prop.lefschetz} and Remark~\ref{rema.lefschetz.green} show that,
when
either of the hypotheses~(4$'$) or~(4$''$) is satisfied,
necessarily $\widehat R_{\Omega^2}=0$.
In particular, $R=0$ and $E=D'$. It follows that~$f$ has degree one, hence  is an isomorphism.
\end{proof}

\begin{lemma}\label{lemma.5.1}
Let~$X$ be a geometrically connected smooth projective
curve over a number field~$K$,
let $D_1$ and~$D_2$ be divisors on~$X$ and let $\Omega_1$
and~$\Omega_2$ be adelic tubes adapted to~$\abs{D_1}$
and~$\abs{D_2}$.
Let us consider a normal projective and flat model~$\mathscr X$ of~$X$
over the ring of integers of~$K$ as well as Arakelov divisors
$\widehat {D_1}_{\Omega_1}$ and~$\widehat {D_2}_{\Omega_2}$
inducing the capacitary metrics on~$\mathscr O_X(D_1)$
and~$\mathscr O_X(D_2)$ relative to the adelic tubes~$\Omega_1$
and~$\Omega_2$.

If $\Omega_{1,v}\cap\Omega_{2,v}=\emptyset$
for any place~$v$ of~$K$,
then
\[ \widehat{D_1}_{\Omega_1}\cdot \widehat{D_2}_{\Omega_2}=0 .\]
\end{lemma}
\begin{proof}
Observe that $D_1$ and~$D_2$ have no common component, since
any point~$P$ common to~$D_1$ and~$D_2$ would belong
to $\Omega_{1,v}\cap\Omega_{2,v}$.

Let~$\mathscr X$ be a normal projective flat model of~$X$ 
adapted to~$\Omega_1$ and~$\Omega_2$, so that
the classes $\widehat{D_i}_{\Omega_i}$ live in~$\hCH^1_\R(\mathscr X)$.
Namely $\widehat {D_i}_{\Omega_i}=(\mathscr D_i,g_i)$, where
$\mathscr D_i$ is the $\Q$-divisor on~$\mathscr X$
extending~$D_i$ defined by Proposition~\ref{prop.metric}
and $g_i=(g_{D_i,\Omega_{i,v}})$
is the family of capacitary Green currents  at archimedean places.
The vertical components of~$\mathscr D_1$ and~$\mathscr D_2$
lying over any finite place~$v$ are distinct the one from the other,
since $\Omega_{1,v}\cap\Omega_{2,v}=\emptyset$.
Consequently, the geometric part of the Arakelov intersection
product is zero.
In view of~\cite[Lemma~5.1]{bost1999},
the contribution of any archimedean place~$v$ is zero too, since $\Omega_{1,v}$
and $\Omega_{2,v}$ are disjoint.
This concludes the proof.
\end{proof}

% (Using~\cite{bost1999}, Lemma~3.2, it is also possible to give an 
% argument similar to that of the ultrametric case.\footnote{Ah bon ?})

The following proposition makes more explicit
the numerical effectivity hypothesis 
in Proposition~\ref{prop.arith-hodge}.
\begin{proposition}
Let~$X$, $\Omega$, $D$, $\mathscr X$,
$\widehat D_\Omega$ be as in the beginning of this Subsection.
\begin{enumerate}\def\theenumi{\alph{enumi}}
\item If~$D$ is effective, then
The Arakelov divisor  $\widehat D_\Omega$ on $\mathscr X,$ attached to the effective divisor $D$ and to the adelic tube $\Omega$, is
effective.
\item Write $D=\sum_i n_i P_i,$ for  some closed points~$P_i$
of~$X$ and positive integers~$n_i$.
Then $\widehat D_\Omega$ is numerically effective
if and only if 
$h_{\widehat D_\Omega}(P_i)\geq 0$ for each~$i$.
\item If~$D$ is a rational point~$P$, then $\widehat D_\Omega$
is numerically effective (resp. $\widehat D_\Omega.\widehat D_\Omega > 0$) if and only if 
the Arakelov degree $\hdeg (T_P X,\norm{\cdot}^\CAP_\Omega)$
is non-negative (resp. positive).
\end{enumerate}
\end{proposition}
\begin{proof}
\emph a)
Let us assume that~$D$ is an effective divisor.
For each archimedean place~$v$ of~$K$, the capacitary
Green function $g_{D,\Omega_v}$ is therefore nonnegative
(\cite[3.1.4]{bost1999}).
Moreover, we have proved in~Proposition~\ref{prop.metric}
that the~$\Q$-divisor~$\mathscr D$ in~$Z^1_\Q(\mathscr X)$
is effective. These two facts together imply that $\widehat D_\Omega$
is an effective Arakelov divisor.

\emph b)
For any archimedean place~$v$,
the definition of the archimedean capacitary Green currents
involved in~$\widehat D_\Omega$ implies that
$\omega(\widehat D_{\Omega_v})$ is a positive
measure on~$X_v$, zero near~$\abs D$ \cite[Theorem~3.1, (iii)]{bost1999}.
By~\cite[Proposition~6.9]{bost1999},
in order to $\widehat D_\Omega$ being numerically effective, 
it is necessary and sufficient that 
$h_{\widehat D_\Omega}(E)\geq 0$
for any irreducible
component~$E$ of~$\mathscr D$. 
This holds by construction if~$E$ is a vertical component
of~$\mathscr X$: according to the conditions of Proposition~\ref{prop.metric},
on has $\mathscr D\cdot V=0$
for any vertical component~$V$ of the support of~$\mathscr D$;
for any other vertical component~$V$, one
has $\mathscr D\cdot V\geq 0$ because the divisor~$\widehat D_\Omega$
is effective.
Consequently, $\widehat D_\Omega$ is nef if and only
if $h_{\widehat D_\Omega}(P_i)\geq 0$ for all~$i$.

\emph c) This follows from b) and from the equality
(Corollary~\ref{coro.height.capa})
\[ h_{\widehat D_\Omega}(P)= \widehat D_\Omega.\widehat D_\Omega=
\hdeg(T_PX,\norm{\cdot}^\CAP_\Omega).\]
\end{proof}

\begin{theorem}\label{theo.harbater} Let $X$ be a geometrically connected smooth projective curve over $K,$ 
$P$ is a rational point in~$X(K)$, and~$\Omega:=(\Omega_v)$ an adelic tube adapted to $P.$

Let $\phi\in\widehat{\mathscr O_{X,P}}$ be any formal function around~$P$
satisfying the following assumptions:
\begin{enumerate}
\item for any $v\in F$, $\phi$ extends to an analytic meromorphic
function on~$\Omega_v$;
\item $\phi$ is algebraic over $\mathscr O_{X,P}$;
\item $\hdeg (T_PX, \norm{\cdot}^\CAP_\Omega)\geq 0$.
\end{enumerate}
If equality holds in the last inequality,
assume moreover that there is an archimedean place~$v$ of~$F$
such that $X_v\setminus\Omega_v$ is not contained in a
real semi-algebraic curve of~$X_v$.

Then $\phi$ is the formal germ at~$P$ of a rational function in~$K(X)$.
\end{theorem}
\begin{proof}
Let $X'$ be the normalization of~$X$ in the field extension of~$K(X)$
generated by~$\phi$. This is a geometrically connected smooth projective 
curve over $K$, which may be identified with the normalization of the Zariski 
closure $Z$ in $X \times \mathbf P^1_K$
of the graph of $\phi.$ It is endowed with
a finite morphism $f\colon X'\ra X$, namely the composite morphism $X'\rightarrow Z \stackrel{\rm pr_1}{\rightarrow} X$. Moreover
the formal function~$\phi$ may be identified with the composition of the formal section~$\sigma$
of~$f$ at~$P$ that lifts the formal section $(\Id_X, \phi)$ of  $Z \stackrel{\rm pr_1}{\rightarrow} X$   and of the rational function $\tilde{\phi}$ in the local ring $\mathscr O_{X',\sigma(P)}$ defined as 
 the composition $X'\rightarrow Z \stackrel{\rm pr_2}{\rightarrow} \mathbf P^1_K.$
 
To show that $\phi$ is the germ at~$P$ of a rational  function,
we want  to show that~$f$ is an isomorphism. 

For any place~$v$,
$\Omega_v$ is a smooth analytic curve in~$X_v$, and
$\sigma$ extends to an analytic section~$\sigma_v\colon\Omega_v\ra X'_v$ of~$f$. Indeed, according to~1), the formal morphism $(\Id_X, \phi)$ extends to an analytic section of
$Z \stackrel{\rm pr_1}{\rightarrow} X$ over $\Omega_v$, which in turn lifts to an analytic section of~$f$ by normality.

By Corollary~\ref{coro.height.capa},
the Arakelov $\Q$-divisor $\widehat P_{\Omega}$
attached to the point~$P$ and the adelic tube~$\Omega$ is
nef. When $\hdeg (T_PX,\norm{\cdot}_{\Omega}^\CAP)$ is positive,
Proposition~\ref{prop.arith-hodge} implies that~$f$ is an isomorphism,
hence $\phi$ is the
formal germ to a rational function on~$X$.
This still holds  when $\hdeg (T_PX,\norm{\cdot}_{\Omega}^\CAP)=0$,
thanks to the supplementary assumption at archimedean places
in that case.
\end{proof}

As an example, this theorem applies when~$X$ is the projective
line, $P$ is the origin and when, for each place~$v$ in~$F$,
$\Omega_v$ is the disk of center~$0$ and of radius~$R_v\in \sqrt{\abs{K_v^*}}$ in the affine
line. Then $(\Omega_v)$ is an adelic tube adapted to $P$ iff almost every $R_v$ equals 1, and
$\hdeg (T_PX,\norm{\cdot}_{\Omega}^\CAP)$ is non-negative iff  $\prod_v R_v\geq 1$. In this special case, Theorem~\ref{theo.harbater} becomes Harbater's 
rationality criterion (\cite[Proposition~2.1]{harbater1988}). 

Actually  Harbater's result is stated without the assumption~$R_v\in \sqrt{\abs{K_v^*}}$ on the non-archimedean radii. The reader will easily check that his rationality criterion may be derived in full generality from Theorem~\ref{theo.harbater}, by shrinking the disks $\Omega_v$ for $v$ non archimedean, and replacing them by larger  simply connected domains for $v$ archimedean.

When $ \hdeg (T_PX, \norm{\cdot}_\Omega^\CAP)= 0$,
some  hypothesis on  the sets
$X_v\setminus\Omega_v$ is really necessary
for a rationality criterion to hold.
As an example, let us consider
the Taylor series of the algebraic function~$\phi(x)=1/\sqrt{1-4x}-1$,
viewed as a formal
function aroung the origin of the projective line~$\mathbf P^1_\Q$.
As shown by 
the explicit expansion
\[ \frac 1{\sqrt{1-4x}}-1 = \sum_{n=1}^\infty (-4)^n \binom{-1/2}{n} x^n
%  = \sum_{n=1}^\infty 4^n\frac{1\cdot 3\cdots (2n-1)}{2^n n!} x^n
 = \sum_{n=1}^\infty \binom {2n}{n} x^n ,\]
the coefficients of this series are rational integers. 
Moreover, the complementary subset~$\Omega$ of the real interval
$[1/4,\infty]$ in~$\P^1(\C)$ is a simply connected open Riemann surface
on which the algebraic function has no ramification. Consequently,
there is a meromorphic function $\phi_\infty$ on~$\Omega$
such that $\phi_\infty(x)=(1-4x)^{-1/2}-1$ around~$0$.
One has $\CAP_0(\Omega)=1$, hence $\hdeg(T_0 \P^1,\norm{\cdot}_\Omega^\CAP))=0$.
However, $\phi$ is obviously not a rational function.

% Let~$D$ be a divisor in~$S$ whose support is contained in~$\Omega_v$
% for any place~$v\in F$; as above capacity theory for $v\in F$
% and the chosen model $\mathscr S$ outside~$v$ allow to define 
% a capacitary metric on~$\mathscr O_S(D)$; the corresponding
% $L_1^2$-Arakelov divisor in~$\hCH^1(\mathscr X)$ will be denoted
% $\widehat D_\Omega$. It is effective if~$D$ is effective,
% meaning that writing $\widehat D_\Omega=(\mathscr D,g)$,
% one has $\ddc  g +\delta_D\geq 0$.
% 
% \begin{lemma}
% Assume moreover that~$D$ is effective; write
% its irreducible components as $\sum_{i=1}^r n_i D_i$. 
% Then the Arakelov
% divisor $\widehat D_\Omega$ is numerically effective
% if and only if for any~$i$, one has $\widehat D_\Omega\cdot D_i\geq 0$.
% \end{lemma}

By combining the algebraicity criterion
of Theorem~\ref{theo.algebraic} and 
the previous corollary, we deduce the following result,
a generalization to curves of any genus of
Borel-Dwork's criterion.

\begin{theorem}\label{theo.borel-dwork}
Let $X$ be a geometrically connected smooth projective curve over $K,$ 
$P$ is a rational point in~$X(K)$, and~$\Omega:=(\Omega_v)$ an adelic tube adapted to $P.$

Let $\phi\in\widehat{\mathscr O_{X,P}}$ be any formal function around~$P$
satisfying the following assumptions:
\begin{enumerate}
\item for any $v\in F$, $\phi$ extends to an analytic meromorphic
function on~$\Omega_v$;
\item the formal graph of~$\phi$ in~$\widehat{X\times\A^1}_{(P,\phi(P))}$
is~$A$-analytic.
\end{enumerate}
If moreover $\hdeg (T_PX, \norm{\cdot}^\CAP_\Omega)>0$, then $\phi$
is the formal germ at~$P$ of a rational function on~$X$ \emph{ (in other words, $\phi$ belongs to $\mathscr O_{X,P}$).}
\end{theorem}
\begin{proof}
In view of Corollary~\ref{theo.harbater},
it suffices to prove prove that~$\phi$ is algebraic.
Let $V=X\times\P^1$ 
and let $\widehat C\subset \widehat V_{(P,\phi(P))}$ 
be the formal graph of~$\phi$.  We need to prove that~$\widehat C$
is algebraic.
Indeed, since at each place~$v$ of~$K$,
the canonical~$v$-adic semi-norms on~$T_P\widehat C$
is smaller than the capacitary one, 
$\hdeg(T_PX,\norm{\cdot}_{X,\widehat C}^\can)\geq
\hdeg(T_PX,\norm{\cdot}_\Omega^\CAP)>0$.
By Theorem~\ref{theo.algebraic}, 
$\widehat C$ is then algebraic, and $\phi$
is algebraic over~$K(X)$.
\end{proof}

Observe that, when condition~1) is satisfied  
in Theorem~\ref{theo.borel-dwork},
%Let us recall from Proposition~\ref{preEisenstein} that 
the~$A$-analyticity condition~2) is implied by the following one:  
\begin{enumerate}
\item[{\textup 2$'$)}] \emph{there exist a positive
integer~$N$ and a smooth 
model~$\mathscr X$ of~$X$ over $\Spec \mathscr O_K[1/N]$
such that $P$ extends to an integral point~$\mathscr P$ in $\mathscr X(\mathscr O_K[1/N])$, and $\phi$ extends to a regular formal function on
the formal completion~$\widehat {\mathscr X}_{\mathscr P}$.}
\end{enumerate}
This follows from Proposition~\ref{preEisenstein}, since  then the formal graph of~$\phi$
extends to a smooth formal curve in $\mathscr X \times \mathbf A^1$ over~$\Spec  \mathscr O_K[1/N]$.

\begin{example}
Theorem~\ref{theo.borel-dwork} may be applied when $X$ is $\mathbf P^1_K$,
$P$ is the origin $0$ in $\mathbf A^1(K) \hookrightarrow \mathbf P^1(K)$, and when, for each place~$v$,
$\Omega_v\subset F_v$ is a disk of center~$0$ and of positive radius~$R_v$ 
in the affine line, provided these radii are almost all equal to 1 and satisfy $\prod R_v>1$. In this case, the rationality of any $\phi$ in
$\widehat{\mathscr O_{X,P}} \simeq K[[X]]$ under the assumptions~1) and~2$'$) 
is precisely Borel-Dwork's rationality criterion
(\cite{eborel1894,dwork1960}). 

More generally, the expression of  capacitary norms in terms of  transfinite diameters and a straightforward approximation
argument\footnote{using the fact that  bounded subsets of~$\C_p$
are contained in affinoids (actually, lemniscates)
with arbitrarily close transfinite diameters.} 
allows one to recover the criterion
of P\'olya-Bertrandias (\cite{polya1928,amice75}) from our Theorem~\ref{theo.borel-dwork} with  
$X=\mathbf P^1_K$.
%at least when, at finite places, the
%formal power series extends to an meromorphic function on
%open subsets~$\Omega_v$ which are 
%complements of affinoids. 
%However, the strict inequality in the hypothesis of the theorem
%and the fact that allow easily 
%to recover the full criterion from Theorem~\ref{theo.borel-dwork}.
\end{example}

\appendix

\section*{Appendix} 

\setcounter{theorem}{0}
\def\thesubsection{\Alph{subsection}}
\def\thetheorem{\thesubsection.\arabic{theorem}}

\subsection{Metrics on line bundles}
\label{app.metrics}

Let~$K$ be a field which is complete
with respect to the topology defined by a discrete absolute value~$\abs\cdot$
on~$K$. Let~$R$ be its valuation ring and $\pi$ be an uniformizing
element of~$R$. We denote by $v=\log\abs\cdot/\log\abs{\pi}$
the corresponding normalized valuation on~$K$.

Let~$X$ be an algebraic variety over~$K$
and let~$L$ be a line bundle on~$X$.
In this appendix, we precise basic facts concerning
the definition of a metric  on the fibres of~$L$.

Let $\bar K$ be an algebraic closure of~$K$; endow it
with the unique absolute value which extends the given one on~$K$.
It might not be complete, however its completion, denoted~$C$,
is a complete field containing $\bar K$ as a  dense subset
on which the absolute value extends uniquely, endowing it
with the structure of a complete valued field.

A metric on the fibres of~$L$ is the data, 
for any $x\in X(C)$, of a norm $\norm{\cdot}$ on the one-dimensional
$C$-vector space $L(x)$. Namely, $\norm{\cdot}$
is a map $L(x)\ra\R_+$ satisfying the following
properties:
\begin{itemize}
\item $\norm{s_1+s_2}\leq \max(\norm{s_1},\norm{s_2})$ 
for all $s_1$, $s_2\in L(x)$ ;
\item $\norm{as}=\abs a\norm{s}$ for all $a\in C$
and $s\in L(x)$;
\item $\norm{s}=0$ implies $s=0$.
\end{itemize}
We also assume that these norms are stable under the
natural action of the Galois group~$\Gal(C/K)$, namely that for any $x\in X(C)$,
$s\in L(x)$ and $\sigma\in\Gal(C/K)$, 
$\norm{\sigma(s)}=\norm{s}$.

We say a metric is continuous if for any open subset $U\subset X$
(for the Zariski topology)
and any section $s\in\Gamma(U,L)$, the
function $x\mapsto \norm{s(x)}$ on~$U(C)$ is continuous.
This definition corresponds to that classical
notion of a Weil function attached to a Cartier divisor on~$X$
and will be sufficient for our purposes;
a better one would be to impose that this function
extends to a continuous function on the  analytic space
attached to~$U$ by Berkovich~\cite{berkovich1990};
see \eg \cite{gubler1998} for this point of view.

Assume that~$X$ is projective and let  $\mathscr X$ be any 
projective and flat~$R$-scheme with generic fibre~$X$,
together with a line bundle $\mathscr L$ on~$\mathscr X$
extending~$L$.
Let $x\in X(C)$; if $C^0$ denotes the valuation ring of~$C$,
there is a unique morphism $\eps_x\colon \Spec C^0\ra \mathscr X$, 
by which the generic point of~$\Spec C^0$ maps to~$x$.
Then, $\eps_x^*\mathscr L$ is a sub-$C^0$-module of~$L(x)$.
For any section $s\in L(x)$, there exists $a\in C^0$
such that  $as\in \eps_x^*\mathscr L$. Define, for any $s\in L(x)$,
\[ \norm{s} = \inf \{ \abs{a}^{-1} \ ,\ as\in\eps_x^*\mathscr L, \quad a\in C\setminus \{0\} \}. \]
This is a continuous metric on the fibres of~$L$,
which we call an algebraic metric.

\medskip
Algebraic metrics are in fact the only metrics 
that we use in this article,
where the language of metrics is just  a convenient way
of comparing various extensions of~$X$ and~$L$ over~$R$.
In that respect, we make the following two remarks:

1) Let~$Y$ be another projective algebraic variety over~$K$
and let $f\colon Y\ra X$ be a morphism. Let $(L,\norm{\cdot}_L)$
be a metrized line bundle on~$X$. Then, the line bundle $f^*L$
on~$Y$ admits a metric~$\norm{\cdot}_{f^*L}$, defined by the formula
$\norm{f^*s(y)}_{f^*L}=\norm{s(f(y))}_{L}$,
where $y\in Y(C)$ and~$s$ is a section of~$L$ in
a neighbourhood of~$f(y)$.
Assume that the metric of~$L$ is algebraic, defined
by a model~$(\mathscr X,\mathscr L)$.
Let $\mathscr Y$ be any projective flat model of~$Y$ over~$R$
such that~$f$ extends to a morphism $\phi\colon\mathscr Y\ra\mathscr X$.
Then, the metric~$\norm{\cdot}_{f^*L}$ is algebraic,
defined by the pair $(\mathscr Y,\phi^*\mathscr L)$.

2) Let $\mathscr X$ be a projective and flat model of~$X$ on~$R$
and let~$\mathscr L$ and~$\mathscr L'$ be two line bundles
on~$\mathscr X$ which induce the same (algebraic) metric on~$L$.
If $\mathscr X$ is \emph{normal} then the identity map
$\mathscr L_K=\mathscr L'_K$ on the generic fibre
extends uniquely to an isomorphism $\mathscr L\simeq\mathscr L'$.

% By definition of the metric defined by~$\mathscr D$, one has,
% for any point~$x\in X(K)$,
% \[ - \log \norm{1_{\mathsf D}} (x) = (\mathscr D,\nu_*(\bar x))
%      = (\nu^*(\mathscr D), \bar x), \]
% where $\bar x$ is the closure of~$x$ in~$\mathscr X_0$.

% In~\cite{zhang95b}, Zhang considers more general metrics
% on the fibres of~$L$: they are limits of a sequence
% a algebraic metrics so that for any open subset $U\subset X$
% and any section $s\in\Gamma(U,L)$, the sequence of continuous functions
% $(x\mapsto \norm{s(x)})$ converges uniformly to a continuous
% function on~$U$.

\subsection{Background on rigid analytic geometry}
\label{app.rigid}

The results of this appendix are basic facts of rigid analytic
geometry: the first one is a version of
the maximum principle, while the second proposition
states that the complementary subsets to an affinoid subspace
in a rigid analytic space has a canonical structure of a rigid space.
They are well known to specialists but,
having been unable to find a convenient reference,
we decided to write them here.

Let~$K$ be a field, endowed with a ultrametric absolute value
for which it is complete.

\begin{proposition}\label{prop.max-principle}
Let~$C$ be a smooth projective connected curve over~$K$,
let $f\in K(C)$ be a non constant rational function
and let~$X$ denote the Weierstrass domain 
$C(f)=\{x\in X\sozat \abs{f(x)}\leq 1\}$ in~$X$. 
Then, any affinoid function~$g$ on~$X$ is bounded; moreover,
there exists $x\in U$ such that
\[ \abs{g(x)} = \sup_{X} \abs g \quad
\text{and}\quad \abs{f(x)}=1. \]
\end{proposition}

The fact that~$g$ is bounded and attains it maximum is the
classical maximum principle; we just want to assure that the
maximum is attained on the ``boundary'' of~$U$.
% \begin{proof}[First proof (Berkovich style)]
% For this proof, we view all rigid spaces in the sense of Berkovich.
% What we have to prove is that the Shilov boundary of~$U$
% is contained in the affinoid defined by $\abs{f(x)}=1$.
% For one-dimensional affinoids, the Shilov boundary is just
% the complement of the set of interior points
% (the Shilov boundary consists of points reducing to the generic points of the
% canonical reduction, the interior points reduce to closed points.)
% The rational function~$f$ defines a finite morphism
% $f\colon X\ra \P^1$ and restricts to a finite
% analytic morphism $f_{U} \colon U \ra \mathbf B$
% where~$\mathbf B$ denotes the unit ball.
% By~\cite{berkovich1990}, 2.5.8, p.~39, the set of interior
% points of~$U$ is given by 
% \[ \Int(U) = f_{U}^{-1}(\Int(\mathbf B)) \cap \Int(U/\mathbf B). \]
% Moreover, since~$f_{U}$ is finite,
% the set $\Int(U/\mathbf B)$
% of interior points of the morphism $U\ra\mathbf B$
% is equal to~$U$.
% Consequently, the Shilov boundary of~$U$ is the preimage
% by~$f_{U}$ of the Shilov boundary of~$\mathbf B$,
% which consists of a single point, corresponding to the Gauss norm
% on the Tate algebra. This is also the Shilov boundary
% of the collar $\abs{x}=1$ in~$\mathbf B$, whence the proposition.
% \end{proof}
 
\begin{proof} % [Second proof (classical rigid style)]
The analytic map $f\colon C^\an\ra(\P^1)^\an$
induced by~$f$ is finite hence restricts to a finite
map $f_X\colon X\ra\mathbf B$ of rigid analytic spaces,
where $\mathbf B=\operatorname{Sp}K\langle t\rangle$ is the unit ball.
It corresponds to~$f_X$ a morphism of affinoid algebras
$K\langle t\rangle \hra \mathscr O(X)$ 
which makes $\mathscr O(X)$ a $K\langle t\rangle$-module of finite type. 
Let $g\in\mathscr O(X)$ be an analytic function.
Then~$g$ is integral over $K\langle t\rangle$, hence there is
a smallest positive integer~$n$, as well as
analytic functions $a_i\in K\langle t\rangle$, for $1\leq i\leq n$,
such that
\[ g(x)^n + a_1(f(x)) g(x)^{n-1} + \dots + a_n(f(x)) = 0 \]
for any~$x\in X$.
Then, (see~\cite[p.~239, Proposition~6.2.2/4]{bosch-g-r1984})
\[ \sup_{x\in X} \abs {g(x)} 
= \max_{1\leq i\leq n} \abs{a_i(t)}^{1/i}. \]
The usual proof of the maximum principle on~$\mathbf B$ shows that
there is for each integer~$i\in\{1,\dots,n\}$
a point $t_i\in \mathbf B$ satisfying $\abs {t_i}=1$
and $\abs{a_i(t_i)}= \norm{a_i}$. (After having reduced
to the case where~$\norm{a_i}=1$, it suffices to lift
any non-zero element of the residue field at which the
reduced polynomial~$\overline{t_i}$ does not vanish.)
Consequently, there is therefore a point $t\in \mathbf B$ such that
$\abs t=1$ and
\[ \max_i \abs{a_i(t)}^{1/i} = \max_i \norm{a_i}^{1/i}. \]
Applying Proposition~3.2.1/2, p.~129, of~\cite{bosch-g-r1984}
to the polynomial $Y^n+a_1(t)Y^{n-1}+\cdots+a_n(t)$,
there is a point $y\in\P^1$ and $\abs y=\max_i\norm{a_i}^{1/i}$. 
Since the morphism $K\langle t\rangle[g]\subset\mathscr O(X)$
is integral, there is a point $x\in X$ such that $f(x)=t$ 
and $g(x)=y$. For such a point, one has $\abs{f(x)}=1$
and $\abs g(x)=\norm g$.
\end{proof}

% \begin{proof}[Third proof (formal style)]
% There is a formal model $\mathfrak X$ of~$X$ over~$R$
% such that~$U$ is the generic fibre of a formal open
% subset $\mathfrak U$ whose special fibre is just the complement
% of a finite number of points. We may moreover
% contract complete curves contained in the special fibre of~$\mathfrak U$.
% Then, it follows that there are horizontal effective
% divisors~$Z$ and~$P$ on~$\mathfrak X$ such that
% $ \div(f) = Z  - P$,
% the special fibre~$\mathsf U$ of~$\mathfrak U$ consisting of the complementary
% subset of~$P$. 
% 
% Let us look at the divisor of~$g$ on~$\mathfrak U$. We may assume
% that $\norm{g}=1$, hence $\div(g)=H+\sum_j b_i V_i$,
% for some nonnegative integers~$b_i$ and some vertical
% components; moreover some~$b_i=0$.
% Any point reducing~$x$ to such a component~$V_i$ but
% not to a point of~$P$ satisfies $\abs{f(x)}=\abs{g(x)}=1$.
% \end{proof}
 
\begin{proposition}
Let~$X$ be a rigid analytic variety over~$K$
and let $A\subset X$ be the union of finitely
many affinoid subsets.

Then $X\setminus A$, endowed with the induced G-topology,
is a rigid analytic variety.
\end{proposition}
\begin{proof}
By~\cite[p.~357, Proposition~9.3.1/5]{bosch-g-r1984}, and the
remark which follows that proposition,
it suffices to prove that $X\setminus A$
is an admissible open subset.

Let $(X_i)$ be an admissible affinoid covering of~$X$;
then, for each~$i$,  $A_i=A\cap X_i$ is a finite union
of affinoid subsets of~$X_i$. Assume that the Proposition
holds when~$X$ is affinoid; then, each $X_i\setminus A_i$
is an admissible open subset of~$X_i$, hence of~$X$.
Then $X\setminus A=\bigcup_i(X_i\setminus A_i)$ is
an admissible open subset of~$X$, by the property~(G\textsubscript 1)
satisfied by the G-topology of rigid analytic  varieties.

We thus may assume that~$X$ is an affinoid variety.
By Gerritzen-Grauert's theorem  
(\cite[p.~309, Corollary~7.3.5/3]{bosch-g-r1984}),
$A$ is a finite union of rational subdomains~$(A_i)_{1\leq i\leq m}$ in~$X$.
For each~$i$, let us consider affinoid functions 
$(f_{i,1},\dots,f_{i,n_i},g_i)$ on~$X$ generating the unit ideal such 
that 
\begin{align*} A_i & =X\bigg(\frac{f_{i,1}}{g_i},\dots,\frac{f_{i,n_i}}{g_i}\bigg)\\
& =\left\{x\in X\sozat \abs{f_{i,1}(x)}\leq\abs{g_i(x)},
\dots,\abs{f_{i,n_i}(x)}\leq \abs{g_i(x)} \vphantom{\bigg)}\right\}. \end{align*} 
We have
\[ X\setminus A = \bigcap_{i=1}^m (X\setminus A_i)
=\bigcap_{i=1}^m \bigcup_{j=1}^{n_i} \left\{ x\in X\sozat
   \abs{f_{i,j}(x)}>\abs{g_i(x)} \vphantom{\bigg)}\right\} . \]
Since any finite intersection of admissible open subsets
is itself admissible open, it suffices to treat the
case where $m=1$, \ie, when~$A$ is a rational subdomain
$X(f_1,\dots,f_n;g)$ of~$X$, which we now assume.

By assumption, $f_1,\dots,f_n,g$ have no common zero. 
By the maximum principle \cite[p.~307, Lemma~7.3.4/7]{bosch-g-r1984},
there is $\delta\in\sqrt{\abs{K^*}}$ such that for any~$x\in X$,
\[ \max(\abs{f_1(x)},\dots,\abs{f_n(x)},\abs{g(x)})\geq\delta .\]
For any $\alpha\in \sqrt{\abs{K^*}}$ with $\alpha>1$,
and any $j\in\{1,\dots,n\}$,
define
\[  X_{j,\alpha}=X\left(\delta\frac1{f_j},\alpha^{-1}\frac g{f_j} \right)
 = \left\{\vphantom{\bigg)}x\in X\sozat \delta\leq\abs{f_j(x)},
            \quad \alpha\abs{g(x)}\leq\abs{f_j(x)}\right\}.\]
This is a rational domain in~$X$. 
For any~$x\in X_{j,\alpha}$, one has $f_j(x)\neq 0$, and 
$ \abs{g(x)}<\abs{f_j(x)}$, hence $x\in X\setminus A$. Conversely,
if $x\in X\setminus A$, there exists $j\in\{1,\dots,n\}$
such that $\max(\abs{f_1(x)},\dots,\abs{f_n(x)},\abs{g(x)})
=\abs{f_j(x)}>\abs{g(x)}$;
it follows that there is $\alpha\in\sqrt{\abs{K^*}}$, $\alpha>1$,
such that $x\in X_{j,\alpha}$. This shows that the affinoid
domains $X_{j,\alpha}$ of~$X$,
for $1\leq j\leq n$ and $\alpha\in\sqrt{\abs{K^*}}$,
$\alpha>1$, form a covering of~$X\setminus A$.
Let us show that this covering is admissible.
Let~$Y$ be an affinoid space and let $\phi\colon Y\ra X$
be an affinoid map such that $\phi(Y)\subset X\setminus A$.
By~\cite[p.~342, Proposition~9.1.4/2]{bosch-g-r1984},
we need to show that  the covering $(\phi^{-1}(X_{j,\alpha)})_{j,\alpha}$
of~$Y$ has a (finite) affinoid covering which refines it.
For that, it is sufficient to prove 
that there are real numbers~$\alpha_1,\dots,\alpha_n$
in~$\sqrt{\abs{K^*}}$, greater than~$1$,
such that $\phi(Y)\subset \bigcup_{j=1}^n X_{j,\alpha_j}$.

For $j\in\{1,\dots,n\}$, define an affinoid subspace $Y_j$ of~$Y$
by
\[ Y_j=\bigg\{y\in Y\sozat \abs{f_i(\phi(y))}\leq\abs{f_j(\phi(y))}\ \text{for $1\leq i\leq n$}\bigg\}. \]
One has $Y=\bigcup_{j=1}^n Y_j$.
Fix some $j\in\{1,\ldots,n\}$.
Since $\phi(Y_j)\subset X\setminus A$,
$\abs {g(x)}<\abs{f_j(x)}$ on~$Y_j$. It follows that 
$f_j\circ\phi$ does not vanish on~$Y_j$, hence
$g\circ\phi/f_j\circ\phi$ is an affinoid function on~$Y_j$ such that
\[ \left|\frac{g\circ \phi}{f_j\circ\phi}(y)\right|<1 \]
for any $y\in Y_j$.
By the maximum principle, there is $\alpha_j\in\sqrt{\abs{K^*}}$
such that $\alpha_j>1$ and
$\left|{\frac{g\circ \phi}{f_j\circ\phi}}\right|<\frac1{\alpha_j}$ on~$Y_j$.
One then has $\phi(Y)\subset \bigcup_{j=1}^nX_{j,\alpha_j}$,
which concludes the proof of the proposition.
\end{proof}

\nocite{cantor1980,polya1928}
\providecommand{\smfedsname}{eds.}
\providecommand{\smfedname}{ed.}

\bibliographystyle{smfplain}
\bibliography{BCL070720}

\providecommand{\noopsort}[1]{}\providecommand{\url}[1]{\textit{#1}}
\providecommand{\bysame}{\leavevmode ---\ }
\providecommand{\og}{``}
\providecommand{\fg}{''}
\providecommand{\smfandname}{\&}
\providecommand{\smfedsname}{\'eds.}
\providecommand{\smfedname}{\'ed.}
\providecommand{\smfmastersthesisname}{M\'emoire}
\providecommand{\smfphdthesisname}{Th\`ese}
\begin{thebibliography}{10}

\bibitem{amice75}
{\scshape Y.~Amice} -- \emph{Les nombres $p$-adiques}, Collection SUP: Le
  Math\'ematicien, vol.~14, Presses Universitaires de France, Paris, 1975.

\bibitem{andre1989}
{\scshape Y.~Andr\'e} -- \emph{${G}$-functions and geometry}, Vieweg,
  Braunschweig, 1989.

\bibitem{badescu2004}
{\scshape L.~B{\u{a}}descu} -- \emph{Projective geometry and formal geometry},
  Mathematics Institute of the Polish Academy of Sciences. Mathematical
  Monographs (New Series), vol.~65, Birkh\"auser Verlag, Basel, 2004.

\bibitem{berkovich1990}
{\scshape V.~G. Berkovich} -- \emph{Spectral theory and analytic geometry over
  non-{A}rchimedean fields}, Mathematical Surveys and Monographs, vol.~33,
  American Mathematical Society, Providence, RI, 1990.

\bibitem{bogomolov-mq2001}
{\scshape F.~A. Bogomolov {\normalfont \smfandname} M.~L. McQuillan} -- {\og
  Rational curves on foliated varieties\fg}, pr\'epublication M/01/07,
  I.H.{\'E}.S., 2001.

\bibitem{eborel1894}
{\scshape {\'E}.~Borel} -- {\og Sur une application d'un th\'eor\`eme de
  {M}.~{H}adamard\fg}, \emph{Bulletin des sciences math\'ematiques} \textbf{18}
  (1894), p.~22--25.

\bibitem{bosch-g-r1984}
{\scshape S.~Bosch, U.~G{\"u}ntzer {\normalfont \smfandname} R.~Remmert} --
  \emph{Non-{A}rchimedean analysis}, Grundlehren der Mathematischen
  Wissenschaften, vol. 261, Springer-Verlag, Berlin, 1984.

\bibitem{bosch-l1993}
{\scshape S.~Bosch {\normalfont \smfandname} W.~L{\"u}tkebohmert} -- {\og
  Formal and rigid geometry. {I}. {R}igid spaces\fg}, \emph{Math. Ann.}
  \textbf{295} (1993), no.~2, p.~291--317.

\bibitem{bosch-l1993b}
\bysame , {\og Formal and rigid geometry. {II}. {F}lattening techniques\fg},
  \emph{Math. Ann.} \textbf{296} (1993), no.~3, p.~403--429.

\bibitem{bosch-l-r90}
{\scshape S.~Bosch, W.~L{\"u}tkebohmert {\normalfont \smfandname} M.~Raynaud}
  -- \emph{N\'eron models}, Ergebnisse der Mathematik und ihrer Grenzgebiete,
  vol.~21, Springer-Verlag, 1990.

\bibitem{bost1999}
{\scshape J.-B. Bost} -- {\og Potential theory and {L}efschetz theorems for
  arithmetic surfaces\fg}, \emph{Ann. Sci. {\'E}cole Norm. Sup.} \textbf{32}
  (1999), no.~2, p.~241--312.

\bibitem{bost2001}
\bysame , {\og Algebraic leaves of algebraic foliations over number fields\fg},
  \emph{Publ. Math. Inst. Hautes {\'E}tudes Sci.} \textbf{93} (2001),
  p.~161--221.

\bibitem{bost2004}
\bysame , {\og Germs of analytic varieties in algebraic varieties: canonical
  metrics and arithmetic algebraization theorems\fg}, (A.~Adolphson,
  F.~Baldassarri, P.~Berthelot, N.~Katz {\normalfont \smfandname} F.~Loeser,
  \smfedsname), vol.~I, Walter de Gruyter GmbH \& Co. KG, Berlin, 2004,
  p.~371--418.

\bibitem{bost2007}
\bysame , {\og Evaluation maps, slopes, and algebraicity criteria\fg}, in
  \emph{Proceedings of the International Congress of Mathematicians (Madrid
  2006)} (M.~Sanz-Sol{\'e}, J.~Soria, J.~L. Varona {\normalfont \smfandname}
  J.~Verdera, \smfedsname), vol.~II, European Mathematical Society, 2007,
  p.~371--418.

\bibitem{bost-g-s94}
{\scshape J.-B. Bost, H.~Gillet {\normalfont \smfandname} C.~Soul\'e} -- {\og
  Heights of projective varieties and positive {G}reen forms\fg}, \emph{J.
  Amer. Math. Soc.} \textbf{7} (1994), p.~903--1027.

\bibitem{cantor1980}
{\scshape D.~G. Cantor} -- {\og On an extension of the definition of
  transfinite diameter and some applications\fg}, \emph{J. Reine Angew. Math.}
  \textbf{316} (1980), p.~160--207.

\bibitem{chambert-loir2002}
{\scshape A.~Chambert-Loir} -- {\og Th\'eor\`emes d'alg\'ebricit\'e en
  g\'eom\'etrie diophantienne (d'apr\`es {J}.-{B}.\ {B}ost, {Y}.\ {A}ndr\'e,
  {D}. \& {G}.\ {C}hudnovsky)\fg}, \emph{Ast{\'e}risque} (2002), no.~282,
  p.~175--209, Exp. No. 886.

\bibitem{chen2006}
{\scshape H.~Chen} -- {\og Positivit\'e en g\'eom\'etrie alg\'ebrique et en
  g\'eom\'etrie d'{A}rakelov : application \`a l'alg\'ebrisation et \`a
  l'\'etude asymptotique des polygones de {H}arder-{N}arasimhan\fg}, Th\`ese,
  \'Ecole polytechnique, 2006.

\bibitem{chudnovsky1985b}
{\scshape D.~V. Chudnovsky {\normalfont \smfandname} G.~V. Chudnovsky} -- {\og
  Applications of {P}ad\'e approximations to the {G}rothendieck conjecture on
  linear differential equations\fg}, in \emph{Number theory} (New York,
  1983--84), Lecture Notes in Math., vol. 1135, 1985, p.~52--100.

\bibitem{chudnovsky1985a}
{\scshape D.~V. Chudnovsky {\normalfont \smfandname} G.~V. Chudnovsky} -- {\og
  Pad\'e approximations and {D}iophantine geometry\fg}, \emph{Proc. Nat. Acad.
  Sci. U.S.A.} \textbf{82} (1985), no.~8, p.~2212--2216.

\bibitem{deligne1972b}
{\scshape P.~Deligne} -- {\og Intersections sur les surfaces
  r\'eguli\`eres\fg}, in \emph{Groupes de monodromie en g\'eom\'etrie
  alg\'ebrique (SGA~7~II)}, Lecture Notes in Math., vol. 340, Springer-Verlag,
  1973, p.~1--38.

\bibitem{dwork1960}
{\scshape B.~Dwork} -- {\og On the rationality of the zeta function of an
  algebraic variety\fg}, \emph{Amer. J. Math.} \textbf{82} (1960), p.~631--648.

\bibitem{eisenstein1852}
{\scshape G.~Eisenstein} -- {\og {\" U}ber eine allgemeine {E}igenschaft der
  {R}eihen-{E}ntwicklungen aller algebraischen {F}unktionen (1852)\fg}, in
  \emph{Mathematische Gesammelte Werke, Band II}, Chelsea Publishing Co., New
  York, 1975, p.~765--767.

\bibitem{franchetta1949}
{\scshape A.~Franchetta} -- {\og Sulle curve riducibili appartenenti ad una
  superficie algebrica\fg}, \emph{Univ. Roma. Ist. Naz. Alta. Mat. Rend. Mat. e
  Appl. (5)} \textbf{8} (1949), p.~378--398.

\bibitem{fresnel-m1986}
{\scshape J.~Fresnel {\normalfont \smfandname} M.~Matignon} -- {\og Sur les
  espaces analytiques quasi-compacts de dimension~1 sur un corps valu\'e
  complet ultram\'etrique\fg}, \emph{Ann. Mat. Pura Appl. (4)} \textbf{145}
  (1986), p.~159--210.

\bibitem{graftieaux2001a}
{\scshape P.~Graftieaux} -- {\og Formal groups and the isogeny theorem\fg},
  \emph{Duke Math. J.} \textbf{106} (2001), no.~1, p.~81--121.

\bibitem{graftieaux2001b}
\bysame , {\og Formal subgroups of abelian varieties\fg}, \emph{Invent. Math.}
  \textbf{145} (2001), no.~1, p.~1--17.

\bibitem{gruson-vdp1974}
{\scshape L.~Gruson {\normalfont \smfandname} M.~van~der Put} -- {\og Banach
  spaces\fg}, \emph{M{\'e}m. Soc. Math. France} (1974), no.~39-40, p.~55--100.

\bibitem{gubler1998}
{\scshape W.~Gubler} -- {\og Local heights of subvarieties over non-archimedean
  fields\fg}, \emph{J. Reine Angew. Math.} \textbf{498} (1998), p.~61--113.

\bibitem{harbater1988}
{\scshape D.~Harbater} -- {\og Galois covers of an arithmetic surface\fg},
  \emph{Amer. J. Math.} \textbf{110} (1988), no.~5, p.~849--885.

\bibitem{hartshorne1968}
{\scshape R.~Hartshorne} -- {\og Cohomological dimension of algebraic
  varieties\fg}, \emph{Ann. of Math. (2)} \textbf{88} (1968), p.~403--450.

\bibitem{hartshorne1969}
\bysame , {\og Curves with high self-intersection on algebraic surfaces\fg},
  \emph{Publ. Math. Inst. Hautes {\'E}tudes Sci.} (1969), no.~36, p.~111--125.

\bibitem{hartshorne1970}
\bysame , \emph{Ample subvarieties of algebraic varieties}, Lecture Notes in
  Mathematics, vol. 156, Springer-Verlag, Berlin, 1970.

\bibitem{hironaka1968}
{\scshape H.~Hironaka} -- {\og On some formal imbeddings\fg}, \emph{Illinois J.
  Math.} \textbf{12} (1968), p.~587--602.

\bibitem{hironaka-m1968}
{\scshape H.~Hironaka {\normalfont \smfandname} H.~Matsumura} -- {\og Formal
  functions and formal embeddings\fg}, \emph{J. Math. Soc. Japan} \textbf{20}
  (1968), p.~52--82.

\bibitem{ihara1994}
{\scshape Y.~Ihara} -- {\og Horizontal divisors on arithmetic surfaces
  associated with {B}ely\u\i\ uniformizations\fg}, in \emph{The Grothendieck
  theory of dessins d'enfants (Luminy 1993)} (L.~Schneps, \smfedname), London
  Math. Soc. Lecture Note Ser., vol. 200, Cambridge Univ. Press, Cambridge,
  1994, p.~245--254.

\bibitem{kani1989}
{\scshape E.~Kani} -- {\og Potential theory on curves\fg}, in \emph{Th\'eorie
  des nombres (Quebec, PQ, 1987)}, de Gruyter, Berlin, 1989, p.~475--543.

\bibitem{manin1985}
{\scshape {\relax Yu}.~I. Manin} -- {\og New dimensions in geometry\fg}, in
  \emph{Workshop Bonn 1984 (Bonn, 1984)}, Lecture Notes in Math., vol. 1111,
  Springer, Berlin, 1985, p.~59--101.

\bibitem{matsumura1980}
{\scshape H.~Matsumura} -- \emph{Commutative algebra}, Mathematics Lecture
  Notes Series, Benjamin/Cummings, 1980.

\bibitem{moret-bailly1989}
{\scshape L.~Moret-Bailly} -- {\og Groupes de {P}icard et probl\`emes de
  {S}kolem. {I}.\fg}, \emph{Ann. Sci. {\'E}cole Norm. Sup.} \textbf{22} (1989),
  no.~2, p.~161--179.

\bibitem{moret-bailly1996}
\bysame , {\og Un probl\`eme de descente\fg}, \emph{Bull. Soc. Math. France}
  \textbf{124} (1996), p.~559--585.

\bibitem{mumford1961}
{\scshape D.~Mumford} -- {\og The topology of normal singularities of an
  algebraic surface and a criterion for simplicity\fg}, \emph{Publ. Math. Inst.
  Hautes {\'E}tudes Sci.} \textbf{9} (1961), p.~5--22.

\bibitem{polya1928}
{\scshape G.~P\'olya} -- {\og {\"U}ber gewisse notwendige
  {D}eterminantenkriterien f\"ur die {F}ortsetzbarkeit einer {P}otenzreihe\fg},
  \emph{Math. Ann.} \textbf{99} (1928), p.~687--706.

\bibitem{vanderput1980}
{\scshape M.~van~der Put} -- {\og The class group of a one-dimensional affinoid
  space\fg}, \emph{Ann. Inst. Fourier (Grenoble)} \textbf{30} (1980), no.~4,
  p.~155--164.

\bibitem{ramanujam1972}
{\scshape C.~P. Ramanujam} -- {\og Remarks on the {K}odaira vanishing
  theorem\fg}, \emph{J. Indian Math. Soc. (N.S.)} \textbf{36} (1972),
  p.~41--51.

\bibitem{randriam2006}
{\scshape H.~Randriambololona} -- {\og M\'etriques de sous-quotient et
  th\'eor\`eme de {H}ilbert-{S}amuel arithm\'etique pour les faisceaux
  coh\'erents\fg}, \emph{J. Reine Angew. Math.} \textbf{590} (2006), p.~67--88.

\bibitem{raynaud94}
{\scshape M.~Raynaud} -- {\og Rev{\^e}tements de la droite affine en
  caract\'eristique $p>0$ et conjecture d'{A}bhyankar\fg}, \emph{Invent. Math.}
  \textbf{116} (1994), p.~425--462.

\bibitem{rumely-l-v2000}
{\scshape R.~Rumely, C.~F. Lau {\normalfont \smfandname} R.~Varley} -- {\og
  Existence of the sectional capacity\fg}, \emph{Mem. Amer. Math. Soc.}
  \textbf{145} (2000), no.~690, p.~1--130.

\bibitem{rumely1989}
{\scshape R.~S. Rumely} -- \emph{Capacity theory on algebraic curves}, Lecture
  Notes in Math., vol. 1378, Springer-Verlag, Berlin, 1989.

\bibitem{serre1997}
{\scshape J.-P. Serre} -- \emph{Lectures on the {M}ordell-{W}eil theorem},
  third \smfedname, Aspects of Mathematics, Friedr. Vieweg \& Sohn,
  Braunschweig, 1997.

\bibitem{thuillier2005}
{\scshape A.~Thuillier} -- {\og Th\'eorie du potentiel sur les courbes en
  g\'eom\'etrie non archim\'edienne. {A}pplications \`a la th\'eorie
  d'{A}rakelov\fg}, Th\`ese, Universit\'e de Rennes 1, 2005.

\end{thebibliography}

\end{document}